\documentclass[reqno,10pt,centertags,draft]{amsart}
\usepackage{amsmath,amsthm,amscd,amssymb,latexsym,upref}

\newcommand{\bbN}{{\mathbb{N}}}
\newcommand{\bbR}{{\mathbb{R}}}

\newcommand{\bbC}{{\mathbb{C}}}

\newcommand{\calD}{{\mathcal D}}

\newcommand{\calR}{{\mathcal R}}
\newcommand{\calS}{{\mathcal S}}

\newcommand{\dott}{\,\cdot\,}
\newcommand{\no}{\nonumber}
\newcommand{\lb}{\label}
\newcommand{\f}{\frac}

\newcommand{\ol}{\overline}

\newcommand{\wti}{\widetilde}
\newcommand{\what}{\widehat}

\newcommand{\oh}{o}

\newcommand{\loc}{\text{\rm{loc}}}
\newcommand{\spec}{\text{\rm{spec}}}
\newcommand{\rank}{\text{\rm{rank}}}

\newcommand{\dom}{\text{\rm{dom}}}
\newcommand{\ess}{\text{\rm{ess}}}

\newcommand{\supp}{\text{\rm{supp}}}
\newcommand{\esssup}{\text{\rm{ess\,sup}}}
\newcommand{\AC}{\text{\rm{AC}}}
\newcommand{\bi}{\bibitem}

\newcommand{\beq}{\begin{equation}}
\newcommand{\eeq}{\end{equation}}
\newcommand{\ba}{\begin{align}}
\newcommand{\ea}{\end{align}}
\newcommand{\abs}[1]{\lvert#1\rvert}

\renewcommand{\Re}{\text{\rm Re}}
\renewcommand{\Im}{\text{\rm Im}}
\renewcommand{\ln}{\text{\rm ln}}

\allowdisplaybreaks
\numberwithin{equation}{section}

\newtheorem{theorem}{Theorem}[section]
\newtheorem{lemma}[theorem]{Lemma}
\newtheorem{corollary}[theorem]{Corollary}
\newtheorem{hypothesis}[theorem]{Hypothesis}
\theoremstyle{definition}
\newtheorem{definition}[theorem]{Definition}
\newtheorem{example}[theorem]{Example}

\theoremstyle{remark}
\newtheorem{remark}[theorem]{Remark}

\begin{document}

\title[M-Function Asymptotics and Borg-type Theorems]{
Weyl-Titchmarsh $M$-Function Asymptotics, Local Uniqueness 
Results, Trace Formulas, and Borg-type
Theorems for Dirac Operators}
\author[Clark and Gesztesy]{Steve Clark and Fritz Gesztesy}
\address{Department of Mathematics and Statistics, University
of Missouri-Rolla, Rolla, MO 65409, USA}
\email{sclark@umr.edu}
\urladdr{http://www.umr.edu/\~{ }clark}
\address{Department of Mathematics,
University of
Missouri, Columbia, MO
65211, USA}
\email{fritz@math.missouri.edu\newline
\indent{\it URL:}
http://www.math.missouri.edu/people/fgesztesy.html}
\dedicatory{Dedicated to F.~V.~Atkinson, one of the pioneers of this subject.}
\thanks{Supported in part by NSF grant INT-9810322.}
\subjclass{Primary 34B20, 34E05, 34L40;  Secondary 34A55.}
\keywords{Weyl-Titchmarsh matrices, high-energy expansions, uniqueness
results, trace formulas, Borg theorems,  Dirac operators.}

\begin{abstract}
We explicitly determine the high-energy asymptotics for
Weyl-Titchmarsh matrices associated with general
 Dirac-type operators on half-lines and on $\bbR$. We 
also prove new local uniqueness results for Dirac-type operators in terms
of exponentially  small differences of Weyl-Titchmarsh matrices. As
concrete  applications of the asymptotic high-energy expansion we derive
a  trace formula for Dirac operators and use it to prove a Borg-type 
theorem.
\end{abstract}

\maketitle

\section{Introduction}\lb{s1}

While the high-energy asymptotics, $|z|\to\infty$, of scalar-valued
Weyl-Titchmarsh functions, $m_+(z,x_0)$, associated with general
half-line Dirac-type differential expressions of the form
\begin{equation}
J\f{d}{dx}-B(x), \quad J=\begin{pmatrix} 0 &-1\\1 &0
\end{pmatrix}
\end{equation}
and $B$ a self-adjoint $2\times 2$ matrix with 
real-valued coefficients, $B^{(n)}\in L^1([x_0,c])^{2\times 2}$ for some
$n\in\bbN_0(=\bbN\cup\{0\})$ and all $c>x_0$, received some attention
over the past two decades as can be inferred, for instance, from
\cite{EHS83}, \cite{Ha85}, \cite{HKS89a}, \cite{HKS89b}, \cite{Mi91} 
(and the literature therein), it may perhaps come as a surprise
that the corresponding matrix extension of this problem,
considering general matrix-valued differential expressions of the
type 
\begin{equation}
J\f{d}{dx}-B(x), \quad J=\begin{pmatrix} 0 &-I_m\\I_m &0\end{pmatrix}
\end{equation}
with $I_m$ the identity matrix in $\bbC^m$, $m\in\bbN$, and $B$ a
self-adjoint
$2m\times 2m$ matrix satisfying $B^{(n)}\in L^1([x_0,c])^{2m\times 2m}$
for some $n\in\bbN_0$ and all $c>x_0$, apparently, received no attention
at all. (It should be noted that this observation discounts papers in the
special scattering theoretic case concerned with short-range coefficients
$B^{(n)}\in L^1([x_0,\infty); (1+|x|)dx)^{2m\times 2m}$, where 
iterations of Volterra-type integral equations yield the asymptotic
high-energy expansion of $M_+(z,x_0)$ as $|z|\to\infty$ to any order, 
cf.~Lemma~\ref{l4.1}.) This is not because of a lack of interest in this
type of problem (we will discuss its relevance below), but simply since 
it is a nontrivial one, which, in many of its aspects, must be regarded as
more difficult than the corresponding matrix-valued
Schr\"odinger operator case, which in turn, was only very recently
settled in \cite{CG99}. The results proven in this paper show that in
leading order (and independently of the self-adjoint boundary condition
chosen at $x_0$),
\begin{equation}\lb{1.1}
M_+(z,x_0) \underset{\substack{\abs{z}\to\infty\\ z\in
C_\varepsilon}}{=} iI_{m} +\oh(1), 
\end{equation}
where $C_{\varepsilon}$ denotes the open sector in the open upper
complex half-plane $\bbC_+$ with vertex at zero, symmetry axis along the
positive imaginary axis, and opening angle $\varepsilon$, with
$0<\varepsilon <\pi/2$. We are interested in proving the
asymptotic expansion \eqref{1.1} and especially in its higher-order
analogs in powers of $1/z$, under optimal smoothness hypotheses on $B$.
Such results are then also derived for the $2m\times 2m$ analog $M(z,x)$
of $M_+(z,x)$ associated with Dirac-type operators on $\bbR$. 
 
Our principal motivation in studying this problem stems from our
general interest in operator-valued Herglotz functions
(cf.~\cite{Ca76}, \cite{GKMT98}, \cite{GKM00}, \cite{GM99}, \cite{GM99a},
\cite{GMN98}, \cite{GMT98}, \cite{GT97},
\cite{Sh71}) and their possible applications in the areas of inverse
spectral theory and completely integrable systems. More precisely,
using higher-order asymptotic expansions of $M_+(z,x)$, one can prove
trace formulas for $B(x)$ and certain higher-order differential
polynomials in $B(x)$ (similar in spirit to an approach pioneered in
\cite{GS96} (see also \cite{GH97}, \cite{GHSZ95}) in connection with
Schr\"odinger operators). These trace formulas, in turn, then can be used
to prove various results in inverse spectral theory for matrix-valued
Dirac-type operators $D=J\f{d}{dx}-B$ in $L^2(\bbR)^{2m}$. For instance,
using one of the principal results of this paper, Theorem~\ref{t4.6}, and
its straightforward application to the asymptotic high-energy
expansion of the diagonal Green's matrix
$G(z,x,x)=(D-z)^{-1}(x,x)$ of $D$, the following matrix-valued
analog of a classical uniqueness result of Borg \cite{Bo46} for
one-dimensional Schr\"odinger operators will be proven in
in the context of Dirac-type operators in Section~\ref{s6}.

\begin{theorem} \lb{t1.1} Suppose that $B$ is of the normal form 
$B(x)=\left(\begin{smallmatrix} B_{1,1}(x)
& B_{1,2}(x)\\B_{1,2}(x) & -B_{1,1}(x)\end{smallmatrix}\right)$, with
$B_{1,1}(x)$ and $B_{1,2}(x)$ self-adjoint for a.e.~$x\in\bbR$, and 
assume that $D$ is reflectionless $($e.g., $B$ is periodic and D has 
uniform spectral multiplicity $2m$$)$. In addition, suppose that $D$ 
has spectrum equal to $\bbR$. Then,
\begin{equation}
B(x)=0 \text{ for a.e.~$x\in\bbR$}. \lb{1.2}
\end{equation}
\end{theorem}

For related results see, for instance, \cite{Am93}, \cite{AG96},
\cite{CJ87}, \cite{Ge91}, \cite{GSS91}, \cite{GJ84},
\cite{GG93}. Incidentally, the higher-order differential polynomials in
$B(x)$ just alluded to represent the Ablowitz-Kaup-Newell-Segur (AKNS) or
Zakharov-Shabat (ZS) invariants (i.e.,  densities associated with the
AKNS-ZS conservation laws) and hence provide a link to
infinite-dimensional completely integrable systems (cf., e.g.,
\cite{AK90}, \cite{Ch96}, \cite{Di91}, \cite{Du77}, \cite{Du83},
\cite{Ma78}, \cite{Ma88}, \cite{Sa88}, \cite{Sa99}, \cite{Sa94a},
\cite{Sa99a}, and the references therein), especially, hierarchies of
matrix-valued (i.e., nonabelian) nonlinear Schr\"odinger equations. 

Although various aspects of inverse spectral theory for scalar 
Schr\"odinger, Jacobi, and Dirac-type operators, and more
generally, for $2\times 2$ Hamiltonian systems, are 
well-understood by now (cf.~the extensive list of references provided in
\cite{GKM00}), the corresponding theory for such  operators and
Hamiltonian systems with $m\times m$ matrix-valued coefficients,
$m\in\bbN$, is still in its infancy. A particular  inverse spectral
theory aspect we have in mind is that of determining isospectral sets
(manifolds) of such systems. It may,  perhaps, come as a surprise that
determining the isospectral set  of Hamiltonian systems with
matrix-valued periodic coefficients  is a completely open problem. It
appears to be no exaggeration to  claim that absolutely
nothing seems to be known  about the corresponding isospectral sets of
periodic Dirac-type operators in the case $m\geq 2$.\ (More or less the
same ignorance applies to Schr\"odinger, Jacobi, and more generally, to 
periodic $2m\times 2m$ Hamiltonian systems with $m\geq 2$.)
Theorem~\ref{t1.1} can be viewed  as a first (and very modest) step
toward the construction of isospectral manifolds of certain classes of
matrix-valued potential coefficients $B$  for Dirac-type operators.

However, asymptotic high-energy expansions for Weyl-Titchmarsh matrices
on half-lines and on $\bbR$, their applications to trace formulas for
$B(x)$, and the derivation of Borg-type theorems for Dirac operators
are not the only topics under consideration in this paper. We also provide a
comprehensive and new treatment of local uniqueness theorems for $B$
in terms of exponentially close Weyl-Titchmarsh matrices. More
precisely, in Section \ref{s5} we will prove the following
result ($\|\cdot\|_{\bbC^{m\times m}}$ denotes a matrix norm on
$\bbC^{m\times m}$).

\begin{theorem} \lb{t1.2} 
Fix $x_0\in\bbR$ and suppose that $B_j\in
L^1([x_0,x_0+R])^{2m\times 2m}$ for all $R>0$, posseses the normal form
given in Theorem~\ref{t1.1}
a.e.~on $(x_0,\infty)$, $j=1,2$. Denote by $M_{j,+}(z,x)$, $x\geq
x_0$, the unique Weyl-Titchmarsh matrices corresponding to the half-line
Dirac-type operators in $L^2([x_0,\infty))^{2m}$ associated with
$B_j$, $j=1,2$ $($fixing some self-adjoint boundary condition at $x_0$$)$.
Then,
\begin{equation}
\text{if for some $a>0$, }\, B_1(x)=B_2(x) \, \text{ for a.e. $x\in
(x_0,x_0+a)$,} \lb{1.3} 
\end{equation}
one obtains 
\begin{equation}
\|M_{1,+}(z,x_0)-M_{2,+}(z,x_0)\|_{\bbC^{m\times
m}}\underset{\substack{\abs{z}
\to\infty\\ z\in \rho_{+}}}{=}
O\big(e^{-2\Im(z)a}\big) \lb{1.4}
\end{equation}
along any ray $\rho_+\subset\bbC_+$ with $0<\arg(z)<\pi$ $($and for all 
self-adjoint boundary condition at $x_0$$)$. Conversely, if $m>1$, assume
in addition that $B_j\in L^\infty([x_0,x_0+a])^{2m\times 2m}$, $j=1,2$.
Moreover, suppose that for all $\varepsilon >0$, 
\begin{equation}
\|M_{1,+}(z,x_0)-M_{2,+}(z,x_0)\|_{\bbC^{m\times
m}}\underset{\substack{|z|\to\infty\\z\in \rho_{+,\ell}}}{=}
O\big(e^{-2\Im(z)(a-\varepsilon)}\big), \quad \ell=1,2, \lb{1.5}
\end{equation}
along a ray $\rho_{+,1}\subset\bbC_+$ with $0<\arg(z)<\pi/2$ and along a
ray $\rho_{+,2}\subset\bbC_+$ with $\pi/2<\arg(z)<\pi$. Then 
\begin{equation}
B_1(x)=B_2(x) \text{ for a.e. } x\in [x_0,x_0+a].  \lb{1.6}
\end{equation}
\end{theorem}
\noindent We also prove the analog of Theorem~\ref{t1.2} for the $2m\times 2m$
Weyl-Titchmarsh matrices $M_j(z,x)$ associated with Dirac-type
operators on $\bbR$ corresponding to $B_j$, $j=1,2$. 

In the context of scalar Schr\"odinger operators, the analog of 
Theorem~\ref{t1.2} was first proved by Simon \cite{Si98}. An alternative
proof, applicable to matrix-valued Schr\"odinger operators was presented
in \cite{GS99} (cf.~also \cite{GKM00}). More recently, yet another proof
was found by Bennewitz \cite{Be00} (following  some ideas in
\cite{Bo52}). In fact, our proof of Theorem~\ref{t1.2} is based on that
of Bennewitz \cite{Be00} with additional modifications necessary to
accomodate  Dirac-type operators. These results extend the
classical (global) uniqueness results due to Borg \cite{Bo52} and
Marchenko \cite{Ma50}, \cite{Ma52} which state
that half-line $m$-functions uniquely determine the corresponding
potential coefficient. The Dirac-type results such as Theorem~\ref{t1.2} 
appear to be new, even in the special case $m=1$. Previous results in the
Dirac case focused on global uniqueness questions only. We refer to
Gasymov and Levitan \cite{GL66} in the case $m=1$ and to Lesch and Malamud
\cite{LM00} in the matrix case $m\in\bbN$.

Next, we briefly sketch the content of each section. Section~\ref{s2} 
provides the necessary background results on 
Dirac-type operators and recalls the basic notions of Weyl-Titchmarsh
theory for Hamiltonian systems on a half-line as well
as on $\bbR$, as developed in detail by Hinton and Shaw in a series of
papers \cite{HS81}--\cite{HS86} (see also \cite{At64}, \cite{Jo87},
\cite{HSH93}, \cite{HSH97}, \cite{JNO00}, \cite{KR74}, \cite{KS88},
\cite{Kr89a}, \cite{Kr89b}, \cite{LM00a}, \cite{Or76}, \cite{Sa94a}). In
fact, most of these references deal with more general singular
Hamiltonian systems and hence we specialize some of  this material to the
Dirac-type operator case at hand. While our treatment of
Weyl-Titchmarsh theory in Section~\ref{s2} is  somewhat 
detailed,  the results presented appear to
be of vital importance for our asymptotic expansions in Sections~\ref{s3}
and \ref{s4}. At any rate, we intended to   present this
material as concisely as possible. 

Section~\ref{s3} is devoted to a proof of the
leading-order for the asymptotic high-energy expansion
\eqref{1.1} of $M_+(z,x)$ for the Dirac case. We follow the strategy 
developed in the context of matrix-valued Schr\"odinger operators in our
joint paper \cite{CG99} by appealing to the theory of Riccati equations.
By doing so, we follow the lead of Atkinson who highlighted the
importance of Riccati equations, in this regard, first in \cite{At81},
subsequently in \cite{At82}, \cite{At88a} and ultimately in the
unpublished manuscript \cite{At88} in which he obtains the leading order
for the  asymptotic  high-energy expansion of $M_+(z,x)$ for the
matrix-valued Schr\"odinger case.

Theorems \ref{t3.6} and \ref{t3.7} contain two characterizations 
of the {\it Weyl disk} (cf.~Definition~\ref{dWD}). These 
characterizations provide an answer in Remark~\ref{r3.3} to a point 
raised in \cite{CG99} concerning  the nature of the Weyl disk. 
{}From these characterizations of the Weyl disk, we obtain a 
realization  of $M_+(z,x)$ as a  differentiable function of $x$ which
satisfies a certain  Riccati equation globally and  whose imaginary part
is strictly positive.   We observe, in Remark~\ref{r3.6a}, that the
totality of Weyl disks, $D_+(z,x)$ (cf.~Defintion~\ref{dLWD}), represents
the  phase space for these solutions.  Thus, the asymptotic expansion we
seek, represents  the asymptotic high-energy  behavior for certain
solutions of a given Riccati equation.

Section~\ref{s4} develops a systematic
higher-order high-energy asymptotic expansion of $M_+(z,x)$ as
$|z|\to\infty$, combining the leading-order asymptotic result in
Section~\ref{s3} with matrix-valued extensions of some methods based 
again on an associated Riccati equation. More precisely, following a
technique in \cite{GS98} in the scalar Schr\"odinger operator context, we
show how to derive the general high-energy asymptotic expansion of
$M_+(z,x)$ as $|z|\to\infty$ by combining Atkinson's
leading-order term in \eqref{1.1} and the corresponding asymptotic
expansion of $M_+(z,x)$ in the special case where $B$ has
compact support. Section~\ref{s5} then contains our new local uniqueness
results for $B(x)$ in terms of exponentially small differences of
Weyl-Titchmarsh matrices as indicated in Theorem~\ref{t1.2}. Finally, in 
Section~\ref{s6} we derive a new trace formula for
Dirac-type operators $D$ in $L^2(\bbR)^{2m}$, using
appropriate Herglotz representation results for the diagonal Green's
matrix $G(z,x,x)$ discussed in Section~\ref{s2}. Moreover, we derive the 
Borg-type Theorem~\ref{t1.1} for  Dirac operators and
close with an application to the case of periodic  potentials
coefficients $B$.

\section{Weyl-Titchmarsh Matrices for Hamiltonian Systems}
\lb{s2}

We now turn to the Weyl-Titchmarsh theory for Hamiltonian
systems as developed by Hinton and Shaw in a series
of papers devoted to the spectral theory of (singular) Hamiltonian systems
\cite{HS81}--\cite{HS86} (see also \cite{HSH93}, \cite{HSH97},
\cite{Kr89a}, \cite{Kr89b}, \cite{Sa92}, \cite{Sa94a}, \cite{Sa99},
\cite{Sa99a}). Throughout this paper all matrices will be considered over
the field of complex numbers $\bbC$. The basic assumptions throughout are
described in the following three hypotheses.

\begin{hypothesis} \lb{h2.1}
Fix $m\in\bbN$ and define the $2m\times 2m$ matrix
\begin{subequations}\lb{2.1}
\begin{equation}\lb{2.1a}
J=\begin{pmatrix}0& -I_m \\ I_m & 0  \end{pmatrix},
\end{equation}
where $I_m$ denotes the identity matrix in $\bbC^{m\times m}$.
Suppose 
\begin{equation}
A_{j,k}, B_{j,k} \in
L_{\loc}^1(\bbR)^{m\times m}, \quad j,k = 1,2
\end{equation}
and assume
\begin{align}
A(x)&=\begin{pmatrix}A_{1,1}(x)&A_{1,2}(x)
\\A_{2,1}(x) & A_{2,2}(x) \end{pmatrix}\ge 0, \lb{2.1c}\\ 
\quad B(x)&=\begin{pmatrix}B_{1,1}(x)&B_{1,2}(x) \\B_{2,1}(x) &
B_{2,2}(x)  \end{pmatrix}=B(x)^*,\lb{2.1d}
\end{align}
\end{subequations}
for a.e.~$x\in \bbR$.
\end{hypothesis}
\noindent
$L_{\loc}^1(\bbR)$ denotes the set of locally integrable
functions on $\bbR$. With $M\in\bbC^{m\times m}$, let
$M^t$ denote the transpose, let $M^*$ denote the adjoint or conjugate
transpose of the matrix $M$ and let $M\ge 0$ and $M\le 0$ denote
nonnegative and nonpositive matrices 
$M$ (i.e., positive and negative semidefinite matrices). Moreover, let
$\Im(M)=(M-M^*)/(2i)$ and
$\Re(M)=(M+M^*)/2$ denote, respectively, the   imaginary and real parts 
of the matrix $M$.

Given Hypothesis~\ref{h2.1}, our Hamiltonian system is  given by
\begin{subequations}\lb{HS}
\begin{equation}\lb{HSa}
J \varPsi'(z,x)=(zA(x)+B(x))\varPsi(z,x), \quad z\in\bbC
\end{equation}
for a.e. $x\in \bbR$, where $z$ plays the role of the spectral
parameter, and where
\begin{equation}\lb{HSb}
\varPsi(z,x) = \begin{pmatrix}\psi_1(z,x)\\
\psi_2(z,x) \end{pmatrix}, \quad \psi_j(z,\cdot)\in
AC_{\loc}(\bbR)^{m\times r}, \,\, j=1,2.
\end{equation}
\end{subequations}
$AC_{\loc}(\bbR)$ denotes the set of locally absolutely
continuous functions on $\bbR$. The parameter $r$ in \eqref{HSb}
will be context dependent and range between $1\leq r \leq m$.

For our discussions of the Weyl-Titchmarsh theory for
the Hamiltonian system \eqref{HS}, we  introduce the
definiteness assumption found in Atkinson~\cite{At64}.
\begin{hypothesis}\lb{h2.2}
For all nontrivial solutions $\varPsi$ of \eqref{HSa} with $r=1$ in
\eqref{HSb},  we assume that
\begin{equation}\lb{2.3}
\int_{a}^b dx \, \varPsi(z,x)^*A(x)\varPsi(z,x) > 0\; ,
\end{equation}
for every interval $(a,b)\subset \bbR$, $a<b$.
\end{hypothesis}
A principal example of such a system is the Dirac-type
system obtained when
\begin{equation}\lb{DS}
A(x) = I_{2m},
\end{equation}
and the subject of the present paper; another example being
the matrix-valued Schr\"odinger system, obtained when
\begin{equation}\lb{SS}
A(x)= \begin{pmatrix}I_m& 0 \\ 0 & 0  \end{pmatrix},
\qquad B(x) = \begin{pmatrix}-Q(x)& 0 \\ 0 & I_m  \end{pmatrix},
\end{equation}
and the subject of \cite{CG99}. When \eqref{SS} holds, we
note that
\eqref{HSa} is equivalent to
\begin{align}
-\psi_1^{\prime\prime}(z,x)+Q(x)\psi_1(z,x)&
=z\psi_1(z,x), \lb{2.6} \\
\psi_2(z,x)&=\psi_1^\prime(z,x) \lb{2.7}
\end{align}
for a.e.~$x\in\bbR$. Hypothesis \ref{h2.2} clearly holds in both 
examples.

Next, we introduce a set of matrices that will serve as
boundary data for separated boundary conditions.
\begin{hypothesis}\lb{h2.3}
Let $\gamma = (\gamma_1\; \gamma_2)$ with $\gamma_j \in \bbC^{m\times m}$,
$j=1,2$.  We assume that $\gamma$ satisfies the following
conditions,
\begin{subequations}\lb{BD}
\begin{equation}\lb{BDa}
\rank (\gamma)  = m,
\end{equation}
and that either
\begin{equation}\lb{BDc}
\Im (\gamma_{2}\gamma_{1}^*) \le 0, 
\quad \text{or}
\quad \Im (\gamma_{2}\gamma_{1}^*) \ge 0,
\end{equation}
where $(2i)^{-1}\, \gamma J\gamma^*=\Im (\gamma_{2}\gamma_{1}^*)$.
Given the rank condition in \eqref{BDa},
we assume, without loss of generality in what follows, the
normalization 
\begin{equation}\lb{BDd}
\gamma\gamma^*  = I_m.
\end{equation}
\end{subequations}
\end{hypothesis}
\begin{remark} \lb{r2.4}
With $\alpha\in\bbC^{m\times 2m}$, the conditions 
\begin{equation}
\alpha\alpha^*=I_m, \quad \alpha J\alpha^*=0 \lb{2.8e}
\end{equation}
imply that $\alpha$ satisfies Hypothesis~\ref{2.3}, and they explicitly
read
\begin{equation}
\alpha_1\alpha_1^* +\alpha_2\alpha_2^*=I_m, \quad \alpha_2\alpha_1^* 
-\alpha_1\alpha_2^*=0. \lb{2.8f}
\end{equation}
In fact, one also has
\begin{equation}
\alpha_1^*\alpha_1 +\alpha_2^*\alpha_2=I_m, \quad \alpha_2^*\alpha_1 
-\alpha_1^*\alpha_2=0, \lb{2.8g}
\end{equation}
as is clear from
\begin{equation}
\begin{pmatrix} \alpha_1 & \alpha_2\\ -\alpha_2 & \alpha_1 \end{pmatrix}
\begin{pmatrix} \alpha_1^* & -\alpha_2^*\\ \alpha_2^* & \alpha_1^*
\end{pmatrix}=I_{2m}=\begin{pmatrix} \alpha_1^* & -\alpha_2^*\\ 
\alpha_2^* & \alpha_1^* \end{pmatrix}\begin{pmatrix} \alpha_1 & \alpha_2\\
-\alpha_2 & \alpha_1 \end{pmatrix}, \lb{2.8h}
\end{equation}
since any left inverse matrix is also a right inverse, and vice versa. 
Moreover, from \eqref{2.8g} we obtain
\begin{equation}\lb{2.8i}
\alpha^*\alpha J + J\alpha^*\alpha =  J.
\end{equation}
\end{remark}

With $\alpha\in\bbC^{m\times 2m}$ satisfying \eqref{2.8e}, 
let $\Psi(z,x,x_0,\alpha)$ be a normalized fundamental system of
solutions of \eqref{HS} at some
$x_0\in\bbR$. That is, $\Psi(z,x,x_0,\alpha)$ satisfies \eqref{HS}
for a.e.\ $x\in\bbR$, and
\begin{subequations}\lb{FS}
\begin{equation}\lb{FSa}
\Psi(z,x_0,x_0,\alpha)=(\alpha^* \; J\alpha^*)=
\begin{pmatrix} \alpha_1^* & -\alpha_2^* \\
\alpha_2^* & \alpha_1^* \end{pmatrix}.
\end{equation}
We partition $\Psi(z,x,x_0,\alpha)$ as follows,
\begin{align}
\Psi(z,x,x_0,\alpha)&=(\Theta(z,x,x_0,\alpha)\;
\Phi(z,x,x_0,\alpha))\lb{FSb}\\
&=\begin{pmatrix}\theta_1(z,x,x_0,\alpha) &
\phi_1(z,x,x_0,\alpha)\\
\theta_2(z,x,x_0,\alpha)& \phi_2(z,x,x_0,\alpha)
\end{pmatrix},\lb{FSc}
\end{align}
\end{subequations}
where $\theta_j(z,x,x_0,\alpha)$ and $\phi_j(z,x,x_0,\alpha)$
for $j=1,2$ are $m\times m$ matrices, entire with
respect to $z\in\bbC$, and normalized according to
\eqref{FSa}.~One can now prove the following result.
\begin{lemma}\lb{l2.4}
Let $ \Theta(z,x,x_0,\alpha)$ and $\Phi(z,x,x_0,\alpha)$
be defined in \eqref{FSb} with $\alpha$ and $\beta$ satisfying
Hypothesis~\ref{h2.3} and with $\Im(\alpha_2\alpha_1^*)=0$. Then, for
$c\ne x_0$, $\beta\Phi(z,c,x_0,\alpha)$ is singular  if and
only if $z$ is an eigenvalue for the  regular  boundary value
problem given by \eqref{HSa} on $[x_0,c]$ if $c>x_0$ and on
$[c,x_0]$ if $c<x_0$ together with the separated boundary conditions
\begin{equation}\lb{BC}
\alpha\varPsi(z,x_0)=0, \quad \beta\varPsi(z,c)=0,
\end{equation}
where $\varPsi (z,x)=(\psi_1(z,x)^t\; \psi_2(z,x)^t)^t$ with
$\psi_j(z,\cdot)\in AC([x_0,c])$ if $c>x_0$ and $\psi_j(z,\cdot)\in
AC([c,x_0])$ if $c<x_0$, $j=1,2$. 
\end{lemma}
\noindent Note that the regular boundary value problem described
in Lemma~\ref{l2.4} is self-adjoint when $\Im(\beta_2\beta_1^*)=0$. 

In light of Lemma~\ref{l2.4}, it is possible to introduce, under appropriate
conditions, the $m\times m$ matrix-valued function,
$M(z,c,x_0,\alpha,\beta)$,  as follows.
\begin{definition}\lb{dMF}
Let $ \Theta(z,x,x_0,\alpha)$, and $\Phi(z,x,x_0,\alpha)$ be
defined in \eqref{FSb} with $\alpha$ and $\beta$ satisfying
Hypothesis~\ref{h2.3} and with $\Im(\alpha_2\alpha_1^*)=0$. For
$c\ne x_0$, and $\beta\Phi(z,c,x_0,\alpha)$ nonsingular let
\begin{equation}\lb{MF}
M(z,c,x_0,\alpha,\beta) =
-[\beta\Phi(z,c,x_0,\alpha)]^{-1}[\beta\Theta(z,c,x_0,\alpha)].
\end{equation}
$M(z,c,x_0,\alpha,\beta)$ is said to be the
{\it Weyl-Titchmarsh $M$-function} for the regular boundary value
problem described in Lemma~\ref{l2.4}.
\end{definition}

The Weyl-Titchmarsh $M$-function  is an
$m\times m$ matrix-valued function with meromorphic entries whose
poles correspond to eigenvalues for the regular boundary value
problem given by \eqref{HSa} and \eqref{BC}. Moreover, if $M\in
\bbC^{m\times m}$, and one defines
\begin{equation}\lb{2.14}
U(z,x,x_0,\alpha)= \begin{pmatrix}
u_1(z,x,x_0,\alpha)\\u_2(z,x,x_0,\alpha) \end{pmatrix}=
\Psi(z,x,x_0,\alpha)\begin{pmatrix}I_m\\M\end{pmatrix},
\end{equation}
with $u_j(z,x,x_0,\alpha)\in \bbC^{m\times m}$, $j=1,2$, then
$U(z,x,x_0,\alpha)$ will satisfy the boundary condition
at $x=c$ in
\eqref{BC} whenever $M=M(z,c,x_0,\alpha,\beta)$. 
Intimately connected with the matrices 
introduced in Definition~\ref{dMF} is the set of
$m\times m$ complex matrices known as the Weyl disk. Several
characterization of this set have appeared in the literature 
(see, e.g., \cite{At64}, \cite{At88a}, \cite{At88}, \cite{HSH93},
\cite{HS81}, \cite{Kr89a},  \cite{Or76}). We now mention two, and 
will introduce two others in Section~\ref{s3}   which we use
in the derivation of the asymptotic expansions that are the subject 
of Sections \ref{s3} and \ref{s4}. 

To describe this set, we first introduce the matrix-valued function 
$E_c(M)$: With $c\ne x_0$, $z\in\bbC\backslash\bbR$, and with
$U(z,c,x_0,\alpha)$  defined by \eqref{2.14} in terms of a  matrix
$M\in\bbC^{m\times m}$, let
\begin{equation}\lb{2.380}
E_c(M) = \sigma(x_0,c,z)U(z,c,x_0,\alpha)^*(iJ)U(z,c,x_0,\alpha),
\end{equation}
where 
\begin{equation}
\sigma(s,t,z)=\frac{(s-t)\Im (z)}{|(s-t)\Im (z)|},\quad
\sigma(s,t)=\sigma(s,t,i),\quad
\sigma(z)= \sigma(1,0,z),
\end{equation}
with $s\ne t$, and  $s,t\in\bbR$.
\begin{definition}\lb{dWD}
Let the following be fixed: Real numbers $x_0$ and $c\ne x_0$,
an   $m\times 2m$ matrix $\alpha$ satisfying \eqref{2.8e},  
and $z\in\bbC\backslash\bbR$.
$\calD(z,c,x_0,\alpha)$ will denote the collection of all $M\in
\bbC^{m\times m}$ for which $E_c(M)\le 0$, where $E_c(M)$ is
defined in
\eqref{2.380}.  $\calD(z,c,x_0,\alpha)$ is said to be a
{\it Weyl disk}.
The set of $M\in \bbC^{m\times m}$ for which $E_c(M) = 0$ is
said to be
a {\it Weyl circle} (even when $m>1$).
\end{definition}
This definition leads to a  presentation  
that is a generalization of the description first given 
by Weyl~\cite{We10}; a presentation which is geometric in nature, 
involves the contractive matrices
$V\in\bbC^{m\times m}$, such that $VV^*\le I_m$, and provides the
justification for the geometric terms of circle and disk 
(cf., e.g., \cite{HS81}, \cite{HSH93}, \cite{Kr89a}, \cite{Or76}).

The disk has also been characterized 
in terms of matrices which statisfy Hypothesis~\ref{h2.3} and which 
serve as boundary data for the regular boundary value problem 
described in Lemma~\ref{l2.4} (cf., e.g., \cite{At88a}, \cite{At88}).
More precisely, one could have used the following alternative
definition. \\

\vspace*{-2mm}
\noindent {\bf Definition 2.7A.} 
$\calD(z,c,x_0,\alpha)$ denotes the collection of all $M\in
\bbC^{m\times m}$ obtained by the construction given in \eqref{MF}
where $c\ne x_0$, $z\in\bbC/\bbR$, where $\alpha$ and $\beta$ are
the $m\times m$ matrices defined in Hypothesis~\ref{h2.3} for which
$\sigma(c,x_0,z)\Im(\beta_2\beta_1^*)\ge 0$, and 
$\Im{(\alpha_2\alpha_1^*)}=0$. \\

\vspace*{-2mm}
\noindent However, in this paper we take 
Definition~\ref{dWD} as our point of departure. 

We note that the Weyl
circle corresponds to the regular boundary value problems in
Lemma~\ref{l2.4} with separated, self-adjoint boundary conditions. For
convenience of the reader, and to achieve a reasonable level of
completeness, we reproduce the corresponding short proof below. 

\begin{lemma}[\cite{HS84}, \cite{HSH93}, \cite{Kr89a}]\lb{l2.11}
Let $M\in\bbC^{m\times m}$, $c\ne x_0$, and $z\in\bbC\backslash\bbR$.
Then, $E_c(M)=0$ if and only if
there is a $\beta\in \bbC^{m\times 2m}$ satisfying
\eqref{2.8e} such that
\begin{equation}\lb{2.24}
0=\beta U(z,c,x_0,\alpha),
\end{equation}
where $U(z,c)=U(z,c,x_0,\alpha)$ is defined in \eqref{2.14} in
terms of $M$.  With $\beta$ so defined,
\begin{equation}\lb{2.25}
M=-[\beta\Phi(z,c,x_0,\alpha)]^{-1}[\beta\Theta(z,c,x_0,\alpha)],
\end{equation}
that is, $M=M(z,c,x_0,\alpha,\beta)$. Moreover, $\beta$ may be chosen
to satisfy \eqref{BDd}, and hence Hypothesis~\ref{2.3}.
\end{lemma}
\begin{proof}
Let $z\in\bbC\backslash\bbR$, and suppose
for a given $M\in\bbC^{m\times m}$ that there is
a $\beta\in\bbC^{m\times 2m}$ which satisfies \eqref{2.8e}
and such that \eqref{2.24} is satisfied.  Given that
$\beta J \beta^* = 2i\Im (\beta_2\beta_1^*)=0$, and given
that $\rank (\beta)=\rank(J\beta^*)=m$, there is a
nonsingular $C\in\bbC^{m\times m}$ such that
$U(z,c) = J\beta^*C$.  Hence,
$E_c(M)=i\sigma(c,x_0,z)C^*\beta J\beta^*C=0$.

Upon showing that $\beta\Phi(z,c)=\beta\Phi(z,c,x_0,\alpha)$ 
is nonsingular, \eqref{2.25}
will then follow from \eqref{2.24}. If $\beta\Phi(z,c)$ is singular,
then there are nonzero vectors
$v, w \in \bbC^{m}$ such that $\beta\Phi(z,c)v=0$, and such that
$\Phi(z,c)v = J\beta^*w$.
Let $\varPsi_j=\varPsi_j(z,x)$, $j=1,2$ denote solutions of
\eqref{HSa} with $z=z_j$, $j=1,2$.  Then,
\begin{equation}\lb{2.19}
(\varPsi_1^*J\varPsi_2)'=(z_2 - \bar{z}_1)\varPsi_1^*A\varPsi_2.
\end{equation}
Using \eqref{2.19}, and recalling that $\Phi(z,x)$ is defined in \eqref{FS}, we obtain
\begin{subequations}\lb{2.230}
\begin{align}
2i\Im(z)\int_{x_0}^c dx\, v^*\Phi(z,x)^* A(x) \Phi(z,x)v&=v^*\Phi(z,c)^* J
\Phi(z,c)v \lb{2.230a} \\
&=w^*\beta J \beta^*w =0.
\end{align}
\end{subequations}
Thus, by Hypothesis~\ref{h2.2}, $\Im(z)=0$. This contradicts the
assumption that $z\in\bbC\backslash\bbR$. 

Conversely, if $E_c(M)=0$ for a given $M\in \bbC^{m\times m}$,
then for $z\in\bbC\backslash\bbR$  let
$\beta =(I_m\; M^*)\Psi(z,c,x_0,\alpha)^*J =
  U(z,c,x_0,\alpha)^* J$. One observes that \eqref{2.24} is
satisfied and that  $\rank (\beta) =m$. Moreover,
$0=\sigma(x_0,c,z)E_c(M)/2=\Im(\beta_2\beta_1^*)$.
If for this choice of $\beta$,  \eqref{BDd} is not yet
satisfied,  one introduces $\Delta = (\beta\beta^*)^{-1/2}\beta$ and
observes that $0=\Delta U(z,c,x_0,\alpha)$,
$\Im(\Delta_2\Delta_1^*)=
(\beta\beta^*)^{-1/2} \Im(\beta_2\beta_1^*)
(\beta\beta^*)^{-1/2}$, and that $\Delta$ satisfies all
requirements of \eqref{2.8e}.
\end{proof}

Next, we recall a fundamental property associated with matrices
in $\calD(z,c,x_0,\alpha)$.
\begin{lemma} \lb{l2.8}
If $M\in\calD(z,c,x_0,\alpha)$, then 
\begin{equation}\lb{Hgz}
\sigma(c,x_0,z) \Im (M)>0.
\end{equation}
Moreover, whenever $\beta\in\bbC^{m\times 2m}$ satisfies \eqref{2.8e},
\begin{equation}\lb{2.270}
M(\bar z,c,x_0,\alpha,\beta)= M(z,c,x_0,\alpha,\beta)^*.
\end{equation}
\end{lemma}
\begin{proof}
Let $\varPsi_j=\varPsi_j(z,x)$, $j=1,2$ denote solutions of 
\eqref{HSa} with $z=z_j$, $j=1,2$.  Then $(\varPsi_1^*J\varPsi_2)'=(z_2
- \bar{z}_1)\varPsi_1^*A\varPsi_2$ as in \eqref{2.19}. This implies
\begin{align}
2i\Im(z)\int_{x_0}^c dx\, U(z,x)^* A(x) U(z,x) &=
U(z,x)^* J U(z,x)\big |_{x_0}^c \no \\
&=2i\Im(M) + U(z,c)^* J U(z,c), \lb{2.280} 
\end{align}
with $U(z,x)=U(z,x,x_0,\alpha)$ defined in \eqref{2.14}. Moreover, by the
definition of $E_c(M)$ given in \eqref{2.380}, one obtains
\begin{align}\lb{2.290}
&2\sigma(c,x_0,z)\Im (M)\\
&= -E_c(M) + 2\sigma(c,x_0)|\Im(z)|\int_{x_0}^c ds\, U(z,s)^*A(s)U(z,s).\no
\end{align}
By Hypothesis~\ref{h2.2} and Definition~\ref{dWD}, one infers that 
$\sigma(c,x_0,z) \Im (M)>0$. To prove \eqref{2.270}, let
$\Psi(z,x)=\Psi(z,x,x_0,\alpha)$, where
$\Psi$ is defined in \eqref{HS}. Then, by \eqref{2.19}, 
\begin{equation}\lb{2.310}
\Psi(\bar{z},x)^*J\Psi(z,x)=J,
\end{equation}
which implies $J\Psi(z.x)(\Psi(\bar{z},x)J)^*=I_{2m}$ 
and hence
\begin{equation}\lb{2.330}
\Psi(z,x)J\Psi(\bar{z},x)^*=J.
\end{equation}
Thus one concludes
\begin{equation}
\beta\Phi(z,c)(\beta\Theta(\bar{z},c))^*-
\beta\Theta(z,c)(\beta\Phi(\bar{z},c))^*=
\beta J\beta^*=0,
\end{equation}
from which  \eqref{2.270} follows immediately by Lemma~\ref{l2.11}.
\end{proof}

For $c>x_0$, the function $M(z,c,x_0,\alpha,\beta)$, 
defined by \eqref{MF}, and satisfying \eqref{Hgz}, is said to be  
a matrix-valued {\it Herglotz} function of rank $m$. Hence, for
$\Im(\beta_2\beta_1^*)=0$, poles of
$M(z,c,x_0,\alpha,\beta)$, $c>x_0$, are at most of first order, are
real, and have nonpositive residues.  Such functions admit a
representation of the form
\begin{align}\lb{NP}
M(z,c,x_0,\alpha,\beta)=& \; C_1(c,x_0,\alpha,\beta) +
zC_2(c,x_0,\alpha,\beta)
\no\\   
&+\int_{-\infty}^\infty
d\Omega(\lambda,c,x_0,\alpha,\beta)\,\left(
\frac{1}{\lambda-z} -\frac{\lambda}{1+\lambda^2} \right), \quad  c>x_0, 
\end{align}
where $C_2(c,x_0,\alpha,\beta)\ge 0$ and $C_1(c,x_0,\alpha,\beta)$ are 
self-adjoint $m\times m$ matrices, and where
$\Omega(\lambda,c,x_0,\alpha,\beta)$ is a nondecreasing $m\times m$
matrix-valued function such that
\begin{subequations}\lb{Mrep}
\begin{align}
&\int_{-\infty}^{\infty}\|
d\Omega(\lambda,c,x_0,\alpha,\beta)\|_{\bbC^{m\times m}}\,
(1+\lambda^2)^{-1}  < \infty, \lb{NPa} \\
&\Omega((\lambda, \mu],c,x_0,\alpha,\beta)= \lim_{\delta\downarrow 0}
\lim_{\epsilon \downarrow 0}\frac{1}{\pi}\int_{\lambda
+ \delta}^{\mu
+ \delta }d\nu\,  \sigma(c,x_0)\Im\left( M(\nu
+i\epsilon,c,x_0,\alpha,\beta)\right). \lb{NPb}
\end{align}
\end{subequations} 
In general, for self-adjoint boundary value problems, 
$\Omega(\lambda,c,x_0,\alpha,\beta)$ is
piecewise constant with jump discontinuities precisely at the eigenvalues
of the boundary value problem, and that in the matrix-valued Schr\"odinger and 
Dirac-type cases $C_2=0$ in \eqref{NP} (and later in \eqref{2.42} and
\eqref{2.64}). Analogous statements apply to
$-M(z,c,x_0,\alpha,\beta)$ if $c<x_0$. For such problems, we note in the 
subsequent lemma that for fixed $\beta$, varying the boundary data
$\alpha$ produces Weyl-Titchmarsh matrices $M(z,c,x_0,\alpha,\beta)$
related  to each other via linear fractional transformations (see also
\cite{GMT98}, \cite{GT97} for a general approach to such linear fractional
transformations).
\begin{lemma}\lb{l2.9}
Suppose $\alpha, \beta, \gamma\in\bbC^{m\times 2m}$ satisfy \eqref{2.8e}.
Let $M_{\alpha}=M(z,c,x_0,\alpha,\beta)$, and 
$M_{\gamma}=M(z,c,x_0,\gamma,\beta)$. Then,
\begin{equation}
M_{\alpha}= [-\alpha J \gamma^* + \alpha\gamma^* M_{\gamma}]
[\alpha\gamma^* +  \alpha J\gamma^*M_{\gamma}]^{-1}. \lb{2.360}
\end{equation}
\end{lemma}
\begin{proof}
Let $U_{\alpha}(z,x)=U(z,x,x_0,\alpha)$  and $U_{\gamma}(z,x)=U(z,x,x_0,\gamma)$
be defined in \eqref{2.14} with $M=M_{\alpha}$  and  $M=M_{\gamma}$
respectively. Then,
\begin{equation}
0=\beta U_{\alpha}(z,c)=\gamma U_{\gamma}(z,c).
\end{equation}
By the rank condition \eqref{BDa},
\begin{equation}
U_{\alpha}(z,c)= J\beta^*C_{\alpha}\, ,\qquad
U_{\gamma}(z,c)= J\beta^*C_{\gamma}
\end{equation}
for nonsingular $C_{\alpha}, \; C_{\gamma}\in \bbC^{m\times m}$.
Thus, by \eqref{FSa}, and by the uniqueness of solution of \eqref{HSa},
there is a nonsingular 
$C\in \bbC^{m\times m}$ for which
\begin{equation}\lb{2.410}
\begin{pmatrix} \alpha^*\hspace{-5pt}&J\alpha^*\end{pmatrix}
\begin{pmatrix} I_m \\ M_{\alpha}\end{pmatrix}=U_{\alpha}(z,x_0)=U_{\gamma}(z,x_0)C=
\begin{pmatrix} \gamma^*\hspace{-5pt}&J\gamma^*\end{pmatrix}
\begin{pmatrix} I_m \\ M_{\gamma}\end{pmatrix}C.
\end{equation}
By \eqref{2.8i},
\begin{equation}
\begin{pmatrix} \alpha^*\hspace{-5pt}&J\alpha^*\end{pmatrix}^{-1}=
\begin{pmatrix} \alpha \\ -\alpha J\end{pmatrix};
\end{equation}
and hence, by \eqref{2.410} we see that
\begin{subequations}
\begin{align}
I_m&=(\alpha \gamma^* + \alpha J \gamma^* M_{\gamma} )C\\
M_{\alpha}&= (-\alpha J\gamma^* + \alpha  \gamma^* M_{\gamma} )C,
\end{align}
\end{subequations}
from which \eqref{2.360} immediately follows.
\end{proof}
\begin{remark}
{}From the proof of the previous lemma one infers, in general,  that 
\begin{equation}
U_{\gamma}(z,x) = 
U_{\alpha}(z,x)(\alpha \gamma^* + \alpha J \gamma^* M_{\gamma} ).
\end{equation}
Moreover,  if $\alpha_0 =(I_m\; 0)$ and $\gamma_0=(0\ I_m)$
one observes, in particular, 
\begin{equation}
M(z,c,x_0,\alpha_0,\beta)=-M(z,c,x_0,\gamma_0,\beta)^{-1}.
\end{equation}
\end{remark}

We  further note that the sets $\calD(z,c,x_0,\alpha)$ are closed, 
and convex, (cf., e.g., \cite{HS84}, \cite{HSH93}, 
\cite{Kr89a}, \cite{Or76}). Moreover, by \eqref{2.290} and Hypothesis~\ref{h2.2},
one concludes that $E_c(M)$ is strictly increasing. This fact together
with Lemma~\ref{l2.11} implies that, as a function of $c$, the sets
$\calD(z,c,x_0,\alpha)$ are strictly nesting in the sense that
\begin{equation}\lb{2.28}
\calD(z,c_2,x_0,\alpha)\subset \calD(z,c_1,x_0,\alpha) \quad
\text{for}\quad x_0<c_1< c_2\quad \text{or}
\quad c_2< c_1<x_0.
\end{equation}
Hence, the intersection of this nested sequence, as $c\to \pm \infty$,
is nonempty, closed and convex.  We say that this intersection is
a limiting set for the nested sequence.
\begin{definition}\lb{dLWD}
Let  $\calD_\pm(z,x_0,\alpha)$ denote the closed, convex set in the space
of $m\times m$ matrices which is the limit, as
$c\to \pm\infty$,  of the nested collection of sets
$\calD(z,c,x_0,\alpha)$ given in Definition~\ref{dWD}.
$\calD_\pm(z,x_0,\alpha)$ is said to be a limiting {\em disk}.
Elements of $\calD_\pm(z,x_0,\alpha)$ are denoted by $M_\pm(z,x_0,\alpha)\in
\bbC^{m\times m}$.
\end{definition}

In light of the containment described in
\eqref{2.28},  for $c\ne x_0$ and $z\in \bbC\backslash\bbR$,
\begin{equation}\lb{2.32}
\calD_\pm(z,x_0,\alpha)\subset \calD(z,c,x_0,\alpha), 
\end{equation}
with emphasis on strict containment of the disks in \eqref{2.32}.
Moreover, by \eqref{2.290},  
\begin{equation}\lb{2.320}
M\in \calD_\pm(z,x_0,\alpha) \text{ precisely when }E_c(M)<0 
\text{ for all }  c \in(x_0, \pm\infty).
\end{equation}
The following Lemma appears to have gone unnoted in the literature.
\begin{lemma}\lb{l2.12}
Let $M\in\bbC^{m\times m}$, $c\ne x_0$, and $z\in\bbC\backslash\bbR$.
Then, $E_c(M)<0$ if and only if
there is a $\beta\in \bbC^{m\times 2m}$   satisfying the condition
\begin{equation}\lb{2.27a}
\sigma(c,x_0,z) \Im(\beta_2\beta_1^*)>0, 
\end{equation}
and such that \eqref{2.24} holds 
with $u_j(z,c)= u_j(z,c,x_0,\alpha)$, $j=1,2$, defined in
\eqref{2.14} in terms of $M$. With $\beta$ so defined,
\eqref{2.25} holds; that is, $M=M(z,c,x_0,\alpha,\beta)$. Moreover,
$\beta$ maybe chosen to satisfy \eqref{BDd}, and hence Hypothesis~\ref{h2.3}.
\end{lemma}
\begin{proof}
Let $z\in \bbC\backslash\bbR$, and for a 
given $M\in \bbC^{m\times m}$  suppose that there is a
$\beta\in\bbC^{m\times 2m}$ satisfying 
\eqref{2.27a} such that \eqref{2.24} holds. The matrices 
$\beta_j$, $j=1,2$, are
invertible by \eqref{2.27a}, and by \eqref{2.24} it follows that
\begin{equation}\lb{2.26}
U(z,c)=\begin{pmatrix}
-\beta_1^{-1}\beta_2\\I_m \end{pmatrix}
u_2(z,c).
\end{equation}
By \eqref{2.380} and \eqref{2.26}, one then concludes that
\begin{equation}\lb{2.27}
E_c(M) = -2\sigma(c,x_0,z) u_2(z,c)^*\beta_1^{-1}
 \Im (\beta_2\beta_1^*)
(\beta_1^*)^{-1} u_2(z,c),
\end{equation}
and hence that $E_c(M)<0$ whenever \eqref{2.27a} holds. 

Upon showing that $\beta\Phi(z,c)$ is nonsingular, \eqref{2.25} will follow from
\eqref{2.24}. If $\beta\Phi(z,c)$ is singular, then there is a nonzero
vector $v\in \bbC^{m}$ such that $\beta\Phi(z,c)v=0$. By the
nonsingularity of $\beta_j$, $j=1,2$, $\phi_1(z,c)v =
-\beta_1^{-1}\beta_2\phi_2(z,c)v$, and as a result, \eqref{2.230a}
yields
\begin{align}
&2\sigma(c,x_0)|\Im(z)|\int_{x_0}^c dx \,v^*\Phi(z,x)^*A(x)\Phi(z,x)v 
\no\\
&= -2\sigma(c,x_0,z)v^*\phi_2(z,c)^*\beta_1^{-1}\Im(\beta_2\beta_1^*)
(\beta_1^*)^{-1}\phi_2(z,c)v,
\end{align}
and hence, a contradiction given \eqref{2.27a} (cf.~\eqref{2.3}). 

Conversely, if $E_c(M)<0$  for  a given $M\in\bbC^{m\times m}$, then for
$z\in\bbC\backslash\bbR$,  $u_j(z,c)$, $j=1,2$,  defined by
\eqref{2.14}, are nonsingular. Indeed, if either $u_1(z,c)$ 
or $u_2(z,c)$
are singular, then there is a  $v\in \bbC^m$, $v\ne 0$,
such that  $v^*E_c(M)v=0$, a contradiction.  Next, let $\beta_1 =I_m$
and let $\beta_2=-u_1(z,c)u_2(z,c)^{-1}$.
Then, for these $\beta_j$, $j=1,2$, \eqref{2.24}
holds. Equation~\eqref{2.27} now implies that
$\sigma(c,x_0,z)\Im(\beta_2\beta_1^*)> 0$ for
$c\ne x_0$ and $z\in\bbC\backslash\bbR$.
For this choice, $\beta$ does not satisfy \eqref{BDd}. However, one
can normalize $\beta$ as described in the proof of
Lemma~\ref{l2.11}.
\end{proof}
Hence by Lemma~\ref{l2.12} and \eqref{2.320}, we see that if 
$M\in\calD_\pm(z,x_0,\alpha)$, then for some
$\beta\in\bbC^{m\times 2m}$ satisfying \eqref{2.27a}
\begin{equation}\lb{2.33}
M_\pm(z,x_0,\alpha)=M(z,c,x_0,\alpha,\beta).
\end{equation}
\begin{remark}\lb{r2.14}
To the reader of \cite{CG99}, our study of the high-energy asymptotics of
the Weyl-Titchmarsh $M$-function for matrix-Schr\"odinger operators, 
we offer this cautionary note: 
In \cite{CG99}, $D(z,c,x_0,\alpha)$ represents the set of matrices
characterized by Lemmas~\ref{l2.11} and \ref{l2.12}.
However, the homeomorphism that exists between the contractive matrices 
$V\in\bbC^{m\times m},\ VV^*\le I_m$, and the Weyl disk, $D(z,c,x_0,\alpha)$,
(cf., \cite{HS84}, \cite{HSH93}, \cite{Kr89a}, \cite{Or76}) shows that
those $M\in\bbC^{m\times m}$ characterized in
Lemma~\ref{l2.11} correspond to the set of unitary matrices while those
characterized in Lemma~\ref{l2.12} correspond
to the contractive matrices for which $VV^*<I_m$. Hence, Lemma~\ref{l2.11}
characterizes part of the boundary while Lemma~\ref{l2.12} characterizes
the interior of the Weyl disk as it is defined in Defintion~\ref{dWD}.
As a result, the closure of the set consisting of those 
$M\in\bbC^{m\times m}$ characterized by these two lemmas (i.e., those $M$
which correspond to $VV^*<I_m$, or to $VV^*=I_m$) is the Weyl disk.
Thus, for deriving high-energy asymptotics
for $M_\pm(z,x_0,\alpha)$, it is sufficient to consider the subset
of the Weyl disk consisting of those matrices, $M\in\bbC^{m\times m}$,
characterized in Lemma~\ref{l2.11} and Lemma~\ref{l2.12}. 
This was the approach taken in \cite{CG99}. 
\end{remark}

When $\calD_\pm(z,x_0,\alpha)$ is a singleton matrix, the
system \eqref{HSa} is said to be in the {\it limit point} (l.p.) case at
$\pm\infty$. If $\calD_\pm(z,x_0,\alpha)$ has nonempty interior, then
\eqref{HSa} is said to be in the {\it limit circle} (l.c.) case at
$\pm\infty$. Indeed, for the case $m=1$, the limit point case
corresponds to a point in $\bbC$, whereas the limit circle case
corresponds to $\calD_\pm(z,x_0,\alpha)$ being a disk in $\bbC$.

These apparent geometric properties for the disk
correspond to  analytic properties for the solutions of the
Hamiltonian system \eqref{HSa}.  To recall this correspondence, we
introduce the following spaces in which we assume that
$ -\infty\le a< b \le \infty$,
\begin{subequations}\lb{2.29}
\begin{align}
L_A^2((a,b))&=\bigg\{\phi:(a,b)\to\bbC^{2m}  \bigg| \int_a^b dx\,
(\phi(x),A\phi(x))_{\bbC^{2m}}<\infty \bigg\}, \lb{2.29a}
\\
N(z,\infty)&=\{\phi\in L_A^2((c,\infty)) \mid J\phi^\prime
=(zA+B)\phi
\text{ a.e. on $(c,\infty)$} \}, \lb{2.29b}
\\
N(z,-\infty)&=\{\phi\in L_A^2((-\infty,c)) \mid
J\phi^\prime=(zA+B)\phi
\text{ a.e. on $(-\infty,c)$} \}, \lb{2.29c}
\end{align}
\end{subequations}
for some $c\in\bbR$ and $z\in\bbC$. (Here
$(\phi,\psi)_{\bbC^n}=\sum_{j=1}^n \overline\phi_j\psi_j$
denotes the standard scalar product in $\bbC^n$, abbreviating
$\chi\in\bbC^n$ by
$\chi=(\chi_1,\dots,\chi_n)^t$.)  Both dimensions of the
spaces in \eqref{2.29b} and \eqref{2.29c},
$\dim_\bbC(N(z,\infty))$ and $\dim_\bbC(N(z,-\infty))$,
are constant for $z\in\bbC_\pm=\{\zeta\in\bbC
\mid \pm\Im(\zeta)> 0 \}$ (see, e.g., \cite{At64},
\cite{KR74}). One then observes that  the Hamiltonian
system \eqref{HSa} is in the limit point case at
$\pm\infty$ whenever
\begin{equation}\lb{2.30}
\dim_\bbC(N(z,\pm\infty))=m \text{  for all
$z\in\bbC\backslash\bbR$}
\end{equation}
and in the limit circle case at $\pm\infty$  whenever
\begin{equation}\lb{2.31}
\dim_\bbC(N(z,\pm\infty))=2m \text{  for all $z\in\bbC$.}
\end{equation}

Next we show that the Dirac-type systems considered in this paper
are always in the limit point case at $\pm\infty$. Results of this
type, under varying sets of assumptions on $B(x)$, are well-known
to experts in the field. For instance, in the case $m=1$ and with
$B_{1,2}(x)=B_{2,1}(x)$ this fact can be found in \cite{We71}. For
$B\in C(\bbR)^{2m\times 2m}$ and a more general constant matrix $A$,
this result is proven in \cite{LM00} (their proof, however,
extends to the current $B\in L^1_{\loc} (\bbR)$ case). More generally,
multi-dimensional Dirac operators with $L^2_{\loc} (\bbR^n)$-type
coefficients (and additional conditions) can be found in \cite{LO82}. A
short proof in the case $m=1$ has recently been sent to us by Don Hinton
\cite{Hi99}. For convenience  of the reader we present its elementary 
generalization to $m\in\bbN$ below (see also \cite{Cl94} for a sketch of
such a proof). After completion of this paper we became aware of a recent
preprint by Lesch and Malamud \cite{LM00a} which provides a thorough
study of self-adjointness questions for more general Hamiltonian systems
than those studied in this paper.

\begin{lemma} \lb{l2.15}
The limit point case holds for Dirac-type systems
{\rm (}i.e., for $A=I_{2m}$ in \eqref{HSa}{\rm )} at $\pm \infty$.
\end{lemma}
\begin{proof}
Let $\{y_\ell(z,x)\}_{\ell=1,\dots,k}$ and
$\{w_j(z,x)\}_{j=1,\dots,k^\prime}$ denote bases for
$N(z,\pm\infty)$
and $N(\ol z,\pm\infty)$, respectively. By Theorem~9.11.1 of
Atkinson \cite{At64}, one has $k,k^\prime\geq m$ for
$z\in\bbC\backslash\bbR$. We now assume that $k>m$.

One observes that $\{y_1(z,c),\dots,y_k(z,c)\}$ and
$w_1(\ol z,c),\dots,w_{k^\prime}(\ol z,c)\}$ are linearly
independent in $\bbC^{2m+1}$, where $k+k^\prime\geq 2m+1$.
Consequently, there is some $s\in\{1,\dots,k\}$ and some
$r\in\{1,\dots,k^\prime\}$ such that
\begin{equation}
w_r(\ol z,c)^*Jy_s(z,c)\neq 0. \lb{2.32a}
\end{equation}
By Lagrange's identity,
\begin{equation}
w_r(\ol z,x)^*Jy_s(z,x)=w_r(\ol z,c)^*Jy_s(z,c) \lb{2.33a}
\end{equation}
is constant with respect to $x$. On the other hand, an application
of Cauchy's inequality shows that the left-hand side of \eqref{2.33a} is in
$L^1 ((c,\pm\infty))$. By \eqref{2.32a} one obtains a contradiction
and hence concludes that
\begin{equation}
\dim_\bbC (N(z,\pm\infty))=m. \lb{2.34}
\end{equation}
The analogous argument then also yields
\begin{equation}
\dim_\bbC (N(\ol z,\pm\infty))=m \lb{2.35}
\end{equation}
and hence the limit point property of Dirac-type systems with
$A(x)=I_{2m}$ in \eqref{HSa}.
\end{proof}

Returning to the general case \eqref{HSa}, in either the limit point or
limit circle cases, $M_\pm(z,x_0,\alpha)\in \partial
\calD_{\pm}(z,x_0,\alpha)$ is said to be a {\em half-line
Weyl-Titchmarsh  matrix}. Each such matrix is associated with the
construction of a self-adjoint operator acting on
$L_A^2([x_0,\pm \infty))\cap \AC([x_0,\pm\infty))^{2m}$ for the
Hamiltonian system \eqref{HSa}.  However, for those intermediate cases
where
$m<\dim_\bbC(N(z,\pm\infty))<2m$, Hinton and Schneider have noted that
not every element of $\partial\calD_{\pm}(z,x_0,\alpha)$ is a
half-line Weyl-Titchmarsh matrix, and have characterized those elements
of the boundary that are (cf.~\cite{HSH93}, \cite{HSH97}).

For convenience of the reader we summarize some of the
principal results on half-line Weyl-Titchmarsh  matrices next.
\begin{theorem}
[\cite{AD56}, \cite{Ca76}, \cite{GT97}, \cite{HS81},
\cite{HS82}, \cite{HS86}, \cite{KS88}] \lb{t2.3}
Suppose Hypotheses \ref{h2.1} and \\ \ref{h2.2}. Let
$z\in\bbC\backslash\bbR$,
$x_0\in\bbR$, and denote by $\alpha, \gamma\in\bbC^{m\times 2m}$ matrices
satisfying \eqref{2.8e}. Then, \\
$(i)$ $\pm M_{\pm}(z,x_0,\alpha)$ is an $m\times m$
  matrix-valued Herglotz
function of maximal rank.
In particular,
\begin{gather}
\Im(\pm M_{\pm}(z,x_0,\alpha)) > 0, \quad z\in\bbC_+, \\
M_{\pm}(\overline z,x_0,\alpha)=M_{\pm}(z,x_0,\alpha)^*, \lb{2.38} \\
\rank (M_{\pm}(z,x_0,\alpha))=m,  \\
\lim_{\varepsilon\downarrow 0} M_{\pm}(\lambda+
i\varepsilon,x_0,\alpha) \text{
exists for a.e.\
$\lambda\in\bbR$},\\
\begin{split}\lb{2.41}
M_\pm(z,x_0,\alpha) &= [-\alpha J \gamma^* +
\alpha\gamma^* M_\pm(z,x_0,\gamma)]\times \\
&\quad \times[ \alpha\gamma^*
+  \alpha J \gamma^*M_\pm(z,x_0,\gamma)]^{-1}.
\end{split}
\end{gather}
Local singularities of $\pm M_{\pm}(z,x_0,\alpha)$ and 
$\mp M_{\pm}(z,x_0,\alpha)^{-1}$ are necessarily real and at most of first
order in the sense that 
\begin{align}
&\mp \lim_{\epsilon\downarrow0}
\left(i\epsilon\,
M_{\pm}(\lambda+i\epsilon,x_0,\alpha)\right) \geq 0, \quad \lambda\in\bbR,
\lb{2.24b} \\ 
& \pm \lim_{\epsilon\downarrow0}
\left(\f{i\epsilon}{M_{\pm}(\lambda+i\epsilon,x_0,\alpha)}\right)
\geq 0, \quad \lambda\in\bbR. \lb{2.24c}
\end{align}
$(ii)$  $\pm M_{\pm}(z,x_0,\alpha)$ admit the representations
\begin{align}
&\pm M_{\pm}(z,x_0,\alpha)=F_\pm(x_0,\alpha)+\int_\bbR
d\Omega_\pm(\lambda,x_0,\alpha) \,
\big((\lambda-z)^{-1}-\lambda(1+\lambda^2)^{-1}\big) \lb{2.42} \\
&=\exp\bigg(C_\pm(x_0,\alpha)+\int_\bbR d\lambda \, \Xi_{\pm}
(\lambda,x_0,\alpha)
\big((\lambda-z)^{-1}-\lambda(1+\lambda^2)^{-1}\big)
\bigg), \lb{2.43}
\end{align}
where
\begin{align}
F_\pm(x_0,\alpha)&=F_\pm(x_0,\alpha)^*, \quad \int_\bbR
\|d\Omega_\pm(\lambda,x_0,\alpha)\|_{\bbC^{m\times m}} \,
(1+\lambda^2)^{-1}<\infty,
  \\
C_\pm(x_0,\alpha)&=C_\pm(x_0,\alpha)^*,
\quad 0\le\Xi_\pm(\dott,x_0,\alpha)
\le I_m \, \text{  a.e.}
\end{align}
Moreover,
\begin{align}
\Omega_\pm((\lambda,\mu],x_0,\alpha)&
=\lim_{\delta\downarrow
0}\lim_{\varepsilon\downarrow 0}\f1\pi
\int_{\lambda+\delta}^{\mu+\delta} d\nu \, \Im(\pm
M_\pm(\nu+i\varepsilon,x_0,\alpha)),  \\
\Xi_\pm(\lambda,x_0,\alpha)&=\lim_{\varepsilon\downarrow 0}
\pi^{-1}\Im(\ln(\pm
M_\pm(\lambda+i\varepsilon,x_0,\alpha))) \text{ for a.e.\
$\lambda\in\bbR$}.
\end{align}
$(iii)$  Define the $2m\times m$ matrices
\begin{align}
U_\pm(z,x,x_0,\alpha)&=\begin{pmatrix}u_{\pm,1}(z,x,x_0,\alpha) \\
u_{\pm,2}(z,x,x_0,\alpha)  \end{pmatrix}
=\Psi(z,x,x_0,\alpha)\begin{pmatrix} I_m \\
M_\pm(z,x_0,\alpha) \end{pmatrix}  \no \\
&=\begin{pmatrix}\theta_1(z,x,x_0,\alpha)
& \phi_1(z,x,x_0,\alpha)\\
\theta_2(z,x,x_0,\alpha)
& \phi_2(z,x,x_0,\alpha)\end{pmatrix}
\begin{pmatrix} I_m \\
M_\pm(z,x_0,\alpha) \end{pmatrix}, \lb{2.52}
\end{align}
with $\theta_j(z,x,x_0,\alpha)$, and
$\phi_j(z,x,x_0,\alpha)$, $j=1,2$, defined by \eqref{FSc}. Then, 
\begin{equation}
\Im(M_\pm(z,x_0,\alpha))=\Im(z) \int_{x_0}^{\pm\infty}ds\,
U_\pm(z,s,x_0,\alpha)^* A(s)
U_\pm(z,s,x_0,\alpha).
\end{equation}
\end{theorem}

In the Dirac-type context, where $A=I_{2m}$, the
$m$ columns of $U_\pm (z,\cdot,x_0,\alpha)$ span $N(z,\pm\infty)$.

Up to this point, we focused exclusively on Hamiltonian systems and
neglected the notion of a linear operator associated with \eqref{HS}. We
did this on purpose as the formalism presented thus far is widely
applicable and goes beyond the prime candidates such as Schr\"odinger and
Dirac-type systems. However, in the remainder of this section and for 
the bulk of the material from  Section~\ref{s3} on, we will focus on 
the Dirac-type case. Thus, in addition to Hypotheses~\ref{h2.1}--\ref{h2.3},
which are assumed throughout this paper, we introduce the following 
hypothesis taylored to these occasions. 

\begin{hypothesis}\lb{h2.4}
Assume Hypotheses~\ref{h2.1} and \ref{h2.3} as well as the 
Dirac-type assumption \eqref{DS}.
\end{hypothesis}

Assuming the Dirac-type Hypothesis~\ref{h2.4}, we now describe 
the associated  Dirac-type operator $D$ on $\bbR$ by first introducing the
Green's matrix associated with \eqref{HS} and \eqref{DS}. Define the
$2m\times 2m$ matrix $G$ by
\begin{align}
G(z,x,x^\prime)=U_\mp(z,x,x_0,\alpha_0)[M_-(z,x_0,\alpha_0) &
-M_+(z,x_0,\alpha_0)]^{-1} U_\pm(\overline z,x^\prime,x_0,\alpha_0)^*, 
\no \\
&  \alpha_0=(I_m\; 0), \quad x\lessgtr x^\prime,\, 
z\in\bbC\backslash\bbR \lb{2.56}
\end{align}
Next, let $\phi\in L^2(\bbR)^{2m}$ and consider
\begin{equation}
J\psi^\prime(z,x)=(zI_{2m}+B(x))\psi(z,x)+ \phi(x), 
\quad z\in\bbC\backslash\bbR \lb{2.58}
\end{equation}
for a.e.\ $x\in\bbR$.  Then, as inferred from \cite{HS81}, 
\cite{HS83}, \eqref{2.58} has a 
unique solution $\psi(z,\dott)\in
L^2(\bbR)^{2m}\cap\AC_{\loc}(\bbR)^{2m}$ given by 
\begin{equation}
\psi(z,x)=\int_\bbR dx^\prime\, G(z,x,x^\prime) 
\phi(x^\prime). \lb{2.59}
\end{equation}
The Dirac-type operator $D$ in $L^2(\bbR)^{2m}$ associated with the
Hamiltonian system \eqref{HS} and \eqref{DS} is then defined by 
\begin{equation}
((D-z)^{-1}\psi)(x)= \int_\bbR dx^\prime\, 
G(z,x,x^\prime)\psi(x^\prime), \quad  \psi\in
L^2(\bbR)^{2m}, \; z\in\bbC\backslash\bbR. \lb{2.60} 
\end{equation}
Explicitly, one obtains
\begin{align}
D&=J \f{d}{dx}-B, \lb{2.61} \\
\dom(D)&=\{\phi\in L^2(\bbR)^{2m}\mid \phi
\in\AC_{\loc}(\bbR)^{2m}; \,(J\phi^\prime-B\phi)\in 
L^2(\bbR)^{2m} \}, \no
\end{align}
taking into account the limit point property of Dirac-type
systems as described in Lemma~\ref{l2.15}. Thus, $D$ is a 
self-adjoint operator in $L^2(\bbR)^{2m}$.

As described in \cite{HS81}--\cite{HS86}, the $2m\times 2m$
Weyl-Titchmarsh matrix $M(z,x_0,\alpha_0)$ associated with $D$ is then
defined by
\begin{align}
M(z,x_0,\alpha_0)
&=\big(M_{j,j^\prime}(z,x_0,\alpha_0)\big)_{j,j^\prime=1,2} \no \\
&=[G(z,x_0,x_0+0)+G(z,x_0,x_0-0)]/2,
\quad  z\in\bbC\backslash\bbR.  \lb{2.62} 
\end{align}
Actually, one can replace $\alpha_0=(I_m\; 0)$ by an arbitrary matrix
$\alpha=[\alpha_1\ \alpha_2]\in\bbC^{m\times 2m}$ satisfying
\eqref{2.8e} and hence introduces
\begin{subequations}\lb{2.620}
\begin{align}
M(z,x_0,\alpha)
&=\big(M_{j,j^\prime}(z,x_0,\alpha)\big)_{j,j^\prime=1,2},  \quad 
z\in\bbC\backslash\bbR,  \lb{2.62A} \\
M_{1,1}(z,x_0,\alpha)&=[M_-(z,x_0,\alpha)-M_+(z,x_0,\alpha)]^{-1},
\lb{2.62B} \\ 
M_{1,2}(z,x_0,\alpha)&=2^{-1}
[M_-(z,x_0,\alpha)-M_+(z,x_0,\alpha)]^{-1}
[M_-(z,x_0,\alpha)+M_+(z,x_0,\alpha)], \no \\
M_{2,1}(z,x_0,\alpha)&=2^{-1} [M_-(z,x_0,\alpha)+M_+(z,x_0,\alpha)]
[M_-(z,x_0,\alpha)-M_+(z,x_0,\alpha)]^{-1},\no \\
M_{2,2}(z,x_0,\alpha)&=M_\pm(z,x_0,\alpha)
[M_-(z,x_0,\alpha)-M_+(z,x_0,\alpha)]^{-1}M_\mp(z,x_0,\alpha). \no 
\end{align}
\end{subequations}
\
The basic results on $M(z,x_0,\alpha)$ then read as follows.

\begin{theorem} [\cite{GT97}, \cite{HS81}, \cite{HS82}, 
\cite{HS86}, \cite{KS88}] \lb{thm2.19} 
Assume Hypothesis~\ref{h2.4} and suppose \, that $z\in\bbC 
\backslash \bbR$, $x_0\in\bbR$, and that $\alpha\in\bbC^{m\times 2m}$
satisfies \eqref{2.8e}.  Then, \\
$(i)$ $M(z,x_0,\alpha)$ is a matrix-valued Herglotz function of rank 
$2m$ with representation
\begin{align}
&M(z,x_0,\alpha)=F(x_0,\alpha)+\int_\bbR d\Omega(\lambda,x_0,\alpha)\,
\big((\lambda-z)^{-1}-\lambda(1+\lambda^2)^{-1}\big), \lb{2.64} \\
&=\exp\bigg(C(x_0,\alpha)+\int_\bbR d\lambda \, \Upsilon
(\lambda,x_0,\alpha)
\big((\lambda-z)^{-1}-\lambda(1+\lambda^2)^{-1}\big) \bigg), \lb{2.65}
\end{align}
where
\begin{align}
F(x_0,\alpha)&=F(x_0,\alpha)^*, \quad \int_\bbR \Vert 
d\Omega(\lambda,x_0,\alpha)
\Vert_{\bbC^{2m\times 2m}} \,(1+\lambda^2)^{-1}<\infty, \lb{2.66} \\
C(x_0,\alpha)&=C(x_0,\alpha)^*, \quad 0\le\Upsilon(\dott,x_0,\alpha)
\le I_{2m} \, \text{  a.e.} \lb{2.67}
\end{align}
Moreover,
\begin{align}
\Omega((\lambda,\mu],x_0,\alpha)&=\lim_{\delta\downarrow
0}\lim_{\varepsilon\downarrow 0}\f1\pi
\int_{\lambda+\delta}^{\mu+\delta} d\nu \, 
\Im(M(\nu+i\varepsilon,x_0,\alpha)), \lb{2.68} \\
\Upsilon(\lambda,x_0,\alpha)&=\lim_{\varepsilon\downarrow 0}
\pi^{-1}\Im(\ln(M(\lambda+i\varepsilon,x_0,\alpha))) \text{ for a.e.\
$\lambda\in\bbR$}. \lb{2.69}
\end{align}
$(ii)$ $z\in\bbC\backslash\spec(D)$ if and only if $M(z,x_0,\alpha)$ is
holomorphic  near $z$. 
\end{theorem}

Here $\spec (T)$ abbreviates the spectrum of a linear operator $T$.

Next, we explicitly discuss the elementary Dirac-type example where
$A=I_{2m}$ and $B=0$.

\begin{example}\lb{e2.20} 
Suppose $A=I_{2m}$, $B=0$ and let $\alpha\in\bbC^{m\times 2m}$ satisfy
\eqref{2.8e}. Then,
\begin{align}
\Theta(z,x,x_0,\alpha)&=\begin{pmatrix}\theta_{1}(z,x,x_0,\alpha) \\
\theta_{2}(z,x,x_0,\alpha) \end{pmatrix}=\begin{pmatrix}
\alpha_1^*\cos(z(x-x_0))+\alpha_2^*\sin(z(x-x_0)) \\
\alpha_2^*\cos(z(x-x_0))-\alpha_1^*\sin(z(x-x_0)) \end{pmatrix}, \no \\
& \hspace*{7.5cm} \quad z\in\bbC, \lb{2.80} \\
\Phi(z,x,x_0,\alpha)&=\begin{pmatrix}\phi_{1}(z,x,x_0,\alpha) \\
\phi_{2}(z,x,x_0,\alpha) \end{pmatrix}=\begin{pmatrix}
-\alpha_2^*\cos(z(x-x_0))+\alpha_1^*\sin(z(x-x_0)) \\
\alpha_1^*\cos(z(x-x_0))+\alpha_2^*\sin(z(x-x_0)) \end{pmatrix}, \no \\
& \hspace*{7.5cm} \quad z\in\bbC, \lb{2.81} \\  
U_\pm (z,x,x_0,\alpha)&=\begin{pmatrix} u_{\pm,1}(z,x,x_0,\alpha) \\
u_{\pm,2}(z,x,x_0,\alpha) \end{pmatrix} 
=\begin{pmatrix} \alpha_1^* \mp i\alpha_2^* \\
\pm i(\alpha_1^* \mp i\alpha_2^*) \end{pmatrix}\exp(\pm iz(x-x_0)), 
\no \\
& \hspace*{7cm} \quad z\in\bbC_+, \lb{2.82} \\ 
M_\pm (z,x,\alpha)&=\pm iI_m, \quad z\in\bbC_+. \lb{2.83} 
\end{align}
\end{example}

Compared to the case of Schr\"odinger operators, it is remarkable that 
$M_\pm(z,x,\alpha)$ in \eqref{2.83} is, in fact, independent of $\alpha$.
Put differently, in Dirac-type situations, $M_\pm(z,x,\alpha)$ may contain
no information on the  boundary condition indexed by
$\alpha\in\bbC^{m\times 2m}$.

In Sections~\ref{s4} and \ref{s5} we will also refer to half-line Dirac
operators $D_+(\alpha)$ in $L^2([x_0,\infty))^{2m}$ associated with a
self-adjoint boundary condition at $x_0$ indexed by
$\alpha\in\bbC^{m\times 2m}$ satisfying \eqref{2.8e},
and hence briefly introduce
\begin{align}
D_+(\alpha)&=J \f{d}{dx}-B, \lb{2.84} \\
\dom(D_+(\alpha))&=\{\phi\in L^2([x_0,\infty))^{2m} \mid \phi
\in\AC([x_0,R])^{2m} \text{ for all $R>0$}; \no \\
& \hspace*{2.4cm} \alpha\phi(x_0)=0; \, (J\phi^\prime-B\phi)\in 
L^2([x_0,\infty))^{2m} \}, \no
\end{align}
taking into account the limit point property of Dirac-type
systems at $+\infty$ as described in Lemma~\ref{l2.15}. Thus,
$D_+(\alpha)$ is a  self-adjoint operator in $L^2([x_0,\infty))^{2m}$. In
complete analogy one introduces $D_-(\alpha)$ in $L^2((-\infty,
x_0])^{2m}$.

Next, we recall a few formulas in connection with Lagrange's identity 
needed in the proof of Theorem~\ref{t4.10} assuming 
$\alpha\in\bbC^{m\times 2m}$ satisfies \eqref{2.8e}. Then,  
explicitly, \eqref{2.310} and \eqref{2.330} read
\begin{align}
\theta_2(\bar z,x,x_0,\alpha)^*\theta_1(z,x,x_0,\alpha)-
\theta_1(\bar z,x,x_0,\alpha)^*\theta_2(z,x,x_0,\alpha)&=0, \lb{2.72}
\\
\phi_2(\bar z,x,x_0,\alpha)^*\phi_1(z,x,x_0,\alpha)-
\phi_1(\bar z,x,x_0,\alpha)^*\phi_2(z,x,x_0,\alpha)&=0, \lb{2.73} \\
\phi_2(\bar z,x,x_0,\alpha)^*\theta_1(z,x,x_0,\alpha)-
\phi_1(\bar z,x,x_0,\alpha)^*\theta_2(z,x,x_0,\alpha)&=I_m, \lb{2.74}
\\
\theta_1(\bar z,x,x_0,\alpha)^*\phi_2(z,x,x_0,\alpha)-
\theta_2(\bar z,x,x_0,\alpha)^*\phi_1(z,x,x_0,\alpha)&=I_m, \lb{2.75}
\end{align}
and 
\begin{align}
\phi_1(z,x,x_0,\alpha)\theta_1(\bar z,x,x_0,\alpha)^*-
\theta_1(z,x,x_0,\alpha_0)\phi_1(\bar z,x,x_0,\alpha)^*&=0, \lb{2.92}
\\
\phi_2(z,x,x_0,\alpha)\theta_2(\bar z,x,x_0,\alpha)^*-
\theta_2(z,x,x_0,\alpha)\phi_2(\bar z,x,x_0,\alpha)^*&=0, \lb{2.93} \\
\phi_2(z,x,x_0,\alpha)\theta_1(\bar z,x,x_0,\alpha)^*-
\theta_2(z,x,x_0,\alpha)\phi_1(\bar z,x,x_0,\alpha)^*&=I_m, \lb{2.94}
\\
\theta_1(z,x,x_0,\alpha)\phi_2(\bar z,x,x_0,\alpha)^*-
\phi_1(z,x,x_0,\alpha)\theta_2(\bar z,x,x_0,\alpha)^*&=I_m. \lb{2.95}
\end{align}
Finally, we note the connection between $\Phi$ defined in \eqref{FSb},
for different boundary value data
$\alpha, \gamma\in\bbC^{m\times 2m}$ satisfying \eqref{2.8e}, namely
\begin{equation}\lb{2.96}
\Phi(z,x,x_0,\gamma)=\Phi(z,x,x_0,\alpha)\alpha\gamma^* + 
\Theta(z,x,x_0,\alpha)\alpha J \gamma^*.
\end{equation}
This connection formula follows by the uniqueness of solutions of 
\eqref{HS} and by the identity given in \eqref{2.8i}. It  
is  needed in the proof of Theorem~\ref{t4.10}.

\section{The Leading Order Term in the Asymptotic \\ Expansion 
of $M_\pm (z,x,\alpha)$} \lb{s3}

Assuming Hypothesis~\ref{h2.4}, the principal result proven in this
section will be the following  leading-order asymptotic result for
half-line Weyl-Titchmarsh matrices $M_\pm(z,x_0,\alpha_0)$ associated
with the Dirac-type operator \eqref{2.61},  
\begin{equation}\lb{3.1}
M_\pm(z,x_0,\alpha_0) \underset{\substack{\abs{z}\to\infty\\ z\in
C_\varepsilon}}{=} \pm iI_{m} +\oh(1). 
\end{equation}
Here $\alpha_0 =(I_m\; 0)\in\bbC^{m\times 2m}$, and 
$C_{\varepsilon} \subset \bbC_+$ denotes the open sector with vertex at
zero, symmetry axis along the positive imaginary axis, and
opening angle $\varepsilon$, with $0<\varepsilon <\pi/2$.  

This particular topic originates with the order result of 
Hille~\cite{Hil63}
and the asymptotic formulas of Everitt~\cite{Ev72} and of Everitt and 
Halvorsen~\cite{EH78}. By appealing to the theory of Riccati equations, 
Atkinson in \cite{At81}, \cite{At82}, and \cite{At88a} obtains results
like those of Hille,  Everitt, and Halvorsen, both for the Schr\"odinger
case as well  as for the  scalar-Dirac ($m=1$) case. Through a deeper
understanding of  the role played by  Riccati theory, Atkinson obtains
the  first order asymptotic expansion of  $M_+(z,x,\alpha_0)$ for the
matrix-valued Schr\"odinger case in an unpublished manuscript \cite
{At88}.  Our strategy of proof for  \eqref{3.1} is patterned after
Atkinson's  approach which also appears in our recent work on the full
asymptotic expansion for $M_+(z,x,\alpha_0)$ in the matrix-valued
Schr\"odinger case \cite{CG99}.

We begin our discussion by noting two additional characterizations
for the Weyl disk, $\calD(z,c,x_0,\alpha)$, for the general Hamiltonian system \eqref{HSa}.
\begin{lemma}\lb{l3.3} 
Assume Hypotheses~\ref{h2.1} and \ref{h2.2}. 
Let $z\in\bbC\backslash\bbR$, $c\ne x_0$, and define
$U(z,x,x_0,\alpha)$, in terms of $M\in\bbC^{m\times m}$ by
\eqref{2.14}. Then $M\in \calD(z,c,x_0,\alpha)$ if and only if
\begin{equation}\lb{3.3}
\sigma(c,x_0,z)\Im (u_1(z,x,x_0,\alpha)^*u_2(z,x,x_0,\alpha)) > 0,
\quad x \in[x_0,c),
\end{equation}
or equivalently, if and only if
\begin{equation}\lb{3.4}
\sigma(c,x_0,z)\Im (u_2(z,x,x_0,\alpha)u_1(z,x,x_0,\alpha)^{-1}) > 0,
\quad x \in[x_0,c).
\end{equation}
Moreover, $M\in \calD_{\pm}(z,x_0,\alpha)$ if and only if \eqref{3.3} and
\eqref{3.4} hold for $c=\pm\infty$.
\end{lemma}
\begin{proof}
Let $U(z,x)=U(z,x,x_0,\alpha)$, and let $u_j(z,x)=u_j(z,x,x_0,\alpha)$, $j=1,2$
with $x\in[x_0,c)$. By \eqref{2.280},
\begin{align}\lb{3.6}
&2\sigma(c,x_0)|\Im(z)|\int_x^c ds\, U(z,s)^*A(s)U(z,s)\no \\
&=\sigma(x_0,c,z)U(z,s)^*(iJ)U(z,s)\Big |_x^c.
\end{align}
By \eqref{2.380}, this yields
\begin{align}\lb{3.5}
&2\sigma(c,x_0,z)\Im (u_1(z,x)^*u_2(z,x))\no\\
&= -E_c(M) + 2\sigma(c,x_0)|\Im(z)|\int_x^c ds\, U(z,s)^*A(s)U(z,s).
\end{align}
The integral expression in \eqref{3.5} is strictly positive  by 
Hypothesis~\ref{h2.2}. This yields the equivalence of $-E_c(M)\ge 0$, and hence
of  $M\in \calD(z,c,x_0,\alpha)$,  with the condition given in \eqref{3.3}.
The equivalence of \eqref{3.3} and \eqref{3.4} follows
from the observation that
\begin{equation}
\Im (u_2(z,x)u_1(z,x)^{-1}) = (u_1(z,x)^{-1})^*\Im (u_1(z,x)^*u_2(z,x))u_1(z,x)^{-1}.
\end{equation}
The analogous characterization of $\calD_{\pm}(z,x_0,\alpha)$ now follows
from Definition~\ref{dLWD}.
\end{proof}

In Lemma~\ref{l3.3},  $u_j(z,c)$, $j=1,2$, are well-defined and 
$E_c(M)=0$ precisely when $\sigma(c,x_0,z)\Im (u_1(z,c)^*u_2(z,c))= 0$.
A similar statement might not hold for \eqref{3.4} since
$u_1(z,c,x_0,\alpha)$ might be singular.  In part, the
latter point motivates the next characterization of the disk.
\begin{lemma}\lb{l3.4} 
Assume Hypotheses~\ref{h2.1} and \ref{h2.2}.
Let $z\in\bbC\backslash\bbR$, $c\ne x_0$, and define 
$u_j(z,x)=u_j(z,x,x_0,\alpha)$, $j=1,2$,  by
\eqref{2.14}.  Then $M\in \calD(z,c,x_0,\alpha)$ if and only if
\begin{equation}\lb{3.8}
u_1(z,x) -i\sigma(c,x_0,z)u_2(z,x)
\end{equation}
is nonsingular for $x\in[x_0,c]$ and
\begin{equation}\lb{3.9}
\begin{split}
\vartheta(z,x)=\vartheta(z,x,x_0,\alpha)&= [u_1(z,x)
+i\sigma(c,x_0,z)u_2(z,x)]\times\\
& \quad \times [u_1(z,x)
-i\sigma(c,x_0,z)u_2(z,x)]^{-1}
\end{split}
\end{equation}
satisfies
\begin{equation}\lb{3.10}
 I_m-\vartheta(z,x)^*\vartheta(z,x)>0,\quad
x\in[x_0,c),
\end{equation}
with nonnegativity holding at $x=c$.
Moreover, $M\in \calD_{\pm}(z,x_0,\alpha)$ if and only if \eqref{3.9}
is well-defined on $[x_0,\pm\infty)$ and
\eqref{3.10} holds for $c=\pm\infty$.
\end{lemma}
\begin{proof}
Let $M\in\calD(z,c,x_0,\alpha)$ and suppose that
$u_1(z,\xi)v=i\sigma(c,x_0,z)u_2(z,\xi)v$ for $\xi\in[x_0,c]$ and
$v\in\bbC^m$, $v\ne 0$. Then,
\begin{equation}
v^*\sigma(c,x_0,z)\Im (u_1(z,\xi)^*u_2(z,\xi))v =
-v^*u_1(z,\xi)^*u_1(z,\xi)v.
\end{equation}
By \eqref{3.3}, an immediate contradiction results if $\xi\ne c$. However,
if $\xi=c$, then either $v^*E_c(M)v>0$  or $u_j(z,c)v=0$, $j=1,2$. In either
case, a contradiction results since $E_c(M)\leq 0$ by Definition~\ref{dWD}
and $U=(u_1^t,u_2^t)^t$ satisfies the first-order system \eqref{HSa}.
Hence, $\vartheta(z,x)$ is well-defined on $[x_0,c]$. For $x\in [x_0,c)$
and $M\in \calD(z,c,x_0,\alpha)$, \eqref{3.3} implies that
\begin{equation}\lb{3.12}
2i\sigma(c,x_0,z)(u_1(z,x)^*u_2(z,x) - u_2(z,x)^*u_1(z,x) )< 0.
\end{equation}
This is equivalent to
\begin{equation}\lb{3.13}
\begin{split}
&[u_1(z,x)^*-i\sigma(c,x_0,z)u_2(z,x)^*]
[u_1(z,x)+i\sigma(c,x_0,z)u_2(z,x)] \\
&<[u_1(z,x)^*+i\sigma(c,x_0,z)u_2(z,x)^*][u_1(z,x)
-i\sigma(c,x_0,z)u_2(z,x)]
\end{split}
\end{equation}
 on $[x_0,c)$. Given the nonsingularity of
$u_1(z,x)-i\sigma(c,x_0,z)u_2(z,x)$ on
$[x_0,c]$,  \eqref{3.13} implies \eqref{3.10}, with nonnegativity holding
at $x=c$. \\
Next, let $M\in\bbC^{m\times m}$, and suppose that
$\vartheta(z,x)$, defined by \eqref{3.9}, is well-defined on $[x_0,c]$,
and satisfies \eqref{3.10}. Then, on $[x_0,c)$, 
\eqref{3.13} and consequently \eqref{3.12} follow, which implies
that \eqref{3.3} holds, and hence that $M\in\calD(z,c,x_0,\alpha)$.
The analogous characterization of $\calD_{\pm}(z,x_0,\alpha)$
follows from Definition~\ref{dLWD}.
\end{proof}

By Lemma~\ref{l3.3} one notes, for $z\in \bbC\backslash\bbR$, that $M\in
\calD(z,c,x_0,\alpha)$ if and only if
\begin{equation}\lb{3.14}
V(z,x,x_0,\alpha)= u_2(z,x,x_0,\alpha)u_1(z,x,x_0,\alpha)^{-1},
\quad x\in [x_0,c),
\end{equation}
is well-defined  while satisfying
\begin{equation}\lb{3.15}
\sigma(c,x_0,z)\Im (V(z,x,x_0,\alpha)) > 0, \quad x\in [x_0,c).
\end{equation}
In terms of $V(z,x,x_0,\alpha)$ and by \eqref{3.9}, one notes that
\begin{equation}\lb{3.16}
\begin{split}
\vartheta(z,x,x_0,\alpha) &=
[I_m + i\sigma(c,x_0,z)V(z,x,x_0,\alpha)]\times\\
& \quad \times[I_m - i\sigma(c,x_0,z)V(z,x,x_0,\alpha)]^{-1},\quad
x\in[x_0,c),
\end{split}
\end{equation}
is a Cayley-type  transformation of
$V(z,x,x_0,\alpha)$. In the scalar context, this
transformation corresponds to a conformal
mapping of the complex upper half-plane to the
unit disk.
Moreover, defined as it is, $V(z,x,x_0,\alpha)$
satisfies a Riccati differential equation
that is associated with the Hamiltonian system
\eqref{HSa} while $\vartheta(z,x,x_0,\alpha)$
satisfies a Riccati  equation obtained by the
Cayley-type transformation \eqref{3.16} applied to
the differential equation satisfied by $V(z,x,x_0,\alpha)$.

For the Dirac-type case of \eqref{HSa}, one observes  by a simple
calculation  that $V(z,x,x_0,\alpha_0)$ is seen to satisfy a particular
initial value problem for a Riccati differential equation.
\begin{lemma} \lb{l3.5} 
Assume Hypotheses~\ref{h2.1}, \ref{h2.2}, and \ref{h2.4}.
Let $\alpha_0=(I_m\; 0)\in\bbC^{m\times 2m}$, let $u_j(z,x)=u_j(z,x,x_0,\alpha_0)$,
$j=1,2$, be defined by \eqref{2.14},  and suppose that  $V(z,x,x_0,\alpha_0)$ 
is well-defined by \eqref{3.14}. Then, $V(z,\cdot)=V(z,\cdot,x_0,\alpha_0)$ satisfies,
\begin{subequations}\lb{3.17}
\begin{align}
&V'(z,x) +zV(z,x)^2 + V(z,x) B_{2,2}(x)V(z,x) +  B_{1,2}(x)V(z,x)
+ V(z,x)B_{2,1}(x) \no \\\
& + B_{1,1}(x) +zI_m =0, \lb{3.17a}\\
&V (z,x_0)=M, \lb{3.17b}
\end{align}
\end{subequations}
where $B_{j,k}\in \bbC^{m\times m}$, $j,k= 1,2$, are defined in
\eqref{2.1d}.
\end{lemma}

Hence,  by Lemma~\ref{l3.3}, the associated relations \eqref{3.14} and
\eqref{3.15}, and the uniqueness of solutions for \eqref{3.17}, we 
obtain the following result for the Dirac-type case.
\begin{theorem}\lb{t3.6} 
Assume Hypotheses~\ref{h2.1}, \ref{h2.2}, and \ref{h2.4}, and let 
$\alpha_0=(I_m\; 0)\in\bbC^{m\times 2m}$. Then,  $M\in
\calD(z,c,x_0,\alpha_0)$ if and only if the initial value problem
given by \eqref{3.17} has a  solution, $V(z,\cdot)$, well-defined and
satisfying
\begin{equation}\lb{3.18}
  \sigma(c,x_0,z)\Im (V(z,x)) > 0,\quad  x\in[x_0,c).
\end{equation}
Moreover, $M\in \calD_{\pm}(z,x_0,\alpha_0)$
if and only if \eqref{3.18} holds for $c=\pm\infty$.
\end{theorem}

\begin{remark} \lb{r3.6a}
An important consequence of Theorem~\ref{t3.6} and
the uniqueness of
solutions for \eqref{3.17} is that solution
trajectories for \eqref{3.17}, which satisfy
\eqref{3.18}, consist of elements of Weyl disks;
that is,
\begin{equation}\lb{3.19}
V(z,x,x_0,\alpha_0)\in \calD(z,c,x,\alpha_0),
\quad  x\in [x_0,c).
\end{equation}
Given the characterization of $\calD(z,c,x_0,\alpha_0)$ in
Defintion~2.7A, for each $x\in [x_0,c)$ there is
a $\beta\in\bbC^{m\times 2m}$ with  
$\sigma(c,x_0,z)\Im (\beta_2\beta_1^*)\ge 0$, such that
\begin{equation}
V(z,x,x_0,\alpha_0)=M(z,c,x,\alpha_0,\beta).
\end{equation}
It is in this sense that we let $M(z,c,x,\alpha_0)$
denote our solution of the initial value problem \eqref{3.17}
that satisfies \eqref{3.18}.  Analogously,
\begin{equation}
V(z,x,x_0,\alpha_0)\in \calD_{\pm}(z,x,\alpha_0),\quad  x\in
[x_0,\pm\infty),
\end{equation}
for trajectories of \eqref{3.17} that satisfy
\eqref{3.18} for $c=\pm\infty$.
Hence, in this sense, we let
$M_{\pm}(z,x,\alpha_0)$ denote those
solutions of \eqref{3.17} that satisfy
\eqref{3.18} for $c=\pm\infty$.
However, by  Lemma~\ref{l2.15}, our Dirac system
is in the limit point case at $\pm\infty$.  Each
$\calD_{\pm}(z,x,\alpha_0)$ consists of a
unique matrix, and thus
$M_{\pm}(z,x,\alpha_0)$ describes {\em
unique} trajectories for \eqref{3.17a}.  This
contrasts with the matrix-valued Schr\"odinger case
considered in \cite{CG99} where there are as many
trajectories, each denoted by either
$M_{+}(z,x,\alpha_0)$ or
$M_{-}(z,x,\alpha_0)$, as there are
matrices in a given initial disk
$\calD_{\pm}(z,x_0,\alpha_0)$.
\end{remark}

Now for the Dirac-type case \eqref{DS}
with $\alpha_0=(I_m\; 0)\in \bbC^{m\times 2m}$, with
$\vartheta(z,x)=\vartheta(z,x,x_0,\alpha_0)$  defined in 
\eqref{3.9} and \eqref{3.16}, and with  $x\in[x_0,c)$, 
one concludes that
\begin{equation}\lb{3.22}
\vartheta(z,x)[u_1(z,x)-i\sigma(c,x_0,z)u_2(z,x)]=u_1(z,x)+
i\sigma(c,x_0,z)u_2(z,x),
\end{equation}
and hence that
\begin{subequations}\lb{3.23}
\begin{align}
I_m + \vartheta(z,x)&=
2u_1(z,x)[u_1(z,x)-i\sigma(c,x_0,z)u_2(z,x)]^{-1},\\ I_m -
\vartheta(z,x)&= -2i\sigma(c,x_0,z)u_2(z,x)[u_1(z,x)-i
\sigma(c,x_0,z)u_2(z,x)]^{-1}.
\end{align}
\end{subequations}
Differentiating \eqref{3.22} one obtains
\begin{equation}
\begin{split}
\vartheta'(u_1-i\sigma u_2)&=(I_m -\vartheta)
(zu_2 + B_{2,1}u_1 + B_{2,2}u_2) \\
& \quad +i\sigma (I_m +\vartheta)(-zu_1 - B_{1,1}u_1 -B_{1,2}u_2).
\end{split}
\end{equation}
By \eqref{3.23} one concludes that $\vartheta(z,\cdot,x_0,\alpha_0)$
satisfies the initial value problem given by
\begin{subequations}\lb{3.25}
\begin{align}
\vartheta'(z,x)&= \frac{1}{2} \begin{pmatrix}I_m +
\vartheta(z,x)^t\\ I_m - \vartheta(z,x)^t\end{pmatrix}^t\times  \no \\
& \quad \times \begin{pmatrix} -i\sigma(c,x_0,z)(zI_m +B_{1,1}(x)) &
B_{1,2}(x)\\B_{2,1}(x)&  i\sigma(c,x_0,z)(zI_m +B_{2,2}(x))
\end{pmatrix}\times\notag \\ 
& \quad \times\begin{pmatrix}I_m +
\vartheta(z,x)\\ I_m - \vartheta(z,x)\end{pmatrix}, \lb{3.25a} \\[5pt]
\vartheta(z,x_0)&=(I_m +i\sigma(c,x_0,z) M)
(I_m -i\sigma(c,x_0,z) M)^{-1},\lb{3.25b}
\end{align}
\end{subequations}
where $B_{j,k}\in \bbC^{m\times m}$,
$j,k= 1,2$, satisfy Hypothesis~\ref{h2.1}.

By Lemma~\ref{l3.4} and the uniqueness of solutions for \eqref{3.25},
one obtains the following result in the Dirac-type case \eqref{DS}.

\begin{theorem}\lb{t3.7}
Assume Hypothesis~\ref{h2.4}. Then 
$M\in \calD(z,c,x_0,\alpha_0)$ if and only if the initial value
problem given by \eqref{3.25} has a  solution, $\vartheta(z,\cdot)$
which is well-defined on $[x_0,c]$ and satisfies
\begin{equation}\lb{3.26}
I_m-\vartheta(z,x)^*\vartheta(z,x)> 0,\quad
x\in[x_0,c).
\end{equation}
Moreover, $M\in \calD_{\pm}(z,x_0,\alpha_0)$ if and only if
\eqref{3.26} holds for  $c=\pm\infty$.
\end{theorem}
Given the positivity present in \eqref{3.26}, we note the exact
correspondence which exists, by \eqref{3.16}, between solutions of
\eqref{3.17} that satisfy \eqref{3.18} and those solutions of
\eqref{3.25} that satisfy \eqref{3.26}. In particular, given  Remark~\ref{r3.6a}, 
we rewrite \eqref{3.16} as
\begin{equation}\lb{3.27}
\begin{split}
\vartheta(z,x,x_0,\alpha_0) &=
[I_m + i\sigma(c,x_0,z)M(z,c,x,\alpha_0)]\times\\
  & \quad \times[I_m - i\sigma(c,x_0,z)M(z,c,x,\alpha_0)]^{-1},
\quad x\in[x_0,c),
\end{split}
\end{equation}
Moreover, our Dirac system is in the limit point case at $\pm\infty$.
Consequently, there are unique solutions
of \eqref{3.25}, $\vartheta_{\pm}(z,\cdot,x_0,\alpha_0)$,
$z\in\bbC\backslash\bbR$,  which satisfy \eqref{3.26} for
$c=\pm\infty$, and  which correspond to the unique solutions of
\eqref{3.17}, $M_{\pm}(z,x,\alpha_0)$,  which satisfy \eqref{3.18}
for $c=\pm\infty$; specifically,
\begin{equation}\lb{3.28}
\vartheta_{\pm}(z,x,x_0,\alpha_0)= [I_m \pm
i\sigma(z)M_{\pm}(z,x,\alpha_0)][I_m \mp
i\sigma(z)M_{\pm}(z,x,\alpha_0)]^{-1}.
\end{equation}

These relationships  form the basis for the
analysis to follow. The asymptotic result \eqref{3.1} is obtained 
by an analysis of the corresponding asymptotic behavior for all
solutions $\vartheta(z,\cdot,x_0,\alpha_0)$
described in  \eqref{3.25}, these include among
them the particular solutions
$\vartheta_{\pm}(z,\cdot,x_0,\alpha_0)$. Thus
asymptotic behavior is deduced for all
corresponding solutions
$M(z,c,\cdot,\alpha_0)$ of \eqref{3.17} which
include among them the solutions
$M_{\pm}(z,\cdot,\alpha_0)$.  The advantage of
this approach comes from the compactification
inherent in the Cayley-type transformation
\eqref{3.27}, and the resulting boundedness of the
solutions as a consequence of \eqref{3.26}.

We pause for a moment to address, in the following remark,
a point raised by us in \cite{CG99} for the matrix-valued
Schr\"odinger case described in \eqref{SS}.
\begin{remark}\lb{r3.3}
With $u_j(z,x)=u_j(z,x,x_0,\alpha)$, $j=1,2$, defined in \eqref{2.14} for
the general Hamiltonian system \eqref{HSa}, an analog to Lemma~\ref{l3.4}
for the characterization of $\calD(z,c,x_0,\alpha)$ is obtained by replacing
the expression in \eqref{3.8} with
\begin{equation}
u_1(z,x) -i|z|^{-1/2}\sigma(c,x_0,z)u_2(z,x),
\end{equation}
and by replacing the definition for $\vartheta(z,x)=\vartheta(z,x,x_0,\alpha)$
given in \eqref{3.9} with
\begin{equation}
\begin{split}
\vartheta(z,x)  =& (u_1(z,x) +i|z|^{-1/2}\sigma(c,x_0,z)u_2(z,x))\times\\
& \times (u_1(z,x) -i|z|^{-1/2}\sigma(c,x_0,z)u_2(z,x))^{-1}.
\end{split}
\end{equation}
Specific to the matrix-valued Schr\"odinger case, we obtain analogs of 
Lemma~\ref{l3.5}, Theorem~\ref{t3.6}, and Theorem~\ref{t3.7} by replacing
equation \eqref{3.17a} with
\begin{equation}
V'(z,x) + V(z,x)^2 - Q(x) +zI_m = 0
\end{equation}
and by replacing the equations in \eqref{3.25} with
\begin{subequations}\lb{3.280}
\begin{align}
\vartheta'(z,x)&= \sigma(c,x_0,z)\frac{1}{2} 
\begin{pmatrix}I_m +\vartheta(z,x)^t\\ 
I_m - \vartheta(z,x)^t
\end{pmatrix}^t 
\begin{pmatrix} -i|z|^{-1/2}(zI_m - Q(x)) & 0 \\
0 &  i|z|^{-1/2}I_m 
\end{pmatrix}\times\notag \\ 
& \quad \times\begin{pmatrix}I_m + \vartheta(z,x)\\ 
I_m - \vartheta(z,x)\end{pmatrix},\lb{3.280a} \\[5pt]
\vartheta(z,x_0)&=(I_m + i|z|^{-1/2}\sigma(c,x_0,z)M)
(I_m - i|z|^{-1/2}\sigma(c,x_0,z)M)^{-1}.\lb{3.280b}
\end{align}
\end{subequations}
$\calD^{\calR}(z,c,x_0,\alpha_0)$ was defined in \cite{CG99} to be the set
of those $M\in \bbC^{m\times m}$ for which the intial value problem given
by \eqref{3.280} has a solution, $\vartheta(z,x)$, which is well-defined 
on $[x_0,c]$ and satisfies \eqref{3.26}.  In \cite{CG99} we  showed that 
$\calD(z,c,x_0,\alpha_0)\subseteq\calD^{\calR}(z,c,x_0,\alpha_0)$. 
This was sufficient for the subsequent analysis
in \cite{CG99}. However, as the analog of Theorem~\ref{t3.7} now shows,
one actually has equality of the two disks in \cite{CG99}, that is, 
\begin{equation}
\calD(z,c,x_0,\alpha_0)=\calD^{\calR}(z,c,x_0,\alpha_0).
\end{equation}
\end{remark}

To obtain a proof of \eqref{3.1} for the Dirac-type case, we adapt an
approach due to Atkinson \cite{At88} for proving a result
analogous to \eqref{3.1} for the matrix-valued Schr\"odinger case 
(cf., e.g., \cite[Theorem 3.1]{CG99}) In light of Remark~\ref{r3.2}, we
begin by restricting our attention to  $z\in\bbC_+$, and as in the
previous discussion, take $\alpha_0=(I_m\; 0)\in
\bbC^{m\times 2m}$.

First we introduce two systems related to \eqref{3.25} by means of a
change of variables. Let
\begin{equation}\lb{3.29}
\varphi(z,t)=\vartheta(z,x),\qquad
  t= (x-x_0)|z|,\qquad  x \in [ x_0  ,  c ).
\end{equation}
With this change,  \eqref{3.25} becomes
\begin{subequations}\lb{3.30}
\begin{align}\lb{3.30a}
\varphi'(z,t)&= \frac{1}{2} |z|^{-1}\begin{pmatrix}I_m
+ \varphi(z,t)^t\\ I_m - \varphi(z,t)^t\end{pmatrix}^t
\begin{pmatrix}\mp i(zI_m +\wti B_{1,1}(t)) & \wti
B_{1,2}(t)\\[1mm] \wti B_{2,1}(t)& \pm i(zI_m +\wti B_{2,2}(t))
\end{pmatrix}\times \notag \\ 
& \quad \times\begin{pmatrix}I_m +
\varphi(z,t)\\ I_m - \varphi(z,t)\end{pmatrix}.
\end{align}
\noindent With $M=M(z,c,x_0,\alpha_0)\in
\calD(z,c,x_0,\alpha_0)$
  \eqref{3.25b} becomes
\begin{equation}\lb{3.30b}
\varphi(z,0)=(iI_m \mp M(z,c,x_0,\alpha_0))
(iI_m \pm M(z,c,x_0,\alpha_0))^{-1},
\end{equation}
\noindent and \eqref{3.26} becomes
\begin{equation}\lb{3.30c}
\varphi(z,t)^*\varphi(z,t)< I_m
\qquad t\in[ 0, (c-x_0)|z|),
\end{equation}
where in \eqref{3.30a},
\begin{equation}\lb{3.30d}
\wti B_{j,k}(t)= B_{j,k}(x_0 +t|z|^{-1}), \quad j,k=1,2.
\end{equation}
\end{subequations}
In the complete system \eqref{3.30}, one now has a set of
conditions equivalent to system \eqref{3.25} and \eqref{3.26}.

We recall that $C_\varepsilon \subset \bbC_+$ represents
the open sector with vertex at zero, symmetry axis along the
positive imaginary axis, and opening angle $\varepsilon$,
with $0<\varepsilon <\pi/2$. Next, consider a sequence, $z_n \in
\bbC_{\varepsilon}$,  $n\in\bbN $, such that $|z_n| \to
\infty$ as $n\to \infty$ and such that
\begin{equation}\lb{3.31}
0< \varepsilon <  \delta_n = \arg{(z_n)}
< \pi - \varepsilon.
\end{equation}
By choosing an appropriate subsequence, we may assume that
\begin{equation}\lb{3.32}
\delta_n \to \delta \in [\varepsilon, \pi - \varepsilon].
\end{equation}
Let $\varphi (z_n ,t)$ denote a corresponding
sequence of functions that satisfy \eqref{3.30a}
and \eqref{3.30c}, with initial data, $\varphi
(z_n ,0)$, defined by \eqref{3.30b} for a sequence
of points  $M(z_n,c,x_0,\alpha_0)$, where each
$M(z_n,c,x_0,\alpha_0)$ is chosen to be an element of the disk
$\calD(z_n,c,x_0,\alpha_0)$. Note that as $z_n\to
\infty$, the intervals described in \eqref{3.30c}
eventually cover all compact subintervals of $
\bbR_+$. Given the uniform boundedness of
$\varphi_n(t)=\varphi (z_n ,t)$ described in
\eqref{3.30c}, we assume, upon passing to an
appropriate subsequence still denoted by
$\varphi_n (0)$, that
\begin{equation}\lb{3.33}
\varphi_n(0) = \varphi (z_n,0) \rightarrow
\varphi_\pm (\delta), \ \text{for}
\ \pm (c-x_0) > 0 \
\text{ as } n\rightarrow \infty ,
\end{equation}
and as a consequence, that
\begin{equation}\lb{3.34}
{\varphi_\pm(\delta)}^* \varphi_\pm(\delta) \le I_m.
\end{equation}

With $\varphi_\pm(\delta)$ defined in \eqref{3.33} as
$|z_n|\to\infty$, we
consider limiting systems associated  with \eqref{3.30}:
\begin{subequations}\lb{3.35}
\begin{align}
\eta_\pm '(t)&= \frac{1}{2}
\begin{pmatrix}
I_m+ \eta_\pm (t)^t \\ I_m-\eta_\pm (t)^t
\end{pmatrix} ^t
\begin{pmatrix}
\mp ie^{i\delta}I_m & 0\\ 0 & \pm ie^{i\delta}I_m
\end{pmatrix}\begin{pmatrix}
I_m+\eta_\pm (t)\\ I_m- \eta_\pm (t)
\end{pmatrix}, \quad \pm t\ge 0, \lb{3.35a} \\
\eta_\pm (0)&= \varphi_\pm(\delta). \lb{3.35b}
\end{align}
\end{subequations}
\begin{theorem}\lb{t3.8} Assume Hypothesis~\ref{h2.4}. Then 
the solution $\eta_\pm$ of \eqref{3.35} satisfies
\begin{equation}\lb{3.36}
\eta_\pm (t)^* \eta_\pm(t) \le I_m,\quad
t\in [0, \pm\infty).
\end{equation}
Moreover, the
solutions $\varphi_n =\varphi (z_n,\cdot)$
of \eqref{3.30} converge to $\eta_\pm$  uniformly on
$[0,\pm T]$ for every $T>0$, as $n\to \infty $.
\end{theorem}
\begin{proof}
In this proof, we consider only the case corresponding
to $t\ge 0$, that is, $\eta_+(0)=\varphi_+(\delta)$ in
\eqref{3.35b}.  The other case follows in a similar manner.
For this reason, we let $\eta( t)= \eta_+( t)$ in the
remaining discussion. We also let $T\in \bbR_+$ be the
greatest value such that
\eqref{3.36} holds for $t\in [0, T] $ and show that
\eqref{3.36} must hold for some $[0, T'] $ with  $ T' > T $,
thus proving $T=\infty.$\\
The solution of  \eqref{3.35}, $\eta$,
presumed to be defined on  $[0, T]$, can be continued
onto some $[0, T']$ with $T' > T$; $\eta$ then satisfies
\begin{equation}\lb{3.37}
\eta (t)^* \eta (t) \le \kappa^2 I_m
\end{equation}
for $0\le t \le T'$ and for some $\kappa\ge 1$. \\
For brevity, let $\varphi'_n (t)= G_n(\varphi_n,t) $ denote
\eqref{3.30a} with
$z=z_n$, and let $\eta'(t)= H(\eta,t) $ denote \eqref{3.35a}
in the following. Integrating \eqref{3.35a} and
\eqref{3.30a}, one obtains
\begin{align}\lb{3.38}
\varphi_n (t) -\eta (t) &= \varphi_n(0) -\varphi_0(\delta)
+ \int_0^t  ds \{
G_n(\eta,s) - H(\eta,s)\}  \no \\
&\quad + \int_0^t ds \{ G_n(\varphi_n,s) - G_n(\eta,s)\}.
\end{align}
We note that
\begin{align}\lb{3.39}
G_n(\eta ,s) - H(\eta ,s) &=
\frac{1}{2}i(e^{i\delta} - e^{i\delta_n})
(I_m +\eta(s) )^2
-\frac{1}{2}i(e^{i\delta} - e^{i\delta_n})
(I_m -\eta(s) )^2 +\no\\
& \quad + \sum_{j,k =1}^2 F_{j,k}(z_n,s),
\end{align}
where,
\begin{subequations}\lb{3.40}
\begin{align}
F_{1,1}(z_n,s)&=-\frac{1}{2}i|z_n|^{-1}
(I_m+\eta(s))\widetilde B_{1,1}(s)(I_m+\eta(s)),\\
F_{2,2}(z_n,s)&=\frac{1}{2}i|z_n|^{-1}
(I_m-\eta(s))\widetilde B_{2,2}(s)(I_m-\eta(s)),\\
F_{1,2}(z_n,s)&=\frac{1}{2}i|z_n|^{-1}
(I_m+\eta(s))\widetilde B_{1,2}(s)(I_m-\eta(s)),\\
F_{2,1}(z_n,s)&=\frac{1}{2}i|z_n|^{-1}
(I_m-\eta(s))\widetilde B_{2,1}(s)(I_m+\eta(s)).
\end{align}
\end{subequations}
Thus, for $t\in [0,T']$,  \eqref{3.37} implies that
as $n\to\infty$
\begin{equation}
|e^{i\delta}-e^{i\delta_n}|\int_0^t
\|I_m \pm \eta (s) \|_{\bbC^{m\times m}}^2 ds= \oh (1),
\end{equation}
and together with \eqref{3.29} and \eqref{3.30d} that
\begin{equation}\lb{3.42}
\int_0^t  \| F_{j,k}(s) \|_{\bbC^{m\times m}}ds =
O\bigg( \int_{x_0}^{x_0 + t|z_n|^{-1}} \| \widetilde
B_{j,k}(s) \|_{\bbC^{m\times m}} ds\bigg) =\oh(1).
\end{equation}
(Here $\|\cdot\|_{\bbC^{m\times m}}$ denotes a norm on $\bbC^{m\times
m}$.) Hence, by \eqref{3.39}--\eqref{3.42}, one infers that for $t\in
[0,T']$ and as $n\to\infty$,
\begin{equation}\lb{3.43}
\int_0^t \{ G_n(\eta,s) -H(\eta,s)\} ds = \oh(1).
\end{equation}
Next, one notes that
\begin{equation}\lb{3.44}
G_n(\varphi_n,s) - G_n(\eta,s) = 2ie^{i\delta_n}
(\eta(s) - \varphi_n(s)) + \sum_{j,k =1}^2 K_{j,k}(z_n,s),
\end{equation}
where
\begin{subequations}
\begin{align}
K_{1,1}(z_n,s) &= \frac{-i}{2}|z_n|^{-1}
\{ (I_m +\varphi_n) B_{1,1}(s) (\varphi_n -\eta)
+ (\varphi_n
-\eta)B_{1,1}(s) (I_m +\eta) \},\\ K_{2,2}(z_n,s) &=
\frac{i}{2}|z_n|^{-1}\{ (I_m -\varphi_n) B_{2,2}(s) (\eta
-\varphi_n) + (\eta -\varphi_n)B_{2,2}(s) (I_m -\eta) \},\\
K_{1,2}(z_n,s) &= \frac{1}{2}|z_n|^{-1}\{(I_m
+\varphi_n)B_{1,2}(s)(\eta -\varphi_n) + (\varphi_n
-\eta)B_{1,2}(s)(I_m -\eta) \},\\ K_{2,1}(z_n,s) &=
\frac{1}{2}|z_n|^{-1}\{(I_m -\varphi_n) B_{2,1}(s) (\varphi_n
-\eta) + (\eta -\varphi_n)B_{2,1}(s) (I_m +\eta) \}.
\end{align}
\end{subequations}
By \eqref{3.34}  and \eqref{3.37}, for $s\in [0,T']$,
\begin{equation}\lb{3.46}
\| I_m \pm \varphi_n(s) \|_{\bbC^{m\times m}}\le 2,
\qquad \| I_m \pm \eta(s)  \|_{\bbC^{m\times m}}\le \kappa +1,
\end{equation}
and hence by \eqref{3.44}--\eqref{3.46},
\begin{align}\lb{3.47}
&\| G_n(\varphi_n,s) - G_n(\eta,s) \|_{\bbC^{m\times m}} \no \\
&\le \|
\eta (s)- \varphi_n (s)\|_{\bbC^{m\times m}} \bigg\{ 2 +
\frac{|z_n|^{-1}}{2}(3+\kappa)\sum_{j,k=1}^2\| \widetilde
B_{j,k}(s) \|_{\bbC^{m\times m}} \bigg\}.
\end{align}
Of course, by \eqref{3.33} as $n\to \infty$,
\begin{equation}\lb{3.48}
\phi_n(0) - \phi_+(\delta) = \oh(1).
\end{equation}
Thus, by \eqref{3.42}, \eqref{3.47} and \eqref{3.48},
one concludes for $t\in [0,T']$ and as $n\to \infty$, that
\begin{align}\lb{3.49}
&\|\varphi_n(t)-\eta(t)  \|_{\bbC^{m\times m}} \le \oh(1) \no \\
&+ \int_0^t \|\varphi_n(s)-\eta(s)  \|_{\bbC^{m\times m}} \bigg\{ 2 +
\frac{|z_n|^{-1}}{2}(3+\kappa)\sum_{j,k=1}^2\| \widetilde
B_{j,k}(s) \|_{\bbC^{m\times m}} \bigg\}ds.
\end{align}

Gronwall's inequality applied to \eqref{3.49}
together with a consideration of the effect of the variable
change \eqref{3.29}, as illustrated in \eqref{3.42}, yields
\begin{equation}\lb{3.50}
\varphi_n(t)-\eta(t)\to 0 \,\text{ as } \, n\to \infty
\end{equation}
uniformly for $t\in [0,T']$.  Thus by \eqref{3.30c},
the contradiction results that for all $t\in [0,T']$,
$\eta$ satisfies \eqref{3.36}.
\end{proof}

What  solutions of \eqref{3.35}  satisfy \eqref{3.36}?
\begin{lemma} Assume Hypothesis~\ref{h2.4}. 
If $\eta_\pm$ is a  solution of \eqref{3.35a}
which satisfies \eqref{3.36}, then
\begin{equation}\lb{3.51}
0=\eta_\pm (t) ,\quad t\in [0, \pm\infty).
\end{equation}
\end{lemma}
\begin{proof}
We note that \eqref{3.35a} is equivalent to \eqref{3.30a}
with $\widetilde B=0$. By the variable change \eqref{3.29},
\eqref{3.35a} is also equivalent to \eqref{3.25a}  with
$B=0$. Next, we recall the connection between the Riccati-type
equations \eqref{3.25a}, and \eqref{3.17a}  by means of the
Cayley transformation \eqref{3.27}. Solution matrices of
\eqref{3.35a} which statisfy \eqref{3.36} at $t=0$ thus
correspond to solution matrices, $V(z,\cdot)$, of \eqref{3.17a} for
which $\Im (V(z,x_0))\ge 0$. Moreover, solutions of \eqref{3.17a}
for which for which $\Im (V(z,x_0))\ge 0$ are obtainable from
solutions of \eqref{HSa}, with $B=0$, by means of
\eqref{2.14} with $\Im (M)\ge 0$.  Thus, by utilizing this
connection between explicit exponential solutions of
\eqref{HSa} with $B=0$ and solutions of the Riccati-type
equation \eqref{3.17a}, and by  performing on the resulting
solution of \eqref{3.17a} the conformal mapping \eqref{3.27}
followed by the variable transformation \eqref{3.29}, one
obtains the following solution for \eqref{3.30a},
\begin{equation}\lb{3.52}
\varphi (z,t)= (iI_m \mp M)(iI_m \pm M)^{-1}
\exp(\mp 2ite^{i\delta}),
\end{equation}
for $\pm t\ge 0$, $\Im (\pm M)\ge 0$, and $z\in \bbC_+$.
By hypothesis, $0<\delta<\pi$.  Thus  the exponential
term in \eqref{3.52} will result in
\begin{equation}
|| \varphi (z,t) ||_{\bbC^{m\times m}} >1 \ \text{ as } t\to \pm\infty
\end{equation}
unless
\begin{equation}
M = \pm i I_m,
\end{equation}
thus implying \eqref{3.51}.
\end{proof}
\noindent One then obtains the following result.
\begin{corollary}\lb{c3.10}
With $\phi_\pm(\delta) $ defined in \eqref{3.33},
$\eta_\pm(0)=\phi_\pm(\delta)=0$.
\end{corollary}
For $M(z_n,c,x_0,\alpha_0)\in\calD(z_n,c,x_0,\alpha_0)$,
it follows by \eqref{3.30b}, \eqref{3.33}, and
Corollary~\ref{c3.10} that
\begin{equation}\lb{3.55}
[iI_m \mp M(z_n,c,x_0,\alpha_0)][iI_m
\pm M(z_n,c,x_0,\alpha_0)]^{-1} = \oh(1), \qquad\pm (c-x_0)>0,
\end{equation}
as $n\to \infty$. Hence one infers,
for elements of $\calD(z_n,c,x_0,\alpha_0)$, that
\begin{equation}\lb{3.56}
M(z_n,c,x_0,\alpha_0)= \pm iI_m + \oh(1), \qquad \pm (c-x_0)>0,
\end{equation}
as $|z|\to \infty$ in $C_\varepsilon$. This proves \eqref{3.1}. 
Actually, \eqref{3.56} is a statement for all elements of
$\calD(z,c,x_0,\alpha_0)$ including the particular element
$M_{\pm}(z,x_0,\alpha_0)$, for $\pm (c-x_0)>0$. 

In \eqref{3.1} an asymptotic expansion is given that
is uniform with respect to $\arg(z)$ for $|z| \to \infty$ in
$C_\varepsilon$.   We now vary the reference point, $x_0$, and
observe that the asymptotic expansion in \eqref{3.1} is also
uniform with respect to $x_0$ whenever $x_0$ is  confined to a
compact subset of $\bbR$.

\begin{theorem} \lb{t3.12} Assume Hypothesis~\ref{h2.4}. 
Let $\alpha_0 =(I_m\; 0)\in\bbC^{m\times 2m}$, and denote by
$C_\varepsilon\subset \bbC_+$  the open sector with vertex at zero,
symmetry axis along the positve imaginary axis and opening
angle $\varepsilon$, with $0<\varepsilon<\pi/2$. Let
$M_\pm(z,x_0,\alpha_0)$ be the unique elements of the limit
disks $\calD_\pm (z,x_0,\alpha_0)$ for the Dirac system given by
\eqref{HS} and \eqref{DS}. Then,
\begin{equation}\lb{3.57}
M_\pm(z,x,\alpha_0) \underset{\substack{|z|
\to\infty\\z\in C_\varepsilon}}{=} \pm iI_m +\oh(1)
\end{equation}
uniformly with respect to $\arg(z)$, for $|z|
\to \infty$ in $C_\varepsilon$, and uniformly with respect to
$x$, as long as $x$ varies in compact subsets of $[x_0,\pm\infty)$.
\end{theorem}
\begin{proof}
We note that the system \eqref{3.35} is independent of
the reference point $x_0$. Next, we recall that $\delta$, defined in
\eqref{3.32} is determined by an apriori choice of the
sequence $z_n$, subject only to $z_n$ being in $C_\varepsilon$
(c.f.~\eqref{3.31}). Moreover, we note that $\varphi_\pm(\delta)$,
defined as a limit in \eqref{3.33}, described explicity in
Corollary~\ref{c3.10}, and which gives solutions of
\eqref{3.35} satisfying \eqref{3.36} for $t\in [0, \pm\infty)$,
is also independent of the reference point
$x_0$.  Thus, had we chosen a different point of reference,
$x_0'\ne x_0$, at the start, the asymptotic analysis begun in
Theorem~\ref{t3.8} and continued through \eqref{3.55}, would
remain the same after the variable change in \eqref{3.29},
except for the integral expression in \eqref{3.42} in which
$x_0$ would be replaced by $x_0'$.  However, given the local
integrability assumption on $B$ in Hypothesis~\ref{h2.1},
one concludes that this integral expression is uniformly continuous
with respect to $x_0$ whenever $x_0$ is confined to a compact
subset of $\bbR$. Thus \eqref{3.42}, and consequently \eqref{3.50}, are
uniform with respect to $t$ and with respect to $x_0$ whenever both are
confined to compact subsets of $\bbR$. Consequently, \eqref{3.55}
holds for elements $\calD(z,c,x_0,\alpha_0)$, that this asymptotic
expansion is uniform with respect to $\arg (z)$  for $|z|\to \infty$ in
$C_\varepsilon$, and that it is uniform with respect to $x_0$ when $x_0$
is confined to compact subsets of $\bbR$.
\end{proof}

\begin{remark}\lb{r3.2} 
(i) In the special case $m=1$, the leading-order asymptotics \eqref{3.57} 
was published by Everitt, Hinton, and Shaw \cite{EHS83} in 1983.
For asymptotic estimates of Weyl solutions in the case
$m=1$ we refer to \cite{Mi91}. \\ 
(ii) A comparison of \eqref{3.57} with
\eqref{2.41} then  proves that the leading-order asymptotic behavior
\eqref{3.57} is in fact independent of the boundary condition  at $x_0$
indexed by
$\alpha$, that is, 
\begin{equation}\lb{3.58}
M_\pm(z,x_0,\alpha) \underset{\substack{\abs{z}\to\infty\\ z\in
C_\varepsilon}}{=} \pm iI_{m} +\oh(1)
\end{equation}
for any $\alpha$ satisfying the conditions stated in
\eqref{2.8e}. In the scalar case $m=1$ this fact had been
noticed in \cite{EHS83}. This boundary condition independence of
the leading-order asymptotic behavior of $M_\pm(z,x_0,\alpha)$ is in
sharp contrast to the case of matrix-valued Schr\"odinger operators 
(see, e.g., \cite{CG99}). Moreover, regarding the conclusion of
Theorem~\ref{t3.12},  no generality is lost by assuming that
$C_\varepsilon \subset \bbC_+$ because of \eqref{2.38}.
\end{remark}

\section{Higher Order Terms in the Asymptotic Expansion of
$M_\pm(z,x,\alpha)$} \lb{s4}

In this section we shall prove one of our principal
results of  this paper, the asymptotic high-energy expansion of 
$M_+(z,x,\alpha_0)$ to arbitrarily high orders in sectors of the type
$C_\varepsilon\subset\bbC_+$ as defined in Theorem~\ref{t3.12}. 

Throughout this section we choose $z\in\bbC_+$. We also recall the
following notion: $x\in [a,b)$ (resp.,
$x\in (a,b]$) is called a  right (resp., left) Lebesgue point of an
element $q\in L^1 ((a,b))$, $a<b$ if$\int_0^\varepsilon dx^\prime \, 
|q(x+x^\prime)-q(x)|=\oh (\varepsilon)$ (resp., $\int_0^\varepsilon
dx^\prime \, |q(x-x^\prime)-q(x)|=\oh (\varepsilon)$)  as
$\varepsilon\downarrow 0$. Similarly, $x\in (a,b)$ is called a Lebesgue
point of $q\in L^1 ((a,b))$ if $\int_{-\varepsilon}^\varepsilon dx^\prime
\,  |q(x+x^\prime)-q(x)|=\oh (\varepsilon)$ as $\varepsilon\downarrow 0$. 
The set of all such points is then denoted the right (resp., left) 
Lebesgue set of $q$ on $[a,b]$ in the former case and simply the Lebesgue
set of $q$ on $[a,b]$ in the latter case. The
analogous notions are applied to $2m\times 2m$ matrices $B\in L^1
((a,b))^{2m\times 2m}$ by simultaneously considering all $4m^2$ entries of
$B$. The right (resp., left) Lebesgue set of $B$ on $[a,b]$ is then
simply the intersection of the right (resp., left) Lebesgue sets of
$B_{j,k}$ for all $1\leq j,k\leq 2m$, and similarly for the Lebesgue 
set of $B$, etc.

Finally, we need one more ingredient, recently proven by Rybkin 
\cite[Lemma~3]{Ry99} using appropriate maximal functions. Let $q\in L^1
((x_0,\infty))$, $\supp(q)\subseteq [x_0,x_0+R]$ for some $R>0$, and
suppose $x\in [x_0,x_0+R]$ is a right Lebesgue point of $q$. Then
\begin{equation}
\int_x^{x_0+R} dx^\prime \, q(x^\prime)\exp(2iz(x^\prime -x))
\underset{\substack{\abs{z}\to\infty\\ z\in
C_\varepsilon}}{=}-\f{q(x)}{2iz} + \oh\big(|z|^{-1}\big). \lb{4.-2} 
\end{equation}
An alternative proof of \eqref{4.-2} follows from 
\cite[Theorem~I.13]{Ti86}, which implies
\begin{equation}
\lim_{\substack{\abs{z}\to\infty\\ z\in
C_\varepsilon}}z^{-1} \int_x^{x_0+R} dx^\prime \,  
|q(x^\prime) -q(x)|\exp(2iz(x^\prime -x)) = 0 \lb{4.-1}
\end{equation}
for any right Lebesgue point $x$ of $q$.
 
We start with the simpler case where $B$ has compact support contained in
some interval $[x_0,y_0]$. Below in \eqref{4.0} and in
analogous formulas in this section, $\|\cdot\|_{\bbC^{\ell\times \ell}}$
denotes a norm in $\bbC^{\ell\times \ell}$.

\begin{lemma}\lb{l4.1} 
Fix $x_0, y_0\in\bbR$ with $y_0>x_0$ and let $x\geq x_0$. Suppose
$A=I_{2m}$, $B\in L^1([x_0,x_0+R])^{2m\times 2m}$ for all $R>0$,
$B=B^*$ a.e.~on $(x_0,\infty)$. In addition, assume that $B$ has compact
support contained in $[x_0,y_0]$, that
$B^{(N-1)}\in L^1([x_0,y_0])^{2m\times 2m}$ for some $N\in\bbN$, that $x$
is a right Lebesgue point of $B^{(N-1)}$, and that
\begin{align}
&\underset{y\in[x_0,y_0]}{\esssup} \, \bigg\|\int_y^{y_0}
dx'\,B^{(N-1)}(x')\exp(2iz(x'-y))
+\f{1}{2iz}B^{(N-1)}(y)\bigg\|_{\bbC^{2m\times 2m}} \no \\
& \underset{\substack{\abs{z}
\to\infty\\ z\in C_\varepsilon}}{=}\oh\big(|z|^{-1}\big). \lb{4.0}
\end{align}
If $N=1$, suppose in addition $B_{k,k'}B_{\ell,\ell'}\in
L^1([x_0,y_0])^{m\times m}$ for all $k,k',\ell,\ell'\in\{1,2\}$. Let 
$\alpha_0=(I_m\; 0)\in\bbC^{m\times 2m}$ and 
denote by $M_+(z,x,\alpha_0)$, $x\geq x_0$, the unique 
Weyl-Titchmarsh matrix associated with the half-line
Dirac-type operator $D_+(\alpha_0)$ in \eqref{2.84}. Then, as
$\abs{z}\to\infty$ in $C_\varepsilon$, $M_+(z,x,\alpha_0)$ has
an asymptotic expansion of the form 
\begin{equation}
M_+(z,x,\alpha_0)\underset{\substack{\abs{z}\to\infty\\ z\in
C_\varepsilon}}{=}
 i I_m +\sum_{k=1}^N m_{+,k}(x,\alpha_0)z^{-k}+
\oh\big(|z|^{-N}\big), \quad N\in\bbN. \lb{4.1}
\end{equation}
The expansion \eqref{4.1} is uniform with respect to $\arg\,(z)$ for $|z|
\to \infty$ in
$C_\varepsilon$ and uniform in $x$ as long as $x$ varies in compact
subintervals of $[x_0,\infty)$ intersected with  the right Lebesgue set
of $B^{(N-1)}$.  The expansion  coefficients $m_{+,k}(x,\alpha_0)$ can be
recursively computed from
\begin{align}
m_{+,1}(x,\alpha_0)&=-\f{1}{2} \big(
B_{1,2}(x)+B_{2,1}(x)\big)  +\f{i}{2} \big(
B_{1,1}(x)-B_{2,2}(x)\big), \no \\
m_{+,k+1}(x,\alpha_0)&=\f{i}2\bigg(m_{+,k}^\prime(x,\alpha_0)+
\sum_{\ell=1}^{k} m_{+,\ell}(x,\alpha_0)
m_{+,k+1-\ell}(x,\alpha_0) \no \\
& \qquad \quad +\sum_{\ell=0}^{k} m_{+,\ell}(x,\alpha_0) 
B_{2,2}(x) m_{+,k-\ell}(x,\alpha_0) \lb{4.2}  \\
& \qquad \quad + B_{1,2}(x) m_{+,k}(x,\alpha_0)
+ m_{+,k}(x,\alpha_0) B_{2,1}(x)\bigg), \no \\
& \hspace*{5.7cm} 1 \leq k\leq N-1. \no
\end{align}
\end{lemma}
\begin{proof}
In the following let $z\in\bbC_+$, and $x\geq x_0$.
The existence of an expansion of the type \eqref{4.1} is shown as
follows.  First one considers a matrix Volterra integral equation of the
type 
\begin{equation}
\widetilde U_+(z,x,\alpha_0)=\begin{pmatrix}I_m\\
iI_m\end{pmatrix}\exp(iz(x-x_0)) +\int_x^\infty dx'\, K(z,x,x') J
B(x')\widetilde U_+(z,x',\alpha_0),  \lb{4.3} 
\end{equation}
where 
\begin{equation}
\widetilde U_+(z,x,\alpha_0)=\begin{pmatrix} \widetilde
u_{+,1}(z,x,\alpha_0)\\
\widetilde u_{+,2}(z,x,\alpha_0)\end{pmatrix}\in 
L^2([x_0,\infty))^{2m\times m}, 
\lb{4.3A} 
\end{equation}
and $K$ abbreviates the $2m\times 2m$ Volterra Green's kernel
\begin{equation}
K(z,x,x')=\begin{pmatrix} \cos(z(x-x'))I_m &\sin(z(x-x'))I_m \\
-\sin(z(x-x'))I_m &\cos(z(x-x'))I_m \end{pmatrix}. \lb{4.3B}
\end{equation}
Clearly, $\wti U_+(z,\cdot,\alpha_0)$ solves the
Dirac-type system \eqref{HS} and \eqref{DS}. In addition, it satisfies 
$\wti U_+(z,\cdot,\alpha_0)\in L^2 ([x_0,\infty))^{2m\times 2m}$.
Thus, up to normalization, $\wti U_+(z,\cdot,\alpha_0)$ represents
the Weyl solution associated with $B$ on the half-line
$[x_0,\infty)$. Next, introducing
\begin{equation}
\widetilde V_+(z,x,\alpha_0)=\begin{pmatrix} \wti
v_{+,1}(z,x,\alpha_0)\\
\wti v_{+,2}(z,x,\alpha_0) \end{pmatrix} =\widetilde
U_+(z,x,\alpha_0)\exp(-iz(x-x_0)), \lb{4.3a}
\end{equation}
one rewrites \eqref{4.3} in the form
\begin{equation}
\widetilde V_+(z,x,\alpha_0)=\begin{pmatrix}I_m\\ iI_m\end{pmatrix}
+\int_x^{y_0} dx'\, \widetilde K(z,x,x') J B(x')\widetilde
V_+(z,x',\alpha_0), \lb{4.3b} 
\end{equation}
where 
\begin{equation}
\widetilde K(z,x,x')=\f{1}{2}\begin{pmatrix} (1+\exp(2iz(x'-x)))I_m
&-i(1-\exp(2iz(x'-x)))I_m  \\[1mm] 
i(1-\exp(2iz(x'-x)))I_m &(1+\exp(2iz(x'-x)))I_m \end{pmatrix}. \lb{4.3C}
\end{equation}
Thus, one infers, 
\begin{equation}
M_+(z,x,\alpha_0)=\wti u_{+,2}(z,x,\alpha_0)\wti
u_{+,1}(z,x,\alpha_0)^{-1}  =\wti v_{+,2}(z,x,\alpha_0)\wti
v_{+,1}(z,x,\alpha_0)^{-1}.
\lb{4.4}
\end{equation}
Introducing
\begin{equation}
R=\begin{pmatrix} C_1 & -iC_2 \\ iC_1 & C_2
\end{pmatrix}, \quad  S=\begin{pmatrix} D_1 & iD_2 \\
-iD_1 & D_2 \end{pmatrix}, \lb{4.4a}
\end{equation}
where
\begin{align}
C_1&=-B_{1,2}^* -iB_{1,1}, \quad C_2=
B_{1,2}-iB_{2,2}, \lb{4.4b} \\ 
D_1&=-B_{1,2}^*+iB_{1,1}, \quad D_2=
B_{1,2}+iB_{2,2}, \lb{4.4c}
\end{align}
\eqref{4.3b} results in 
\begin{align}
\widetilde V_+(z,x,\alpha_0)&=\begin{pmatrix}I_m\\ iI_m\end{pmatrix}
+\int_x^{y_0} dx'\, \big(R(x')+S(x')\exp(2iz(x'-x)) 
\big) \widetilde V_+(z,x',\alpha_0) \lb{4.4d} \\
&=\bigg(I_{2m}+\sum_{k=1}^\infty 2^{-k}\int_x^{y_0}  dx_1 \, 
\big(R(x_1)+S(x_1)e^{2iz(x_1-x)} \big) \times \no \\
& \hspace*{2.5cm} \times \int_{x_1}^{y_0}  dx_2
\, \big(R(x_2)+S(x_2)e^{2iz(x_2-x_1)} \big)\dots \lb{4.4e} \\
& \hspace*{2.3cm} \dots \int_{x_{k-1}}^{y_0}  dx_k
\, \big(R(x_k)+S(x_k)e^{2iz(x_k-x_{k-1})} \big)
\bigg)\begin{pmatrix}I_m\\ iI_m\end{pmatrix}. \no
\end{align}
This yields
\begin{equation}
\|\wti v_{+,j}(z,x,\alpha_0) \| \leq C_j, \quad z\in\bbC_+, \; 
\Im(z)> 0, \; x\geq x_0, \; j=1,2 \lb{4.4ca}
\end{equation}
for some $C_j>0$, $j=1,2$, depending on $\|B\|_1$. Integrating
by parts in \eqref{4.4e}, repeatedly applying \eqref{4.-2} and \eqref{4.0}
to $q(x)=(S(x))_{j,k}$ for all $1\leq j,k\leq 2m$ then results in the
existence of an asymptotic expansion for $\widetilde V_+(z,x,\alpha_0)$
of the type
\begin{equation}
\widetilde V_+(z,x,\alpha_0)=\begin{pmatrix} \wti v_{+,1}(z,x,\alpha_0)\\
\wti v_{+,2}(z,x,\alpha_0) \end{pmatrix}=\sum_{k=0}^{N} \widetilde
V_{+,k}(x,\alpha_0)\, z^{-k} +\oh\big(|z|^{-N}\big). \lb{4.3f}
\end{equation}
Inserting the expansions for $\wti v_{+,2}(z,x,\alpha_0)$ and
$\wti  v_{+,1}(z,x,\alpha_0)^{-1}$ into \eqref{4.4} (using a geometric
series expansion for $\wti v_{+,1}(z,x,\alpha_0)^{-1}$) then yields the
existence of an expansion of the type \eqref{4.1} for
$M_+(z,x,\alpha_0)$. The actual expansion coefficients and the
associated recursion relation \eqref{4.2} then follow upon inserting
expansion \eqref{4.1} into the Riccati-type equation \eqref{3.17a}. The
stated uniformity assertions concerning the asymptotic expansion
\eqref{4.1} then follow from iterating the system of Volterra integral
integral equations
\eqref{4.3b}.
\end{proof}

\begin{remark} \lb{r4.2}
The analogous solution $\wti U_-(z,\cdot,\alpha_0)$ of the Dirac-type
operator \eqref{2.61} on the interval $(-\infty,x_0]$ satisfies 
\begin{align}
\widetilde U_-(z,x,\alpha_0)&=\begin{pmatrix}I_m\\
-iI_m\end{pmatrix}\exp(-iz(x-x_0)) \no \\
& \quad -\int_{-\infty}^x dx'\, K(z,x,x') J
\widetilde B(x')\widetilde U_-(z,x',\alpha_0),  \lb{4.14A} 
\end{align}
with integral kernel $K$ given by \eqref{4.3B}. (Again $\wti U_-$
coincides with the Weyl solution $U_-$ up to normalization.) A closer look
at the system of Volterra integral equations
\eqref{4.3},
\eqref{4.4d},
\eqref{4.4e}, and similarly in connection with \eqref{4.14A}, then 
reveals that $\wti U_\pm(z,\cdot,\alpha_0)$ have the asymptotic behavior
\begin{equation}
\wti U_\pm (z,x,\alpha_0)\underset{\substack{\abs{z}\to\infty\\ z\in
C_\varepsilon}}{=} \left(\sum_{k=0}^N \begin{pmatrix} \wti
v_{\pm,k,1}(x,\alpha_0) \\ \wti v_{\pm,k,2}(x,\alpha_0)
\end{pmatrix}z^{-k} +\oh\big(|z|^{-N}\big)\right)\exp(\pm iz(x-x_0)),
\lb{4.14B}
\end{equation}
with leading asymptotics determined as follows.
\begin{align} 
\begin{split} \lb{4.14C} 
\wti v_{\pm,0,1}(x,\alpha_0)&=I_m+\wti w_{\pm,0,1}(x,\alpha_0), \\ 
\wti v_{\pm,0,2}(x,\alpha_0)&=\pm i\big(I_m+\wti
w_{\pm,0,1}(x,\alpha_0)\big),
\end{split} 
\end{align}
where $\wti w_{\pm,0,1}(x,\alpha_0)$ satisfies
\begin{equation}
\wti w_{\pm,0,1}'(x,\alpha_0)=\f{1}{2}\big[\wti B_{2,1}(x) -\wti
B_{1,2}(x) \pm i\wti B_{2,2}(x) \pm i\wti B_{1,1}(x)\big]\big(I_m+\wti
w_{\pm,0,1}(x,\alpha_0)\big), \lb{4.14D}
\end{equation}
and
\begin{equation}
\lim_{x\to\pm\infty} \wti w_{\pm,0,1}(x,\alpha_0)=0  
\lb{4.14F}
\end{equation}
(in fact, $\wti v_{\pm,0,1}(\cdot,\alpha_0)=I_m$, $\wti
v_{\pm,0,2}(\cdot,\alpha_0)=\pm iI_m$, and $\wti
v_{\pm,k,j}(\cdot,\alpha_0)=0$, $j=1,2$, $1\leq k\leq N$ outside the
support of $\wti B$). In particular,
\begin{equation}
\wti w_{\pm,0,1}(x,\alpha_0)=0 \lb{4.14G}
\end{equation}
and hence 
\begin{equation}
\wti U_\pm (z,x,\alpha_0)\underset{\substack{\abs{z}\to\infty\\ z\in
C_\varepsilon}}{=} \left(\begin{pmatrix} I_m \\ \pm iI_m
\end{pmatrix} +\oh(1)\right)\exp(\pm iz(x-x_0)), \lb{4.14H}
\end{equation}
if and only if $\wti B$ is in the normal form
\begin{equation}
\wti B(x)=\begin{pmatrix} \wti B_{1,1}(x) & \wti B_{1,2}(x)\\
\wti B_{1,2}(x) &-\wti B_{1,1}(x)\end{pmatrix}, \quad
\wti B_{1,1}^*(x)=\wti B_{1,1}(x),
\; \wti B_{1,2}^*(x)=\wti B_{1,2}(x) \text{ a.e.} \lb{4.14I}
\end{equation} 
\end{remark}

For more details we refer to Lemma~\ref{l4.9}.

Next we recall an elementary result on finite-dimensional
evolution equations essentially taken from \cite{MPS90} (cf.~also 
\cite[Lemma~4.2]{CG99}).

\begin{lemma} \mbox{\rm (\cite{MPS90}.)}
\lb{l4.2} Let $\Gamma_j\in L^1_{\loc}(\bbR)^{m\times m}$, $j=1,2$. Then
any
$m\times m$ matrix-valued solution $X$ of
\begin{equation}
X'(x)=\Gamma_1(x)X(x)+X(x)\Gamma_2(x) \text{ for ~a.e. } x\in\bbR,
\lb{4.5}
\end{equation}
is of the type
\begin{equation}
X(x)=Y(x)CZ(x), \lb{4.6}
\end{equation}
where $C$ is a constant $m\times m$ matrix and $Y$ is a
fundamental system of solutions of
\begin{equation}
\Psi'(x)=\Gamma_1(x)\Psi(x) \lb{4.7}
\end{equation}
and $Z$ is a fundamental system of solutions of
\begin{equation}
\Phi'(x)=\Phi(x)\Gamma_2(x). \lb{4.8}
\end{equation}
\end{lemma}

The next result provides the proper extension of Lemma~4.3 in 
\cite{CG99} in the context of matrix-valued Schr\"odinger operators 
(which in turn extended Proposition~2.1 in the scalar context in
\cite{GS98} to the matrix-valued case) to the Dirac-type case under
consideration.

\begin{lemma} \lb{l4.3} 
Fix $x_0,y_0\in\bbR$ with $y_0>x_0$. Suppose
$A_j=I_{2m}$, $B_j\in L^1([x_0,x_0+R])^{2m\times 2m}$ for all $R>0$, 
$B_j=B_j^*$ a.e.~on $[x_0,\infty)$, $j=1,2$, and $B_1=B_2$ a.e.~on
$[x_0,y_0]$. Let $\alpha_0=(I_m\; 0)\in\bbC^{m\times 2m}$ and denote by
$M_{j,+}(z,x,\alpha_0)$, $x\geq x_0$, the unique Weyl-Titchmarsh matrix 
corresponding to the half-line Dirac operators $D_{+,j}(\alpha_0)$,
$j=1,2$, in \eqref{2.84}. Then,
\begin{align}
&[M_{1,+}'(z,x,\alpha_0)-M_{2,+}'(z,x,\alpha_0)] \no \\
&=-(z/2)[M_{1,+}(z,x,\alpha_0)+M_{2,+}(z,x,\alpha_0)]
[M_{1,+}(z,x,\alpha_0)-M_{2,+}(z,x,\alpha_0)] \no \\
& \quad -(z/2)[M_{1,+}(z,x,\alpha_0)-M_{2,+}(z,x,\alpha_0)]
[M_{1,+}(z,x,\alpha_0) +M_{2,+}(z,x,\alpha_0)]  \no \\ 
& \quad
-[M_{1,+}(z,x,\alpha_0)+M_{2,+}(z,x,\alpha_0)]B_{2,2}(x)
[M_{1,+}(z,x,\alpha_0)-M_{2,+}(z,x,\alpha_0)]/2 \no \\ 
& \quad
-[M_{1,+}(z,x,\alpha_0)-M_{2,+}(z,x,\alpha_0)]B_{2,2}(x)
[M_{1,+}(z,x,\alpha_0)+M_{2,+}(z,x,\alpha_0)]/2  \no \\ 
& \quad -B_{1,2}(x)[M_{1,+}(z,x,\alpha_0)-M_{2,+}(z,x,\alpha_0)] 
\no \\
& \quad -[M_{1,+}(z,x,\alpha_0)-M_{2,+}(z,x,\alpha_0)]B_{2,1}(x) 
\; \text{ for~a.e. $x\in [x_0,y_0]$,} \lb{4.14}
\end{align}
where we denoted $B_1=B_2=\left(\begin{smallmatrix}B_{1,1}
&B_{1,2}\\ B_{2,1} &B_{2,2} \end{smallmatrix}\right)$ a.e.~on $(x_0,y_0)$.
\end{lemma}
\begin{proof}
This is obvious from \eqref{3.17a}.
\end{proof}
\begin{lemma} \lb{l4.4} 
Fix $x_0,y_0\in\bbR$ with $y_0>x_0$. Suppose
$A_j=I_{2m}$, $B_j\in L^1([x_0,x_0+R])^{2m\times 2m}$ for all
$R>0$, and $B_j=B_j^*$ a.e.~on $[x_0,\infty)$, $j=1,2$. Let
$\alpha_0=(I_m\; 0)\in\bbC^{m\times 2m}$ and denote by
$M_{j,+}(z,x,\alpha_0)$, $x\geq x_0$, the unique Weyl-Titchmarsh matrix 
corresponding to the half-line Dirac operators $D_{+,j}(\alpha_0)$,
$j=1,2$, in \eqref{2.84}. Define
\begin{align}
\Gamma_1(z,x)&=-(z/2)[M_{1,+}(z,x,\alpha_0)
+M_{2,+}(z,x,\alpha_0)] \no \\
& \quad -(1/2)[M_{1,+}(z,x,\alpha_0)
+M_{2,+}(z,x,\alpha_0)]B_{2,2}(x)-B_{1,2}(x), \lb{4.15} \\
\Gamma_2(z,x)&=-(z/2)[M_{1,+}(z,x,\alpha_0)
+M_{2,+}(z,x,\alpha_0)] \no \\
& \quad -(1/2)B_{2,2}(x)[M_{1,+}(z,x,\alpha_0)
+M_{2,+}(z,x,\alpha_0)]-B_{2,1}(x), \lb{4.15a}
\end{align}
for a.e.~$x\in [x_0,y_0]$. In addition, assume $Y_+(z,\cdot)$
and $Z_+(z,\cdot)$ to be fundamental matrix solutions of
\begin{equation}
\Psi'(z,x)=\Gamma_1(z,x)\Psi(z,x) \text{ and }
\Phi'(z,x)=\Phi(z,x)\Gamma_2(z,x) \lb{4.16}
\end{equation}
on $[x_0,y_0]$, respectively, with
\begin{equation}
Y_+(z,y_0)=I_m, \quad Z_+(z,y_0)=I_m. \lb{4.17}
\end{equation}
Then, as $|z|\to\infty$, $z\in C_\varepsilon$,
\begin{equation}
\|Y_+(z,x_0)\|_{\bbC^{m\times m}}, \|Z_+(z,x_0)\|_{\bbC^{m\times m}}
\leq \exp(-\Im(z)(y_0-x_0)(1+\oh (1))). \lb{4.18}
\end{equation}
\end{lemma}
\begin{proof}
Define $\wti\Gamma_j(z,x)$, $j=1,2$, by
\begin{equation}
\wti\Gamma_j (z,x)=\Gamma_j (z,x)+izI_m, \quad j=1,2, \lb{4.19}
\end{equation}
then
\begin{equation}
\int_{x_0}^{y_0} dx\, \|\wti\Gamma_j (z,x)\|_{\bbC^{m\times m}}
\underset{\substack{\abs{z}\to\infty\\ z\in
C_\varepsilon}}{=}  \oh (z), \quad j=1,2 \lb{4.20}
\end{equation}
due to the uniform nature of the asymptotic expansion
\eqref{3.57} for $x$ varying in compact intervals. Next,
introduce
\begin{equation}
E_+(z,x,y_0)=I_m\exp(iz(y_0-x)), \quad x\leq y_0,
\lb{4.21}
\end{equation}
then
\begin{align}
Y_+(z,x)&=E_+(z,x,y_0)-\int_x^{y_0}dx'\,E_+(z,x,x')
\wti\Gamma_1(z,x')Y_+(z,x'), \lb{4.22} \\
Z_+(z,x)&=E_+(z,x,y_0)-\int_x^{y_0}dx'\,Z_+(z,x')
\wti\Gamma_2(z,x')E_+(z,x,x'). \lb{4.23}
\end{align}
Using
\begin{equation}
\|E_+(z,x_0,y_0)\|_{\bbC^{m\times m}}\leq\exp(-\Im(z)(y_0-x_0)),
\lb{4.24}
\end{equation}
a standard Volterra-type iteration argument in
\eqref{4.22}, \eqref{4.23} then yields
\begin{align}
\|Y_+(z,x_0)\|_{\bbC^{m\times m}}& \leq
\exp\left(-\Im(z)(y_0-x_0)+\int_{x_0}^{y_0}dx
\,\|\wti\Gamma_1(z,x)\|\right), \lb{4.25} \\
\|Z_+(z,x_0)\|_{\bbC^{m\times m}}& \leq
\exp\left(-\Im(z)(y_0-x_0)+\int_{x_0}^{y_0}dx
\,\|\wti\Gamma_2(z,x)\|\right), \lb{4.25a}
\end{align}
and hence \eqref{4.18}.
\end{proof}
\begin{theorem} \lb{t4.5} 
Fix $x_0,y_0\in\bbR$ with $y_0>x_0$. Suppose
$A_j=I_{2m}$, $B_j\in L^1([x_0,x_0+R])^{2m\times 2m}$ for all $R>0$, 
$B_j=B_j^*$ a.e.~on $[x_0,\infty)$, $j=1,2$, and $B_1=B_2$ a.e.~on
$[x_0,y_0]$. Let
$\alpha_0=(I_m\; 0)\in\bbC^{m\times 2m}$ and denote by
$M_{j,+}(z,x,\alpha_0)$, $x\geq x_0$, the unique Weyl-Titchmarsh matrix 
corresponding to the half-line Dirac operators $D_{+,j}(\alpha_0)$,
$j=1,2$, in \eqref{2.84}. Then, as $\abs{z}\to\infty$ in $C_\varepsilon$,
\begin{equation}
\|M_{1,+}(z,x_0,\alpha_0)-M_{2,+}(z,x_0,\alpha_0)\|_{\bbC^{m\times m}}
\leq C\exp(-2\Im(z)(y_0-x_0)(1+\oh (1))) \lb{4.26}
\end{equation}
for some constant $C>0$.
\end{theorem}
\begin{proof}
Define for $z\in\bbC\backslash\bbR$, $x\in [x_0,y_0]$, 
\begin{equation}
X_+(z,x)=M_{1,+}(z,x,\alpha_0)-M_{2,+}(z,x,\alpha_0), \lb{4.26a} 
\end{equation}
and for $z\in\bbC\backslash\bbR$ and a.e.~$x\in [x_0,y_0]$,
\begin{align}
\Gamma_1(z,x)&=-(z/2)[M_{1,+}(z,x_0,\alpha_0)+M_{2,+}(z,x_0,\alpha_0)] 
\no \\
& \quad -(1/2)[M_{1,+}(z,x_0,\alpha_0)+M_{2,+}(z,x_0,\alpha_0)]B_{2,2}(x)
-B_{1,2}(x), \lb{4.26b} \\
\Gamma_2(z,x)&=-(z/2)[M_{1,+}(z,x_0,\alpha_0)+M_{2,+}(z,x_0,\alpha_0)] 
\no \\
& \quad -(1/2)B_{2,2}(x)[M_{1,+}(z,x_0,\alpha_0)+M_{2,+}(z,x_0,\alpha_0)]
-B_{2,1}(x). \lb{4.26c}
\end{align}
By Lemma~\ref{l4.3},
\begin{equation}
X_+'=\Gamma_1X_++X_+\Gamma_2
\end{equation}
and hence by Lemma~\ref{l4.2},
\begin{equation}
X_+(z,x)=Y_+(z,x)X_+(z,x_1)Z_+(z,x), \lb{4.30}
\end{equation}
where $Y_+(z,x)$ and $Z_+(z,x)$ are fundamental solution matrices of
\begin{equation}
\Psi'(z,x)=\Gamma_1(z,x)\Psi(z,x) \text{ and }
\Phi'(z,x)=\Phi(z,x)\Gamma_2(z,x),
\end{equation}
respectively, with
\begin{equation}
Y_+(z,y_0)=I_m, \quad Z_+(z,y_0)=I_m.
\end{equation}
By Lemma~\ref{l4.4},
\begin{equation}
\|Y_+(z,x_0)\|_{\bbC^{m\times m}}, \|Z_+(z,x_0)\|_{\bbC^{m\times m}} 
\leq \exp(-\Im(z)(y_0-x_0))(1+\oh (1))) \lb{4.33}
\end{equation}
as $|z|\to\infty$, $z\in C_\varepsilon$. Thus, as
$|z|\to\infty$, $z\in C_\varepsilon$,
\begin{align}
\|X_+(z,x_0)\|_{\bbC^{m\times m}}&\leq
\|X_+(z,y_0)\|_{\bbC^{m\times m}}\,\|Y_+(z,x_0)\|_{\bbC^{m\times
m}}\,\|Z_+(z,x_0)\|_{\bbC^{m\times m}} \no \\ 
&\leq C\exp(-2\Im(z)(y_0-x_0)(1+\oh(1)))
\lb{4.34}
\end{align}
for some constant $C>0$ by \eqref{3.57}, \eqref{4.30}, and
\eqref{4.33}.
\end{proof}

Given these preparations we can now drop the compact
support assumption on $B$ in Lemma~\ref{l4.1} and hence arrive at one 
of the principal results of this paper.

\begin{theorem} \lb{t4.6} 
Fix $x_0,y_0\in\bbR$ with $y_0>x_0$ and suppose $A=I_{2m}$, $B\in
L^1([x_0,x_0+R])^{2m\times 2m}$ for all $R>0$, and $B=B^*$ a.e.~on
$(x_0,\infty)$. In addition, assume that for some $N\in\bbN$,
$B^{(N-1)}\in L^1([x_0,c])^{2m\times 2m}$ for all $c>x_0$, that $x_0$ is
a right Lebesgue point of $B^{(N-1)}$, and that
\begin{align}
&\underset{y\in [x_0,y_0]}{\esssup} \, \bigg\|\int_y^{y_0}
dx'\,B^{(N-1)}(x')\exp(2iz(x'-y))
+\f{1}{2iz}B^{(N-1)}(y)\bigg\|_{\bbC^{2m\times 2m}} \no \\
&\underset{\substack{\abs{z}
\to\infty\\ z\in C_\varepsilon}}{=}\oh\big(|z|^{-1}\big). \lb{4.34a}
\end{align}
If $N=1$, suppose in addition 
$B_{k,k'}B_{\ell,\ell'}\in L^1([x_0,y_0])^{m\times m}$ for all
$k,k',\ell,\ell'\in\{1,2\}$. Let $\alpha_0=(I_m\; 0)\in\bbC^{m\times
2m}$ and denote by $M_+(z,x_0,\alpha_0)$ the unique element of the
limit disk $\calD_+ (z,x_0,\alpha_0)$ for the half-line Dirac operator
$D_+(\alpha_0)$ in \eqref{2.84}. Then, as
$\abs{z}\to\infty$ in $C_\varepsilon$, $M_+(z,x_0,\alpha_0)$ has an
asymptotic expansion of the form 
\begin{equation}
M_+(z,x_0,\alpha_0)\underset{\substack{\abs{z}\to\infty\\ z\in
C_\varepsilon}}{=} i I_m +\sum_{k=1}^N
m_{+,k}(x_0,\alpha_0)z^{-k}+ o\big(|z|^{-N}\big), \quad N\in\bbN.
\lb{4.35}
\end{equation}
The expansion \eqref{4.35} is uniform with respect to $\arg\,(z)$ 
for $|z|\to \infty$ in $C_\varepsilon$. The expansion  coefficients
$m_{+,k}(x_0,\alpha_0)$ can be recursively computed from \eqref{4.2}.
\end{theorem}
\begin{proof}
Define
\begin{equation} 
\widetilde B(x)=\begin{cases} B(x) &\text{ for }
x\in [x_0,y_0], \; x_0<y_0  \\
0 &\text{ otherwise} \lb{4.36} \end{cases} 
\end{equation}
and apply Theorem~\ref{t4.5} with $B_1=B$, $B_2=\widetilde
B$. Then (in obvious notation)
\begin{equation}
\|M_+(z,x_0,\alpha_0)-\widetilde M_+(z,x_0,\alpha_0)\|_{\bbC^{m\times
m}}\leq C \exp(-2\Im(z)(y_0-x_0)(1+\oh(1))) \lb{4.36a}
\end{equation}
as $|z|\to\infty$, $z\in C_\varepsilon$, and hence the
asymptotic expansion \eqref{4.1} for $\widetilde M_+(z,x_0,\alpha_0)$
in Lemma~\ref{l4.1} coincides with that of $M_+(z,x_0,\alpha_0)$.
\end{proof}

In analogy to Theorem~\ref{t3.12}, the asymptotic expansion
\eqref{4.35} extends to one for $M_+(z,x,\alpha_0)$ valid uniformly
with respect to $x$ as long as $x$ varies in compact subintervals
of $[x_0,\infty)$ intersected with the right Lebesgue set of 
$B^{(N-1)}$.

\begin{theorem} \lb{t4.7} 
Fix $x_0\in\bbR$ and let $x\geq x_0$. Suppose $A=I_{2m}$, $B\in
L^1([x_0,x_0+R])^{2m\times 2m}$ for all $R>0$, and $B=B^*$ a.e.~on
$(x_0,\infty)$. In addition, assume that for some $N\in\bbN$,
$B^{(N-1)}\in L^1([x_0,c))^{2m\times 2m}$ for all $c>x_0$, that $x$ is
a right Lebesgue point of $B^{(N-1)}$, and that for all $R>0$,
\begin{align}
&\underset{y\in [x_0,x_0+R]}{\esssup} \,
\bigg\|\int_y^{x_0+R} dx'\,B^{(N-1)}(x')\exp(2iz(x'-y))
+\f{1}{2iz}B^{(N-1)}(y)\bigg\|_{\bbC^{2m\times 2m}} \no \\
& \underset{\substack{\abs{z}
\to\infty\\ z\in C_\varepsilon}}{=}\oh\big(|z|^{-1}\big). \lb{4.36b}
\end{align}
If $N=1$, suppose in addition 
$B_{k,k'}B_{\ell,\ell'}\in L^1([x_0,x_0+R])^{m\times m}$ for all $R>0$
and all $k,k',\ell,\ell'\in\{1,2\}$. Let
$\alpha_0=(I_m\; 0)\in\bbC^{m\times 2m}$  and denote by
$M_+(z,x,\alpha_0)$, $x\geq x_0$, the unique element of the limit disk
$\calD_+(z,x,\alpha_0)$ for the half-line Dirac operator
$D_+(\alpha_0)$ in \eqref{2.84}. Then, as $\abs{z}\to\infty$ in
$C_\varepsilon$, $M_+(z,x,\alpha_0)$ has an asymptotic expansion of the
form
\begin{equation}
M_+(z,x,\alpha_0)\underset{\substack{\abs{z}\to\infty\\ z\in
C_\varepsilon}}{=}  i I_m +\sum_{k=1}^N
m_{+,k}(x,\alpha_0)z^{-k}+ o\big(|z|^{-N}\big), \quad N\in\bbN. \lb{4.37}
\end{equation}
The expansion \eqref{4.37} is uniform with respect to $\arg\,(z)$ for 
$|z|\to \infty$ in $C_\varepsilon$ and uniform in $x$ as long as $x$
varies in compact subsets of $\bbR$ intersected with the right Lebesgue
set of $B^{(N-1)}$. The expansion coefficients $m_{+,k}(x,\alpha_0)$ can
be recursively computed from \eqref{4.2}.
\end{theorem}
\begin{proof}
To see that uniformity holds for this expansion, first recall the role of
Theorem~\ref{t3.12} in providing uniformity in the asymptotic expression
\eqref{4.20} which then leads to \eqref{4.18} holding uniformly with
respect to $x_0$ varying within compact subsets of $\bbR$ and with respect
to $\arg\,(z)$ for $|z|\to \infty$ in $C_\varepsilon$. This in turn leads
to a similar uniformity holding for \eqref{4.26} which is the key to
\eqref{4.35} holding with respect to $x_0$ varying within compact subsets 
of $\bbR$ and with respect to $\arg\,(z)$ for $|z|\to \infty$ in 
$C_\varepsilon$. 
\end{proof}
\begin{remark} \lb{r4.9}
For simplicity, we focused thus far on the expansion of 
$M_+(z,x_0,\alpha_0)$ as $|z|\to\infty$. Of course, 
Theorem~\ref{t4.7} holds also for $M_-(z,x_0,\alpha_0)$ replacing the
hypotheses concerning right Lebesgue points by those of left 
Lebesgue points, etc. For convenience we just state the corresponding
expansion and associated nonlinear recursion formula which covers both
cases. 
\begin{equation}
M_\pm (z,x,\alpha_0)\underset{\substack{\abs{z}\to\infty\\ z\in
C_\varepsilon}}{=}  \sum_{k=0}^N
m_{\pm,k}(x,\alpha_0)z^{-k}+ o\big(|z|^{-N}\big), \quad N\in\bbN.
\lb{4.100}
\end{equation}
\begin{align}
m_{\pm,0}(x,\alpha_0)&=\pm iI_m, \no \\
m_{\pm,1}(x,\alpha_0)&=-\f{1}{2} \big(
B_{1,2}(x)+B_{2,1}(x)\big)  \pm \f{i}{2} \big(
B_{1,1}(x)-B_{2,2}(x)\big), \no \\
m_{\pm,k+1}(x,\alpha_0)&=\pm\f{i}2\bigg(m_{\pm,k}^\prime(x,\alpha_0)+
\sum_{\ell=1}^{k}m_{\pm,\ell}(x,\alpha_0)
m_{\pm,k+1-\ell}(x,\alpha_0) \no \\
& \qquad \quad +\sum_{\ell=0}^{k}m_{\pm,\ell}(x,\alpha_0)B_{2,2}(x)
m_{\pm,k-\ell}(x,\alpha_0) \lb{4.101}  \\
& \qquad \quad +B_{1,2}(x)m_{\pm,k}(x,\alpha_0)
+m_{\pm,k}(x,\alpha_0)B_{2,1}(x)\bigg), \no \\
& \hspace*{5.8cm} 1\leq k\leq N-1. \no
\end{align} 
\end{remark}

Combining Theorem~\ref{t4.7} and \eqref{2.620} then yields the
analogous asymptotic expansion for $M(z,x,\alpha_0)$.

\begin{theorem} \lb{t4.10a}
Assume Hypothesis~\ref{h2.1} with $A=I_{2m}$, and let 
$\alpha_0=(I_m\; 0)\in\bbC^{m\times 2m}$. Fix $x_0\in\bbR$ and let
$x\in\bbR$. Suppose that for some $N\in\bbN$,
$B^{(N-1)}\in L^1_{\loc}(\bbR)^{2m\times 2m}$, that
$x$ is a right and a left Lebesgue point of $B^{(N-1)}$, and that for
all $R>0$,
\begin{align}
& \quad \underset{y\in [x_0,x_0+R]}{\esssup} \,
\bigg\|\int_y^{x_0+R} dx'\,B^{(N-1)}(x')\exp(2iz(x'-y))
+\f{1}{2iz}B^{(N-1)}(y)\bigg\|_{\bbC^{2m\times 2m}} \no \\
&+ \underset{y\in [x_0-R,x_0]}{\esssup} \,
\bigg\|\int_{x_0-R}^{y} dx'\,B^{(N-1)}(x')\exp(2iz(x'-y))
-\f{1}{2iz}B^{(N-1)}(y)\bigg\|_{\bbC^{2m\times 2m}} \no \\
& \underset{\substack{\abs{z}
\to\infty\\ z\in C_\varepsilon}}{=}\oh\big(|z|^{-1}\big). \lb{4.102}
\end{align}
If $N=1$, assume in addition 
$B_{k,k'}B_{\ell,\ell'}\in L^1_{\loc}(\bbR)^{m\times m}$
for all $k,k',\ell,\ell'\in\{1,2\}$.  Let $M(z,x,\alpha_0)$ be defined as
in \eqref{2.62} $($see also
\eqref{2.620}$)$.  Then, as $\abs{z}\to\infty$ in
$C_\varepsilon$, $M(z,x,\alpha_0)$ has an asymptotic expansion of the
form
\begin{equation}
M(z,x,\alpha_0)\underset{\substack{\abs{z}\to\infty\\ z\in
C_\varepsilon}}{=} (i/2) I_{2m} +\sum_{k=1}^N
M_{k}(x,\alpha_0)z^{-k}+ o\big(|z|^{-N}\big), \quad N\in\bbN, \lb{4.103}
\end{equation}
where
\begin{align}
M_1(x,\alpha_0)&=-\f{i}8 \begin{pmatrix}B_{1,1}(x+0)-B_{2,2}(x+0) &
B_{1,2}(x+0)+B_{2,1}(x+0) \\
 B_{1,2}(x+0)+B_{2,1}(x+0)& B_{2,2}(x+0)-B_{1,1}(x+0)\end{pmatrix} \no \\
& \quad -\f{i}8 \begin{pmatrix}B_{1,1}(x-0)-B_{2,2}(x-0) &
B_{1,2}(x-0)+B_{2,1}(x-0) \\
 B_{1,2}(x-0)+B_{2,1}(x-0)& B_{2,2}(x-0)-B_{1,1}(x-0)\end{pmatrix},
\text{ etc.} \lb{4.104}
\end{align}
The expansion \eqref{4.103} is uniform with respect to $\arg\,(z)$ for 
$|z|\to \infty$ in $C_\varepsilon$ and uniform in $x$ as long as $x$
varies in compact subsets of $\bbR$ intersected with the right and
left Lebesgue set of $B^{(N-1)}$. \\ 
\noindent If one merely assumes Hypothesis~\ref{h2.1} with $A=I_{2m}$,  
$\alpha_0=(I_m\; 0)$, and $B\in L^1_{\loc}(\bbR)^{2m\times 2m}$, then
\begin{equation}
M(z,x,\alpha_0)\underset{\substack{\abs{z}\to\infty\\ z\in
C_\varepsilon}}{=} (i/2)I_{2m}+\oh(1). \lb{4.105}
\end{equation}
Again the asymptotic expansion \eqref{4.105} is uniform with respect to
$\arg\,(z)$  for $|z|\to \infty$ in $C_\varepsilon$ and uniform in
$x\in\bbR$ as long as $x$ varies in compact intervals.
\end{theorem}
\noindent The higher-order coefficients in \eqref{4.103} can be derived
upon inserting \eqref{4.100} into \eqref{3.17a}, taking into account
\eqref{2.620}. 

Theorems~\ref{t4.6} and \ref{t4.7} (with $N\in\bbN$) are new even in the
scalar case $m=1$ with respect to the regularity assumptions on $B$. For
previous results in the case $m=1$ under stronger hypotheses on $B$ we
refer to \cite{EHS83}, \cite{Ha85}, \cite{HKS89a}, \cite{HKS89b}, 
\cite{Mi91}. In particular, \cite{Ha85}, \cite{HKS89a}, and \cite{HKS89b}
derived alternative high-energy expansions for the Weyl-Titchmarsh 
$m$-function in the case $m=1$.

Throughout this section we fixed $\alpha$ to be $\alpha_0=(I_m\; 0)$. The
case of general $\alpha\in\bbC^{2m\times m}$ satisfying \eqref{2.8e} then
follows from \eqref{2.41}.

\section{A Local Uniqueness Result} \lb{s5}

In this section we assume that $B$ is in the normal form given in 
Theorem~\ref{t1.1}, 
\begin{equation}\lb{4.115}
B(x)=\begin{pmatrix} B_{1,1}(x)
& B_{1,2}(x)\\B_{1,2}(x) & -B_{1,1}(x)  \end{pmatrix},
\end{equation}
with $B_{1,1}$ and $B_{1,2}$ self-adjoint a.e.
We prove fundamental new local uniqueness results  for $B$ in terms of 
exponentially small differences of Weyl-Titchmarsh 
matrices $M_+(z,x,\alpha)$ and $M(z,x,\alpha)$. These results, in turn, yield  
new global ramifications. We start with an auxiliary result concerning
asymptotic expansions.

\begin{lemma} \lb{l4.9} Suppose $\alpha=(\alpha_1 \;
\alpha_2)\in\bbC^{m\times 2m}$ satisfies \eqref{2.8e}, fix
$x_0,y_0\in\bbR$ with
$y_0>x_0$, and let
$x\geq x_0$. Assume $A=I_{2m}$, $B\in L^1([x_0,\infty))^{2m\times 2m}$, 
$\supp(B)\subseteq [x_0,y_0]$, with $B$ in the normal form given in 
\eqref{4.115}  a.e.~on
$(x_0,y_0)$. Then, the following asymptotic expansions hold for
$\Theta(z,x,x_0,\alpha)$, $\Phi(z,x,x_0,\alpha)$, and 
$U_{+}(z,x,x_0,\alpha)$ associated with \eqref{HSa},
\begin{align}
\Theta(z,x,x_0,\alpha)&\underset{\substack{\abs{z}
\to\infty\\ z\in \bbC_+}}{=}\f{1}{2}\begin{pmatrix}\alpha_1^*+i\alpha_2^*\\
-i(\alpha_1^*+i\alpha_2^*)\end{pmatrix}\exp(-iz(x-x_0))\big(1+\oh(1)\big),
\quad x>x_0,  \lb{4.48} \\
\Phi(z,x,x_0,\alpha)&\underset{\substack{\abs{z}
\to\infty\\ z\in \bbC_+}}{=}\f{i}{2}\begin{pmatrix}-\alpha_2^*+i\alpha_1^*\\
-i(-\alpha_2^*+i\alpha_1^*)\end{pmatrix}\exp(-iz(x-x_0))\big(1+\oh(1)\big),
\quad x>x_0, \lb{4.54} \\ 
U_{+}(z,x,x_0,\alpha)&\underset{\substack{\abs{z}
\to\infty\\ z\in \bbC_+}}{=}\begin{pmatrix}\alpha_1^*-i\alpha_2^*\\
i(\alpha_1^*-i\alpha_2^*)\end{pmatrix}\exp(iz(x-x_0))\big(1+\oh(1)\big),
\quad x\geq x_0. \lb{4.60} 
\end{align}
Next, we introduce the abbreviation
\begin{equation}
C=-B_{1,2}-iB_{1,1}, \quad C^*=-B_{1,2}+iB_{1,1}, \lb{4.38f}
\end{equation}
and suppose in addition that   
\begin{equation}
\underset{y\in [x_0,y_0]}{\esssup} \, \bigg\|\int_y^{y_0}
dx'\,B(x')\exp(2iz(x'-y)) +\f{1}{2iz}B(y)\bigg\|_{\bbC^{2m\times 2m}}
\underset{\substack{\abs{z}
\to\infty\\ z\in\rho_+}}{=}\oh\big(|z|^{-1}\big), \lb{4.39a} 
\end{equation}
along a ray $\rho_+\subset\bbC_+$, and that 
\begin{equation}
B_{1,1}^2,\ B_{1,2}^2,\ B_{1,1}B_{1,2},\ B_{1,2}B_{1,1} \in
L^1([x_0,y_0])^{m\times m}. \lb{4.39f}
\end{equation}
Then, 
\begin{align}
&\Theta(z,x,x_0,\alpha)\underset{\substack{\abs{z}
\to\infty\\ z\in
\rho_+}}{=}\bigg(\f{1}{2}\begin{pmatrix}\alpha_1^*+i\alpha_2^*\\
-i(\alpha_1^*+i\alpha_2^*)
\end{pmatrix}-\f{i}{4z}\begin{pmatrix}C(x_0)^*(\alpha_1^*-i\alpha_2^*)\\
-iC(x_0)^*(\alpha_1^*-i\alpha_2^*)\end{pmatrix} \no \\
& -\f{i}{4z}\begin{pmatrix}C(x)(\alpha_1^*+i\alpha_2^*)\\
iC(x)(\alpha_1^*+i\alpha_2^*)\end{pmatrix} \no \\  
& +\f{i}{4z}\int_{x_0}^x dx' \,
\begin{pmatrix}C(x')^*C(x')(\alpha_1^*+i\alpha_2^*)\\
-iC(x')^*C(x')(\alpha_1^*+i\alpha_2^*)\end{pmatrix}\bigg)
e^{-iz(x-x_0)}\big(1+\oh\big(|z|^{-1}\big)\big), \no \\
& \hspace*{9.5cm} x>x_0,  \lb{4.49} \\
&\Phi(z,x,x_0,\alpha)\underset{\substack{\abs{z}
\to\infty\\ z\in \rho_+}}{=}\bigg(\f{i}{2}\begin{pmatrix}
-\alpha_2^*+i\alpha_1^*\\
-i(-\alpha_2^*+i\alpha_1^*)\end{pmatrix}-\f{1}{4z}
\begin{pmatrix}C(x_0)^*(-\alpha_2^*-i\alpha_1^*)\\ 
-iC(x_0)^*(-\alpha_2^*-i\alpha_1^*)\end{pmatrix} \no \\
&+\f{1}{4z}\begin{pmatrix}C(x)(-\alpha_2^*+i\alpha_1^*)\\
iC(x)(-\alpha_2^*+i\alpha_1^*)\end{pmatrix}
\no \\  
& -\f{1}{4z}\int_{x_0}^x dx' \,
\begin{pmatrix}C(x')^*C(x')(-\alpha_2^*+i\alpha_1^*)\\
-iC(x')^*C(x')(-\alpha_2^*+i\alpha_1^*)\end{pmatrix}\bigg)
e^{-iz(x-x_0)}\big(1+\oh\big(|z|^{-1}\big)\big), \no \\
& \hspace*{9.8cm} x>x_0,  \lb{4.55}
\end{align} 
whenever $x_0$ is a right Lebesgue point of $B$ and $x$ is a left
Lebesgue point of $B$, and 
\begin{align}
&U_{+}(z,x,x_0,\alpha)\underset{\substack{\abs{z}
\to\infty\\ z\in \rho_+}}{=}\bigg(\begin{pmatrix} \alpha_1^*-i\alpha_2^*\\
i(\alpha_1^*-i\alpha_2^*)\end{pmatrix}+\f{i}{2z}
\begin{pmatrix}(C(x)^*-C(x_0)^*)(\alpha_1^*-i\alpha_2^*)\\
-i(C(x)^*+C(x_0)^*)(\alpha_1^*-i\alpha_2^*)\end{pmatrix}  \no \\  
& -\f{i}{2z}\int_{x_0}^{x} dx' \,
\begin{pmatrix}C(x')C(x')^*(\alpha_1^*-i\alpha_2^*)\\
iC(x')C(x')^*(\alpha_1^*-i\alpha_2^*) \end{pmatrix}\bigg)
e^{iz(x-x_0)}\big(1+\oh\big(|z|^{-1}\big)\big), \quad 
x\geq x_0, \lb{4.61}
\end{align}
whenever $x$ is a right Lebesgue point of $B$.
\end{lemma}
\begin{proof}
Since $x_0$ and $\alpha$ are fixed throughout this proof, we will
temporarily suppress these variables whenever possible to simplify
notations. Introducing 
\begin{equation}
\what \Theta(z,x)= 2\Theta(z,x)\exp(iz(x-x_0)), \lb{4.41}
\end{equation}
the Volterra integral equation for $\Theta$ (cf. \eqref{4.3B}), 
\begin{align}
\Theta(z,x)&=\begin{pmatrix}\alpha_1^*\cos(z(x-x_0))
+\alpha_2^*\sin(z(x-x_0))\\
\alpha_2^*\cos(z(x-x_0))-\alpha_1^*\sin(z(x-x_0))\end{pmatrix} \no \\
& \quad -\int_{x_0}^x dx'\, K(z,x,x') J B(x')\Theta(z,x'), \lb{4.42} 
\end{align}
can be rewritten in terms of that of $\what \Theta$ in the form
\begin{align}
\what\Theta(z,x)&=\begin{pmatrix}\alpha_1^*+i\alpha_2^*\\
-i(\alpha_1^*+i\alpha_2^*)\end{pmatrix}
+\begin{pmatrix}\alpha_1^*-i\alpha_2^*\\
i(\alpha_1^*-i\alpha_2^*)\end{pmatrix}\exp(2iz(x-x_0)) \no \\
& \quad -\f{1}{2}\int_{x_0}^x dx'
\big(R(x')\exp(2iz(x-x'))+S(x')\big)\what\Theta(z,x'), \lb{4.43}
\end{align} 
where we abbreviated
\begin{equation}
R=\begin{pmatrix} C&iC\\ iC &-C \end{pmatrix}, \quad 
S=\begin{pmatrix} C^* & -iC^*\\ -iC^* & -C^*\end{pmatrix}. \lb{4.44}
\end{equation}
Using the elementary algebraic facts
\begin{equation}
R\begin{pmatrix}a\\ ia\end{pmatrix}=0, \quad 
R\begin{pmatrix}b\\ -ib\end{pmatrix}
=2\begin{pmatrix} Cb\\ iCb\end{pmatrix}, \quad 
S\begin{pmatrix} a\\ ia\end{pmatrix}
=2\begin{pmatrix} C^* a\\ -iC^* a\end{pmatrix}, \quad 
S\begin{pmatrix}b\\ -ib\end{pmatrix}=0  \lb{4.46} 
\end{equation}
for any $a,b \in \bbC^{m\times m}$, iterating \eqref{4.43} yields
\begin{align}
&\what\Theta(z,x)=\begin{pmatrix}\alpha_1^*+i\alpha_2^*\\
-i(\alpha_1^*+i\alpha_2^*)\end{pmatrix}+
\begin{pmatrix}\alpha_1^*-i\alpha_2^*\\
i(\alpha_1^*+i\alpha_2^*)\end{pmatrix}e^{2iz(x-x_0)} \no \\
& \quad +\sum_{m=1}^\infty (-2)^{-m}\int_{x_0}^x d\xi_1 \,
\big(R(\xi_1)e^{2iz(x-\xi_1)}+S(\xi_1)\big)\times \no \\
& \quad \times \int_{x_0}^{\xi_1} d\xi_2 \,
\big(R(\xi_2)e^{2iz(\xi_1-\xi_2)}+S(\xi_2)\big)\dots \lb{4.47} \\
& \quad \;\, \dots \int_{x_0}^{\xi_{m-2}}
d\xi_{m-1} \,
\big(R(\xi_{m-1})e^{2iz(\xi_{m-2}-\xi_{m-1})}+S(\xi_{m-1})\big)
\times \no \\ 
& \quad \times \int_{x_0}^{\xi_{m-1}} d\xi_{m} \,
\bigg(R(\xi_{m})\begin{pmatrix}\alpha_1^*-i\alpha_2^*\\
-i(\alpha_1^*-i\alpha_2^*)\end{pmatrix}e^{2iz(\xi_{m-1}-\xi_{m})} 
\no \\
& \hspace*{2.8cm} +S(\xi_{m}) \begin{pmatrix}\alpha_1^*+i\alpha_2^*\\
i(\alpha_1^*+i\alpha_2^*)\end{pmatrix}e^{2iz(\xi_m-x_0)}\bigg).  \no
\end{align}
Applying the Riemann-Lebesgue lemma to \eqref{4.47} then proves
\eqref{4.48} assuming $B\in L^1([x_0,\infty))^{2m\times 2m}$, 
only. Assuming also \eqref{4.39a} and \eqref{4.39f} one can compute the
next term in the asymptotic expansion \eqref{4.48} and then obtains
\eqref{4.49} using \eqref{4.47} and the finite-interval variant of
\eqref{4.-2}, whenever $x_0$ is a right Lebesgue point of $B$ and $x$ is a
left Lebesgue point of $B$. \\ 
Exactly the same arguments apply to $\Phi$. Introducing  
\begin{equation}
\what \Phi(z,x)= 2\Phi(z,x)\exp(iz(x-x_0)), \lb{4.50}
\end{equation}
the Volterra integral equation for $\Phi$, 
\begin{align}
\Phi(z,x)&=\begin{pmatrix}-\alpha_2^*\cos(z(x-x_0))
+\alpha_1^*\sin(iz(x-x_0))\\
\alpha_1^*\cos(iz(x-x_0))+\alpha_2^*\sin(z(x-x_0))\end{pmatrix} \no \\
& \quad -\int_{x_0}^x dx'\, K(z,x,x') J B(x')\Phi_j(z,x'), \lb{4.51} 
\end{align}
can be rewritten in terms of that of $\what \Phi$ in the form
\begin{align}
\what\Phi(z,x)&=i\begin{pmatrix}-\alpha_2^*+i\alpha_1^*\\
-i(-\alpha_2^*+i\alpha_1^*)\end{pmatrix}
-i\begin{pmatrix}-\alpha_2^*-i\alpha_1^*\\
i(-\alpha_2^*-i\alpha_1^*)\end{pmatrix}\exp(2iz(x-x_0)) \no \\
& \quad -\f{1}{2}\int_{x_0}^x
\big(R(x')\exp(2iz(x-x'))+S(x')\big)\what\Phi(z,x'). \lb{4.52}
\end{align} 
Iterating \eqref{4.52}, taking into account \eqref{4.46}, yields
\begin{align}
&\what\Phi(z,x)=i\begin{pmatrix}-\alpha_2^*+i\alpha_1^*\\
-i(-\alpha_2^*+i\alpha_1^*)\end{pmatrix}
-i\begin{pmatrix}-\alpha_2^*-i\alpha_1^*\\
i(-\alpha_2^*-i\alpha_1^*)\end{pmatrix}e^{2iz(x-x_0)} \no \\
& \quad +\sum_{m=1}^\infty (-2)^{-m}\int_{x_0}^x d\xi_1 \,
\big(R(\xi_1)e^{2iz(x-\xi_1)}+S(\xi_1)\big)\times \no \\
& \quad \times \int_{x_0}^{\xi_1} d\xi_2 \,
\big(R(\xi_2)e^{2iz(\xi_1-\xi_2)}+S(\xi_2)\big)\dots \lb{4.53} \\
& \quad \;\, \dots \int_{x_0}^{\xi_{m-2}}
d\xi_{m-1} \,
\big(R(\xi_{m-1})e^{2iz(\xi_{m-2}-\xi_{m-1})}+S(\xi_{m-1})\big)
\times \no \\ 
& \quad \times \int_{x_0}^{\xi_{m-1}} d\xi_{m} \,
\bigg(iR(\xi_{m})\begin{pmatrix}-\alpha_2^*+i\alpha_1^*\\
-i(-\alpha_2^*+i\alpha_1^*)\end{pmatrix}
e^{2iz(\xi_{m-1}-\xi_{m})} \no \\
& \hspace*{2.8cm} -iS(\xi_{m})
\begin{pmatrix}-\alpha_2^*-i\alpha_1^*\\ 
i(-\alpha_2^*-i\alpha_1^*)\end{pmatrix}e^{2iz(\xi_m-x_0)}\bigg). \no
\end{align}
Applying the Riemann-Lebesgue lemma to \eqref{4.53} the proves
\eqref{4.54} assuming $B\in L^1([x_0,\infty))^{2m\times 2m}$, 
only. Assuming also \eqref{4.39a} and \eqref{4.39f} one can compute the
next term in the asymptotic expansion \eqref{4.54} and then obtains 
\eqref{4.55} using \eqref{4.53} and the finite-interval variant of
\eqref{4.-2}, whenever $x_0$ is a right Lebesgue point of $B$ and $x$ is
a left Lebesgue point of $B$. \\ 
Finally, we turn to $U_{+}(z,x)$. Introducing  
\begin{equation}
{\wti V}_{+}(z,x)= \wti U_{+}(z,x)\exp(-iz(x-x_0)), \lb{4.56}
\end{equation}
the Volterra integral equation for $\wti U_{+}$, 
\begin{equation}
\wti U_{+}(z,x)=\begin{pmatrix} \alpha_1^*-i\alpha_2^*\\
i(\alpha_1^*-i\alpha_2^*)\end{pmatrix}\exp(iz(x-x_0)) +\int_x^\infty
dx'\, K(z,x,x') J B(x')\wti U_{+}(z,x'),  \lb{4.57} 
\end{equation}
can be rewritten in terms of that of ${\wti V}_{+}$ in the form
\begin{equation}
{\wti V}_{+}(z,x)=\begin{pmatrix} \alpha_1^*-i\alpha_2^*\\
i(\alpha_1^*-i\alpha_2^*)\end{pmatrix} 
 +\f{1}{2}\int_{x}^{y_0} dx'\,
\big(R(x')+S(x')\exp(2iz(x'-x))\big){\wti V}_{+}(z,x'). \lb{4.58}
\end{equation} 
Iterating \eqref{4.58}, taking into account \eqref{4.46}, yields
\begin{align}
{\wti V}_{+}(z,x)&=\begin{pmatrix} \alpha_1^*-i\alpha_2^*\\
i(\alpha_1^*-i\alpha_2^*)\end{pmatrix}+\sum_{k=1}^\infty
2^{-2k}\int_{x}^{y_0} d\xi_1 \, R(\xi_1) \int_{\xi_1}^{y_0} d\xi_2 \,
S(\xi_2)e^{2iz(\xi_2-\xi_1)}\times
\no \\ 
& \quad \;\, \times \int_{\xi_2}^{y_0} d\xi_3 \,
R(\xi_3) \dots \int_{\xi_{2k-2}}^{y_0} d\xi_{2k-1} \, R(\xi_{2k-1})
\times \no \\ 
& \quad \times \int_{\xi_{2k-1}}^{y_0} d\xi_{2k} \,
S(\xi_{m})\begin{pmatrix} \alpha_1^*-i\alpha_2^*\\
i(\alpha_1^*-i\alpha_2^*)\end{pmatrix}e^{iz(\xi_{2k}-\xi_{2k-1})} \no \\
& \quad +\sum_{\ell=0}^\infty 2^{-2\ell+1}
\int_x^{y_0} d\xi_1 \, S(\xi_1)e^{2iz(\xi_1-x)} \int_{\xi_1}^{y_0} d\xi_2
\, R(\xi_2)\times \no \\
& \quad \;\, \times \int_{\xi_2}^{y_0} d\xi_3 \,
S(\xi_3)e^{2iz(\xi_3-\xi_2)} \dots \int_{\xi_{2\ell-1}}^{y_0}
d\xi_{2\ell} \, R(\xi_{2\ell}) \times \no \\ 
& \quad \times \int_{\xi_{2\ell}}^{y_0} d\xi_{2\ell+1} \,
S(\xi_{2\ell+1})\begin{pmatrix} \alpha_1^*-i\alpha_2^*\\
i(\alpha_1^*-i\alpha_2^*)\end{pmatrix}e^{iz(\xi_{2\ell+1}-\xi_{2\ell})}.
\lb{4.59}
\end{align}
Next, we take into account the different normalizations of $U_+$ and 
$\wti U_+$. Using $U_+(z,x_0)=[I_m \; M_+(z,x_0)^t]^t$ (cf., \eqref{2.52}
and $\Psi(z,x_0,x_0,\alpha_0)=I_{2m}$), one readily verifies the
relationship 
\begin{equation}
u_{+,1}(z,x)={\wti u}_{+,1}(z,x){\wti u}_{+,1}(z,x_0)^{-1}, 
\quad 
u_{+,2}(z,x)={\wti u}_{+,2}(z,x){\wti u}_{+,1}(z,x_0)^{-1}. 
\lb{4.59a}
\end{equation}
Thus, applying the Riemann-Lebesgue lemma to \eqref{4.59} then proves
\eqref{4.60} (in agreement with \eqref{4.14H}), assuming $B\in
L^1([x_0,\infty))^{2m\times 2m}$, only. Assuming also
\eqref{4.39a} and \eqref{4.39f} one can compute the next term in the
asymptotic expansion \eqref{4.60} and then obtains \eqref{4.61} using
\eqref{4.59} and \eqref{4.-2}, whenever $x$ is a right Lebesgue point of
$B$. 
\end{proof}

In the special case $m=1$ (and for $\alpha=(1\; 0)$), the expansion 
\eqref{4.61} was stated in \cite{Gr92}.  

Next, we note an elementary result concerning the boundary data
independence of exponentially close  Weyl-Titchmarsh matrices.

\begin{lemma} \lb{l4.9a}
Fix $x_0\in\bbR$ and suppose $A_j=I_{2m}$, $B_j=B_j^*\in
L^1([x_0,x_0+R])^{2m\times 2m}$ for all $R>0$. Denote by
$M_{+,j}(z,x,\alpha)$, $x\geq x_0$, the unique Weyl-Titchmarsh matrices
corresponding to the half-line Dirac-type operators $D_{+,j}(\alpha)$,
$j=1,2$, in \eqref{2.84}.  Fix an $\hat\alpha\in\bbC^{m\times 2m}$
satisfying
\eqref{2.8e} and assume that for all $\varepsilon >0$, 
\begin{equation}
\|M_{+,1}(z,x_0,\hat\alpha)-M_{+,2}(z,x_0,\hat\alpha)\|_{\bbC^{m\times
m}}\underset{\substack{|z|\to\infty\\z\in \rho_{+}}}{=}
O\big(e^{-2\Im(z)(a-\varepsilon)}\big) \lb{4.24a}
\end{equation}
along some ray $\rho_{+}\subset\bbC_+$. Then, for all 
$\alpha\in\bbC^{m\times 2m}$ satisfying \eqref{2.8e} and for all
$\varepsilon >0$, 
\begin{equation}
\|M_{+,1}(z,x_0,\alpha)-M_{+,2}(z,x_0,\alpha)\|_{\bbC^{m\times
m}}\underset{\substack{|z|\to\infty\\z\in \rho_{+}}}{=}
O\big(e^{-2\Im(z)(a-\varepsilon)}\big) \lb{4.25aa}
\end{equation}
along the ray $\rho_{+}$.
\end{lemma}
\begin{proof}
Using \eqref{2.38} and \eqref{2.41} one estimates
\begin{align}
&\|M_{+,1}(z,x_0,\alpha)-M_{+,2}(z,x_0,\alpha)\|_{\bbC^{m\times m}} \no \\
&=\|M_{+,1}(\bar z,x_0,\alpha)^* -M_{+,2}(z,x_0,\alpha)\|_{\bbC^{m\times
m}} \no \\ 
&\leq \|[\hat\alpha\alpha^* -M_{+,1}(\bar z,x_0,\hat\alpha)^*\hat\alpha
 J\alpha^*]^{-1}\|_{\bbC^{m\times m}} \times \no \\
& \quad \times
\|M_{+,1}(z,x_0,\hat\alpha)-M_{+,2}(z,x_0,\hat\alpha)\|_{\bbC^{m\times
m}}\times \no \\ & \quad \times \|[\alpha\hat\alpha^* +\alpha J\hat\alpha^*
 M_{+,2}(z,x_0,\hat\alpha)]^{-1}\|_{\bbC^{m\times m}},
\lb{4.26aa}
\end{align}
since by \eqref{2.8i}
\begin{equation}
\hat\alpha J\alpha^*\alpha\hat\alpha^*+\hat\alpha\alpha^*\alpha
J\hat\alpha^*=0, \quad \hat\alpha\alpha^*\alpha\hat\alpha^*- \hat\alpha
J\alpha^*\alpha J\hat\alpha^*=I_m. \lb{4.27a}
\end{equation}
Moreover, since
\begin{equation}
[\hat\alpha\alpha^*-i\hat\alpha J\alpha^*][\hat\alpha\alpha^*
-i\hat\alpha J\alpha^*]^*=I_m, \lb{5.29a}
\end{equation}
by \eqref{4.27a}, one infers \eqref{4.25aa} from 
\eqref{4.26aa} and $M_{+,j}(z,x_0,\alpha)=iI_m +\oh(1)$ as $|z|\to\infty$,
$z\in\bbC_+$, $j=1,2$ (cf.~\eqref{3.1}).
\end{proof}

Our principal new local uniqueness result for Dirac-type operators in
terms of Weyl-Titchmarsh matrices then reads as follows. 

\begin{theorem} \lb{t4.10} 
Fix $x_0\in\bbR$ and suppose $A_j=I_{2m}$, $B_j\in
L^1([x_0,x_0+R])^{2m\times 2m}$ for all $R>0$. Suppose also that $B_j$
is in the normal form given in \eqref{4.115} a.e.~on $(x_0,\infty)$,
$j=1,2$. Denote by $M_{j,+}(z,x,\alpha)$, $x\geq
x_0$, the unique Weyl-Titchmarsh matrices corresponding to the half-line
Dirac-type operators $D_{+,j}(\alpha)$, $j=1,2$, in \eqref{2.84}. Then,
\begin{equation}
\text{if for some $a>0$, }\, B_1(x)=B_2(x) \, \text{ for a.e. $x\in
(x_0,x_0+a)$,} \lb{4.38AA} 
\end{equation}
one obtains 
\begin{equation}
\|M_{1,+}(z,x_0,\alpha)-M_{2,+}(z,x_0,\alpha)\|_{\bbC^{m\times
m}}\underset{\substack{\abs{z}
\to\infty\\ z\in \rho_{+}}}{=}
O\big(e^{-2\Im(z)a}\big) \lb{4.38}
\end{equation}
along any ray $\rho_+\subset\bbC_+$ with $0<\arg(z)<\pi$ and for all
$\alpha\in\bbC^{m\times 2m}$ satisfying \eqref{2.8e}. Conversely,
fix an $\hat\alpha\in\bbC^{m\times 2m}$ satisfying
\eqref{2.8e} and if $m>1$, assume in addition that
$B_j\in L^\infty([x_0,x_0+a])^{2m\times 2m}$, $j=1,2$. Moreover, suppose
that for all $\varepsilon >0$, 
\begin{equation}
\|M_{1,+}(z,x_0,\hat\alpha_1)-M_{2,+}(z,x_0,
\hat\alpha_1)\|_{\bbC^{m\times
m}}\underset{\substack{|z|\to\infty\\z\in \rho_{+,\ell}}}{=}
O\big(e^{-2\Im(z)(a-\varepsilon)}\big), \quad \ell=1,2, \lb{4.39}
\end{equation}
along a ray $\rho_{+,1}\subset\bbC_+$ with $0<\arg(z)<\pi/2$ and along a
ray $\rho_{+,2}\subset\bbC_+$ with $\pi/2<\arg(z)<\pi$. Then 
\begin{equation}
B_1(x)=B_2(x) \text{ for a.e. } x\in [x_0,x_0+a].  \lb{4.40}
\end{equation}
\end{theorem}
\begin{proof}
Since \eqref{4.38} follows from Theorem~\ref{t4.5} and Lemma~\ref{l4.9a},
it suffices to focus on the proof of \eqref{4.40}. Moreover, applying
Theorem~\ref{t4.5}, we may without loss of generality assume for the rest
of the proof that 
\begin{equation}
\supp(B_j)\subseteq [x_0,x_0+a], \quad j=1,2. \lb{4.40a}
\end{equation}
In the following, we will adapt the principal ingredients of a recent
proof of the local Borg-Marchenko uniqueness theorem for scalar
Schr\"odinger operators (i.e., for $m=1$) by Bennewitz \cite{Be00}, to the
current Dirac-type situation. First we recall that by
Lemma~\ref{l4.9a}, \eqref{4.39} holds along the rays $\rho_{+,j}$,
$j=1,2$ for all $\alpha=(\alpha_1\; \alpha_2)\in\bbC^{m\times 2m}$
satisfying \eqref{2.8e}. To simplify notations in the following
we will again suppress $x_0$ and $\alpha$ whenever possible and hence 
abbreviate, $\Theta(z,x,x_0,\alpha)$, $\Phi(z,x,x_0,\alpha)$, and
$U_{j,+}(z,x,x_0,\alpha)$ by $\Theta(z,x)$, $\Phi(z,x)$, and
$U_{j,+}(z,x)$, respectively. Next, denoting in obvious notation by 
\begin{align}
&\Theta_j(z,x)=\begin{pmatrix}\theta_{j,1}(z,x)\\
\theta_{j,2}(z,x)\end{pmatrix}, \quad
\Phi_j(z,x)=\begin{pmatrix}\phi_{j,1}(z,x)\\
\phi_{j,2}(z,x)\end{pmatrix}, \quad
U_{j,+}(z,x)=\begin{pmatrix} u_{j,+,1}(z,x)\\
u_{j,+,2}(z,x)\end{pmatrix}, \no \\
& \hspace*{8.5cm} j=1,2, \; x\geq x_0, \lb{4.62}
\end{align}
the solutions  associated with $B_j$, $j=1,2$, which are defined in 
\eqref{FSb} and \eqref{2.14}, we introduce 
\begin{equation}
g_{j,k}(z,x)=\phi_{j,k}(z,x) u_{j,+,k}(\bar z,x)^*, \quad j,k
\in\{1,2\}, \; x\geq x_0. \lb{4.63}
\end{equation} 
Using the asymptotic expansions \eqref{4.48}--\eqref{4.60} for 
$\Theta_j(z,x)$, $\Phi_j(z,x)$, and $U_{j,+}(z,x)$, and the analogous
ones for $\Theta_j(\bar z,x)^*$, $\Phi_j(\bar z,x)^*$, and $U_{j,+}(\bar
z,x)^*$, one  verifies for each fixed $x>x_0$,
\begin{equation}
g_{j,k}(z,x) \underset{\substack{\abs{z}
\to\infty\\ z\in \bbC_+}}{=}(i/2)I_m + \oh(1), \quad j,k \in\{1,2\},
\lb{4.65} 
\end{equation}
assuming $B_j\in L^1([x_0,x_0+R])^{2m\times 2m}$ for all $R>0$, $j=1,2$,
only.  Next, using the fact that for each fixed $x>x_0$,
\begin{align}
\phi_{1,k}(z,x)^{-1}\phi_{2,k}(z,x)&\underset{\substack{\abs{z}
\to\infty\\ z\in \bbC_+}}{=}I_m + \oh(1), \quad k=1,2, \lb{4.65a} \\ 
(u_{1,+,k}(\bar z,x)^*)^{-1} u_{2,+,k}(\bar
z,x)^*&\underset{\substack{\abs{z}
\to\infty\\ z\in \bbC_+}}{=}I_m + \oh(1), \quad k=1,2, \lb{4.66}
\end{align}
by \eqref{4.54}, \eqref{4.60}, one concludes 
\begin{align}
&\phi_{1,k}(z,x) u_{2,+,j}(\bar z,x)^*-u_{1,+,k}(z,x)
\phi_{2,k}(\bar z,x)^* \no \\
&=\phi_{1,k}(z,x) \theta_{2,k}(\bar z,x)^*-\theta_{1,k}(z,x)
\phi_{2,k}(\bar z,x)^* \no \\
& \quad +\phi_{1,k}(z,x) \big(M_{2,+}(z)-M_{1,+}(z) \big)
\phi_{2,k}(\bar z,x)^* \underset{\substack{\abs{z}\to\infty\\ z\in
\bbC_+}}{=}\oh(1), \lb{4.67}
\end{align}
using \eqref{4.65}, \eqref{4.66}, and 
$M_{2,+}(\bar z)^*=M_{2,+}(z)$. Combining hypothesis \eqref{4.39} and 
\eqref{4.54}, one infers
\begin{equation}
\big\|\phi_{1,k}(z,x) \big(M_{2,+}(z)-M_{1,+}(z) \big)
\phi_{2,k}(\bar z,x)^*\big\| \underset{\substack{\abs{z}\to\infty\\ z\in
\rho_{+,\ell}}}{=}\oh(1), \ \ x\in (x_0,x_0+a) \lb{4.68}
\end{equation}
along the rays $\rho_{+,\ell}$, $\ell=1,2$. Thus, \eqref{4.67} implies
\begin{equation}
\big\|\phi_{1,k}(z,x) \theta_{2,k}(\bar z,x)^*-\theta_{1,k}(z,x)
\phi_{2,k}(\bar z,x)^*\big\| \underset{\substack{\abs{z}\to\infty\\ z\in
\rho_{+,\ell}}}{=}\oh(1), \ \ x\in (x_0,x_0+a) \lb{4.69}
\end{equation}
along the rays $\rho_{+,\ell}$, $\ell=1,2$. The analogous estimate
\eqref{4.69} holds along  the complex conjugate rays $\bar \rho_{+,\ell}$,
$\ell=1,2$, in the lower complex half-plane $\bbC_-$. To simplify
notations we denote the open sector generated by $\rho_{+,1}$ and its
complex conjugate $\bar \rho_{+,1}$ by $\calS_1$, the open sector 
generated by the $\rho_{+,2}$ and its complex conjugate $\bar \rho_{+,2}$
by
$\calS_2$, the remaining sector in $\bbC_+$ is denoted by $\calS_3$, and
its complex conjugate sector in $\bbC_-$ is denoted by $\calS_4$. Thus,
one obtains a partition of $\bbC$ into
\begin{equation}
\bbC=\bigcup_{\ell=1}^4 \ol {\calS_\ell}, \lb{4.70}
\end{equation}
where each sector $\calS_\ell$, $1\leq \ell\leq 4$, has opening angle
strictly less than $\pi$. Since (each matrix element of) the expression
under the norm in \eqref{4.69} is entire and of order less or equal to
one, one can apply the Phragm\'en-Lindel\"of principle (cf., e.g.,
\cite[No.~322, p.~166--167, 379]{PS72}) to each sector $\calS_\ell$,
$1\leq \ell\leq 4$, and obtains that each matrix element under the norm in
\eqref{4.69} is uniformly bounded in each sector and hence on all of
$\bbC$. By Liouville's theorem, these matrix elements are all equal to
certain constants. By the right-hand side of \eqref{4.69}, these
constants all vanish. Thus, we proved
\begin{equation}
\phi_{1,k}(z,x) \theta_{2,k}(\bar z,x)^*=\theta_{1,k}(z,x)
\phi_{2,k}(\bar z,x)^* \, \text{ for all $x\in (x_0,x_0+a)$}
\lb{4.71}
\end{equation}
and hence
\begin{equation}
\phi_{1,k}(z,x)^{-1}\theta_{1,k}(z,x)=\theta_{2,k}(\bar z,x)^*
(\phi_{2,k}(\bar z,x)^*)^{-1} \, \text{ for all $x\in (x_0,x_0+a)$.}
\lb{4.72}
\end{equation}
Differentiating $\phi_{j,k}(z,x)^{-1}\theta_{j,k}(z,x)$, $j,k=1,2$, with
respect to $x$ yields
\begin{align} 
&\big(\phi_{j,1}(z,x)^{-1}\theta_{j,1}(z,x)\big)' \no \\
&=\phi_{j,1}(z,x)^{-1}((B_j)_{1,1}(x)-z)(\phi_{j,2}(z,x)
\phi_{j,1}(z,x)^{-1}
\theta_{j,1}(z,x)-\theta_{j,2}(z,x)), \lb{4.73} \\
&\big(\phi_{j,2}(z,x)^{-1}\theta_{j,2}(z,x)\big)' \no \\
&=\phi_{j,2}(z,x)^{-1}((B_j)_{1,1}(x)+z)(\phi_{j,1}(z,x)
\phi_{j,2}(z,x)^{-1}\theta_{j,2}(z,x)-\theta_{j,1}(z,x)). \lb{4.74} 
\end{align}
Multiplying \eqref{4.73} by $\phi_{j,1}(\bar z,x)^*(\phi_{j,1}(\bar
z,x)^*)^{-1}$  and using \eqref{2.92}, \eqref{2.94}, and similarly,
multiplying
\eqref{4.74} by $\phi_{j,2}(\bar z,x)^*(\phi_{j,2}(\bar z,x)^*)^{-1}$ 
and using \eqref{2.93}, \eqref{2.95} then yields
\begin{align} 
\big(\phi_{j,1}(z,x)^{-1}\theta_{j,1}(z,x)\big)' &=
\phi_{j,1}(z,x)^{-1}((B_j)_{1,1}(x)-z)(\phi_{j,1}(\bar z,x)^*)^{-1},
\lb{4.75}
\\
\big(\phi_{j,2}(z,x)^{-1}\theta_{j,2}(z,x)\big)' &=
\phi_{j,2}(z,x)^{-1}((B_j)_{1,1}(x)+z)(\phi_{j,2}(\bar z,x)^*)^{-1}.
\lb{4.76} 
\end{align}
In exactly the same way one derives 
\begin{align} 
&\big(\theta_{j,1}(\bar z,x)^*(\phi_{j,1}(\bar z,x)^*)^{-1}\big)' \no \\
&=
(\theta_{j,1}(\bar z,x)^*(\phi_{j,1}(\bar z,x)^*)^{-1}\phi_{j,2}(\bar
z,x)^*-\theta_2(\bar z,x)^*)((B_j)_{1,1}(x)-z)(\phi_{j,1}(\bar
z,x)^*)^{-1} \no \\ 
&= \phi_{j,1}(z,x)^{-1}((B_j)_{1,1}(x)-z)(\phi_{j,1}(\bar
z,x)^*)^{-1}, \lb{4.77}
\\ &\big(\theta_{j,2}(\bar z,x)^*(\phi_{j,2}(\bar z,x)^*)^{-1}\big)' \no
\\ &=(\theta_{j,2}(\bar z,x)^*(\phi_{j,2}(\bar z,x)^*)^{-1}\phi_{j,1}(\bar
z,x)^*-\theta_{j,1}(\bar z,x)^*)((B_j)_{1,1}(x)+z)(\phi_{j,2}(\bar
z,x)^*)^{-1}  \no \\
&=\phi_{j,2}(z,x)^{-1}((B_j)_{1,1}(x)+z)(\phi_{j,2}(\bar z,x)^*)^{-1}, 
\lb{4.78} 
\end{align}
using \eqref{2.72}--\eqref{2.75}. Thus, \eqref{4.72} implies
\begin{align}
\phi_{1,1}(\bar z,x)^*((B_1)_{1,1}(x)-z)^{-1}\phi_{1,1}(z,x)&=
\phi_{2,1}(\bar z,x)^*((B_2)_{1,1}(x)-z)^{-1}\phi_{2,1}(z,x), \lb{4.79} \\
\phi_{1,2}(\bar z,x)^*((B_1)_{1,1}(x)+z)^{-1}\phi_{1,2}(z,x)&=
\phi_{2,2}(\bar z,x)^*((B_2)_{1,1}(x)+z)^{-1}\phi_{2,2}(z,x), \lb{4.80} \\
\theta_{1,1}(\bar z,x)^*((B_1)_{1,1}(x)-z)^{-1}\theta_{1,1}(z,x)&=
\theta_{2,1}(\bar z,x)^*((B_2)_{1,1}(x)-z)^{-1}\theta_{2,1}(z,x),
\lb{4.81}
\\
\theta_{1,2}(\bar z,x)^*((B_1)_{1,1}(x)+z)^{-1}\theta_{1,2}(z,x)&=
\theta_{2,2}(\bar z,x)^*((B_2)_{1,1}(x)+z)^{-1}\theta_{2,2}(z,x)
\lb{4.82} 
\end{align}
for a.e.~$x\in (x_0,x_0+a)$. Thus far we only used $B_j\in
L^1([x_0,x_0+R])^{2m\times 2m}$ for all $R>0$, $j=1,2$ and \eqref{4.40a}.

In the special case $m=1$, each of the equations
\eqref{4.79}--\eqref{4.82} allows for the completion of the proof of \eqref{4.40}.
Indeed, using the fact that
\begin{equation}
\overline{ \phi_{j,k}(\bar z,x)}=\phi_{j,k}(z,x), \quad 
\overline{ \theta_{j,k}(\bar z,x)}=\theta_{j,k}(z,x), \quad j,k\in\{1,2\}, 
\lb{4.83}
\end{equation}
and taking for instance 
\eqref{4.79}, one infers for a.e.~$x\in (x_0,x_0+a)$, that
\begin{equation}
\f{\phi_{1,1}(z,x)^2}{\phi_{2,1}(z,x)^2}=\f{(B_1)_{1,1}(x)-z}
{(B_2)_{1,1}(x)-z}. \lb{4.84}
\end{equation}
Since all zeros (and poles) of the left-hand side of \eqref{4.84} have
even multiplicity, while all zeros (and poles) of the right-hand side of
\eqref{4.83} are simple, one concludes, assuming only that $B_j\in
L^1([x_0,x_0+R])^{2\times 2}$ for all $R>0$, $j=1,2$, that
\begin{equation} \lb{4.84a}
(B_1)_{1,1}(x)=(B_2)_{1,1}(x) \, \text{  for a.e.~$x\in (x_0,x_0+a)$.}
\end{equation}   
Thus for the case $m=1$, we see by \eqref{4.79},
and \eqref{4.80}, \eqref{4.83}, and \eqref{4.84a},  for a.e.~$x\in (x_0,x_0+a)$, that
\begin{equation}\lb{4.84b}
\phi_{1,k}^2(z,x) = \phi_{2,k}^2(z,x), \ \ k=1,2.
\end{equation}
Now, \eqref{2.75}, \eqref{4.72}, and \eqref{4.83} show, for a.e.~$x\in (x_0,x_0+a)$,
that
\begin{align}\lb{4.84c}
\phi_{1,1}(z,x)\phi_{1,2}(z,x)&= \frac{\phi_{1,1}(z,x)}{\theta_{1,1}(z,x)}-
\frac{\phi_{1,2}(z,x)}{\theta_{1,2}(z,x)} \no \\
&=\frac{\phi_{2,1}(z,x)}{\theta_{2,1}(z,x)}-
\frac{\phi_{2,2}(z,x)}{\theta_{2,2}(z,x)}=\phi_{2,1}(z,x)\phi_{2,2}(z,x).
\end{align}
By \eqref{HSa} we see that 
\begin{equation}\lb{4.84d}
(\phi_{j,1}^2(z,x))' = 2(z-(B_j(x))_{1,1} )\phi_{j,1}(z,x)\phi_{j,2}(z,x)
 + (B_j(x))_{1,2} \phi_{j,1}^2(z,x), \quad j=1,2. 
\end{equation}
Thus, by \eqref{4.84a}, \eqref{4.84b}, and \eqref{4.84c}, 
\begin{equation}\lb{4.84e}
(B_1(x))_{1,2} = (B_2(x))_{1,2} \,\text{  for a.e.~$x\in (x_0,x_0+a)$.}
\end{equation}
Together, \eqref{4.84a} and \eqref{4.84e} imply \eqref{4.40} in the
special case $m=1$.

Unfortunately, the case $m>1$ appears to be quite a bit more
involved. To deal with this case we first note that taking determinants
in \eqref{4.79} yields
\begin{equation}
\f{\det(\phi_{1,1}(\bar z,x,x_0,\alpha)^*)\det(\phi_{1,1}(z,x,x_0,\alpha))}
{\det(\phi_{2,1}(\bar z,x,x_0,\alpha)^*)\det(\phi_{2,1}(z,x,x_0,\alpha)}=
\f{\det((B_1)_{1,1}(x)-zI_m))}{\det((B_2)_{1,1}(x)-zI_m))} \lb{4.85}
\end{equation}
for a.e.~$x\in (x_0,x_0+a)$.  
Next, we intend to prove that 
\begin{equation}
\det((B_1)_{1,1}(x)-zI_m))=\det((B_2)_{1,1}(x)-zI_m)) \, 
\text{  for a.e.~$x\in (x_0,x_0+a)$.} \lb{4.86}
\end{equation}
Given the fact that $(B_j)_{1,1}(x)$, $j=1,2$, is self-adjoint, showing
\eqref{4.86}  is equivalent to showing that $B_1(x)$ and $B_2(x)$ are
unitarily equivalent for a.e.~$x\in (x_0,x_0+a)$. Arguing by
contradiction, we assume that  at least one pair of eigenvalues of
$B_1(x)$ and $B_2(x)$ differs. Thus,  fixing $x_1\in (x_0,x_0+a)$, let
$\lambda(x_1)$ be an eigenvalue of
$B_1(x_1)$ but not of $B_2(x_1)$. Then 
\eqref{4.85} implies, for all $\alpha\in \bbC^{m\times 2m}$ satisfying
\eqref{2.8e}, that
\begin{equation}
\det(\phi_{1,1}(\lambda(x_1),x_1,x_0,\alpha))=0.  \lb{4.87}
\end{equation}
Next, for $\lambda\in\bbR$ and $x>x_0$ define 
\begin{align}\lb{4.88}
N(\lambda,x,\alpha)=&\theta_{1,2}(\lambda,x,x_0,\alpha)\theta_{1,2}
(\lambda,x,x_0,\alpha)^*  \no \\ 
&+\phi_{1,2}(\lambda,x,x_0,\alpha)\phi_{1,2}(\lambda,x,x_0,\alpha)^*. 
\end{align}
Then, $N(\lambda,x,\alpha)$ is strictly positive definite,
\begin{equation}
N(\lambda,x,\alpha)>0. \lb{4.89}
\end{equation}
Indeed, suppose $Nf=0$ for some $f\in\bbC^m$, then 
\begin{equation}
\theta_{1,2}(\lambda)\theta_{1,2}(\lambda)^*f+
\phi_{1,2}(\lambda)\phi_{1,2}(\lambda)^*f=0 \lb{4.90}
\end{equation}
implies  
\begin{equation}
\theta_{1,2}(\lambda)\theta_{1,2}(\lambda)^*f=0, 
\; \phi_{1,2}(\lambda)\phi_{1,2}(\lambda)^*f=0 \lb{4.91}
\end{equation}
and hence
\begin{equation}
\theta_{1,2}(\lambda)^*f=0, \; \phi_{1,2}(\lambda)^*f=0. \lb{4.92}
\end{equation}
Thus,
\begin{equation}
f=(\theta_{1,1}(\lambda)\phi_{1,2}(\lambda)^*-
\phi_{1,1}(\lambda)\theta_{1,2}(\lambda)^*)f=0 \lb{4.93} 
\end{equation}
by \eqref{2.95}, and hence $f=0$ proves \eqref{4.89}. Introducing
$\alpha_0=(I_m\; 0)\in\bbC^{m\times 2m}$ and $\gamma=(\gamma_1 \;
\gamma_2)\in\bbC^{m\times 2m}$ defined by 
\begin{align}
\gamma_1 &=[\theta_{1,2}(\lambda(x_1),x_1,x_0,\alpha_0)
\theta_{1,2}(\lambda(x_1),x_1,x_0,\alpha_0)^* \lb{4.94} \\
& \quad +\phi_{1,2}(\lambda(x_1),x_1,x_0,\alpha_0)
\phi_{1,2}(\lambda(x_1),x_1,x_0,\alpha_0)^*]^{-1/2}
\theta_{1,2}(\lambda(x_1),x_1,x_0,\alpha_0), \no \\
\gamma_2 &=[\theta_{1,2}(\lambda(x_1),x_1,x_0,\alpha_0)
\theta_{1,2}(\lambda(x_1),x_1,x_0,\alpha_0)^* \lb{4.95} \\
& \quad +\phi_{1,2}(\lambda(x_1),x_1,x_0,\alpha_0)
\phi_{1,2}(\lambda(x_1),x_1,x_0,\alpha_0)^*]^{-1/2}
\phi_{1,2}(\lambda(x_1),x_1,x_0,\alpha_0), \no
\end{align}
one verifies $\gamma\gamma^*=I_m$ (by \eqref{4.94} and \eqref{4.95}) and 
$\gamma J\gamma^*=0$ (by \eqref{2.93}). Thus, $\gamma$ satisfies 
\eqref{2.8e}. Next, since
\begin{equation}
\phi_{1,1}(\lambda(x_1),x_1,x_0,\gamma)=
\phi_{1,1}(\lambda(x_1),x_1,x_0,\alpha_0)\gamma_1^*-
\theta_{1,1}(\lambda(x_1),x_1,x_0,\alpha_0)\gamma_2^* \lb{4.96}
\end{equation}
as a special case of \eqref{2.96}, one derives
\begin{align}
\phi_{1,1}(\lambda(x_1),x_1,x_0,\gamma)&=
[\phi_{1,1}(\lambda(x_1),x_1,x_0,\alpha_0)
\theta_{1,2}(\lambda(x_1),x_1,x_0,\alpha_0)^* \no \\
& \quad -\theta_{1,1}(\lambda(x_1),x_1,x_0,\alpha_0)
\phi_{1,2}(\lambda(x_1),x_1,x_0,\alpha_0)^*]\times \no \\
& \quad \times [\theta_{1,2}(\lambda(x_1),x_1,x_0,\alpha_0)
\theta_{1,2}(\lambda(x_1),x_1,x_0,\alpha_0)^* \no \\
& \quad +\phi_{1,2}(\lambda(x_1),x_1,x_0,\alpha_0)
\phi_{1,2}(\lambda(x_1),x_1,x_0,\alpha_0)^*]^{-1/2} \no \\
& =-[\theta_{1,2}(\lambda(x_1),x_1,x_0,\alpha_0)
\theta_{1,2}(\lambda(x_1),x_1,x_0,\alpha_0)^* \lb{4.97} \\
& \quad +\phi_{1,2}(\lambda(x_1),x_1,x_0,\alpha_0)
\phi_{1,2}(\lambda(x_1),x_1,x_0,\alpha_0)^*]^{-1/2}<0. \no
\end{align}
using \eqref{2.95}. This contradiction to \eqref{4.87} proves
\eqref{4.86}.  Hence for $\lambda\in \bbR$ and for a.e.~$x\in (x_0,x_0+a)$
\begin{equation}
\abs{\det(\phi_{1,1}(\lambda,x,x_0,\alpha))}=
\abs{\det(\phi_{2,1}(\lambda,x,x_0,\alpha))},
\lb{4.98}
\end{equation}
by \eqref{4.85}. Equation \eqref{4.98} implies that for 
a.e.~$x_1\in(x_0,x_0+a)$, the family of Dirac
operators $D_+(\alpha,\alpha_0)$ in $L^2([x_0,x_1])^{2m}$, defined by 
\begin{align}
D_+(\alpha,\alpha_0)&=J \f{d}{dx}-B, \lb{4.99} \\
\dom(D_+(\alpha,\alpha_0))&=\{\phi\in L^2([x_0,x_1])^{2m} \mid \phi
\in\AC([x_0,x_1])^{2m}; \no \\
& \qquad \alpha\phi(x_0)=0, \, \alpha_0\phi(x_1)=0; \,
(J\phi^\prime-B\phi)\in  L^2([x_0,x_1])^{2m} \}, \no
\end{align}
with $\alpha_0 = (I_m\;0)$, have identical spectra for all boundary data 
$\alpha\in\bbC^{m\times 2m}$ satisfying \eqref{2.8e}. Hence, 
assuming $B_j\in L^\infty([x_0,x_0+a])^{2m\times 2m}$, $j=1,2$, one can
apply Theorem~2.3 of Malamud \cite{Ma99} and obtains \eqref{4.40}. 
\end{proof}

We should note that Malamud's Theorem~2.3 in \cite{Ma99} only requires
the equality of $m^2+1$ spectra (associated with linearly independent
boundary data indexed by $\alpha\in\bbC^{m\times 2m}$) in order to
conclude \eqref{4.40}. 

There is no particular significance of the rays $\rho_\ell$, 
$\ell=1,2$, in Theorem~\ref{t4.10}. Any non-selfintersecting Jordan 
arc that tends to infinity in the sectors $\varepsilon\leq\arg(z)\leq
(\pi/2)-\varepsilon$ and $(\pi/2)+\varepsilon\leq\arg(z)\leq
\pi-\varepsilon$ for some $0<\varepsilon<\pi/4$ will do.

\begin{remark} \lb{r4.10}
We were not able to prove \eqref{4.40} directly from
\eqref{4.79}--\eqref{4.82}, without resorting to the arguments involving
\eqref{4.98} and \eqref{4.99}. To conclude the proof according to the
Borg-type Theorem~2.3 of Malamud \cite{Ma99} (cf.~also Theorem~4 in
\cite{Ma99a}), requires the introduction of the extra hypothesis 
$B_j\in L^\infty([x_0,x_0+a])^{2m\times 2m}$, $j=1,2$ in
the matrix context $m>1$, since the construction of transformation
operators for Dirac-type systems, to date, uses such an additional
hypothesis on $B$. This extra hypothesis is clearly superfluous in the
case $m=1$. Obviously, one conjectures that this extra hypothesis on $B_j$
should also be redundant in Theorem~\ref{t4.10}, but this appears
to require nontrivial future efforts. In this context it might be
interesting to note that the higher-order expansions
\eqref{4.49}--\eqref{4.61} do not determine $B$ uniquely. An explicit
analysis shows that while they do determine $B_{1,2}$, they only determine
$B_{1,1}^2$, not $B_{1,1}$ itself. So that approach does not aide in
proving \eqref{4.40} (besides, it would require the additional hypotheses
\eqref{4.39a} on $B$).
\end{remark}

The corresponding local uniqueness result in terms of $M(z,x_0,\alpha)$
then reads as follows. 

\begin{theorem} \lb{t4.11}
Fix $x_0\in\bbR$ and suppose $A_j=I_{2m}$, $B_j\in
L^1_{\loc}(\bbR)^{2m\times 2m}$, and $B_j=B_j^*$
a.e.~on $\bbR$, $j=1,2$. Suppose also that $B_j$
is in the normal form given in \eqref{4.115} a.e.~on $(x_0,\infty)$,
$j=1,2$. Denote by $M_{j}(z,x_0,\alpha)$,  the unique
Weyl-Titchmarsh matrices \eqref{2.62} corresponding to the Dirac-type
operators $D_{j}$, $j=1,2$, in \eqref{2.61}. Then,
\begin{equation}
\text{if for some $a>0$, }\, B_1(x)=B_2(x) \, \text{ for a.e. $x\in
(x_0-a,x_0+a)$,} \lb{4.110} 
\end{equation}
one obtains 
\begin{equation}
\|M_{1}(z,x_0,\alpha)-M_{2}(z,x_0,\alpha)\|_{\bbC^{2m\times
2m}}\underset{\substack{\abs{z} \to\infty\\ z\in \rho_{+}}}{=}
O\big(e^{-2\Im(z)a}\big) \lb{4.111}
\end{equation}
along any ray $\rho_+\subset\bbC_+$ with $0<\arg(z)<\pi$ and for all
$\alpha\in\bbC^{m\times 2m}$ satisfying \eqref{2.8e}. Conversely,
fix a $\hat\alpha\in\bbC^{m\times 2m}$ satisfying
\eqref{2.8e} and if $m>1$, assume in addition that
$B_j\in L^\infty([x_0-a,x_0+a])^{2m\times 2m}$, $j=1,2$. Moreover, suppose
that for all $\varepsilon >0$, 
\begin{equation}
\|M_{1}(z,x_0,\hat\alpha_1)-M_{2}(z,x_0,\hat\alpha_1)\|_{\bbC^{2m\times
2m}}\underset{\substack{|z|\to\infty\\z\in \rho_{+,\ell}}}{=}
O\big(e^{-2\Im(z)(a-\varepsilon)}\big), \quad \ell=1,2, \lb{4.112}
\end{equation}
along a ray $\rho_{+,1}\subset\bbC_+$ with $0<\arg(z)<\pi/2$ and along a
ray $\rho_{+,2}\subset\bbC_+$ with $\pi/2<\arg(z)<\pi$. Then 
\begin{equation}
B_1(x)=B_2(x) \text{ for a.e. } x\in [x_0-a,x_0+a].  \lb{4.113}
\end{equation}
\end{theorem}
\begin{proof}
\eqref{4.111} is proved by combining \eqref{2.620},  and
\eqref{4.38AA}, \eqref{4.38}, and \eqref{4.113} then follows by  combining
\eqref{2.620},  and \eqref{4.39}, \eqref{4.40}, taking into
account  the asymptotic expansions 
\begin{equation}
M_\pm (z,x_0)\underset{|z|\to\infty}{=}\pm iI_m + o(1) \lb{4.114}
\end{equation} 
along any ray with $\varepsilon<\arg(z)<\pi-\varepsilon$ in the case
of Dirac-type  operators (cf.~\eqref{3.1}).
\end{proof}

\begin{remark} \lb{r4.12}
Theorem~\ref{t4.10} and Theorem~\ref{t4.11} yield new global 
uniqueness theorems for half-line and full-line Dirac-type operators, 
extending the classical Borg-Marchenko-type results. Indeed,
if \eqref{4.39} (resp., \eqref{4.112}) holds for all $a>0$, then
\eqref{4.40} (resp. \eqref{4.113}) holds for a.e.~$x\in[x_0,\infty)$
(resp., for a.e.~$x\in\bbR$).
\end{remark}

In the case of scalar Schr\"odinger operators, the analog of 
Theorem~\ref{t4.10} is due to Simon \cite{Si98}. An alternative proof,
applicable to matrix-valued Schr\"odinger operators was presented in
\cite{GS99} (cf.~also \cite{GKM00}). More recently, yet another proof was
found by Bennewitz \cite{Be00} (following  some ideas in \cite{Bo52}).
These results extend the classical (global) uniqueness results due to
Borg \cite{Bo52} and Marchenko \cite{Ma50}, \cite{Ma52} (cf. also
\cite{Be00a}), which state that half-line $m$-functions uniquely
determine the corresponding potential coefficient. The Dirac-type results
presented in this section (especially, all local considerations) appear
to be new, even in the special case $m=1$. Previous results in the Dirac
case focused on global uniqueness questions only. We refer to Gasymov and
Levitan \cite{GL66} in the case $m=1$ and to Lesch and Malamud
\cite{LM00} in the matrix case $m\in\bbN$. Most recently, Alexander 
Sakhnovich kindly informed us that his integral representation of the
Weyl-Titchmarsh matrix in \cite{Sa88a} can be used to derive asymptotic
expansions for the Weyl-Titchmarsh matrix and its associated
matrix-valued spectral function, and also yields a result analogous to
Theorem~\ref{t4.10}\,(i) for a certain class of canonical systems.
Moreover, in the case of skew-adjoint Dirac-type systems, similar results
are discussed in \cite{Sa90} and applied to the nonlinear Schr\"odinger
equation on a half-axis.

Although not directly used in this paper, it should be pointed out that 
inverse monodromy problems for canonical systems received a lot of
attention (some of it very recently). The interested reader is referred
to \cite{AD97}, \cite{AD00}, \cite{AD00a}, \cite{Ma95}, \cite{Ma99},
\cite{Ma99a}, \cite{Sa94}, \cite{Sa99a} and the extensive literature
cited therein. Moreover, inverse spectral theory associated with
canonical systems is discussed in \cite{MST01}, \cite{Sa90}, \cite{Sa97b},
\cite{Sa01}, \cite{Sa94}, \cite{Sa94a}, \cite{Sa99}, \cite{Sa99a} (see
also the extensive literature cited in \cite{GKM00}).

\section{Trace Formulas and Borg-Type Theorems}\lb{s6}

In our final section we derive a trace formula for $B$ and then discuss
its application to Borg-type uniqueness theorems for Dirac-type 
operators.

\begin{theorem} \lb{t5.1}
Assume Hypothesis~\ref{h2.1} with $A=I_{2m}$, and let 
$\alpha_0=(I_m\; 0)\in\bbC^{m\times 2m}$. Fix $x_0\in\bbR$ and suppose 
that for all $R>0$,
\begin{align}
& \quad \underset{y\in [x_0,x_0+R]}{\esssup} \,
\bigg\|\int_y^{x_0+R} dx'\,B(x')\exp(2iz(x'-y))
+\f{1}{2iz}B(y)\bigg\|_{\bbC^{2m\times 2m}} \no \\
& + \underset{y\in [x_0-R,x_0]}{\esssup} \,
\bigg\|\int_{x_0-R}^{y} dx'\,B(x')\exp(2iz(x'-y))
-\f{1}{2iz}B(y)\bigg\|_{\bbC^{2m\times 2m}} \no \\
& \underset{\substack{\abs{z}
\to\infty\\ z\in \rho_+}}{=}\oh\big(|z|^{-1}\big) \lb{5.0}
\end{align}
along a ray $\rho_+\subset\bbC_+$. In addition, assume
$B_{k,k'}B_{\ell,\ell'}\in 
L^1_{\loc}(\bbR)^{m\times m}$ for all
$k,k',\ell,\ell'\in\{1,2\}$. Then, with $\Upsilon(\lambda,x,\alpha_0)$ 
defined in \eqref{2.69},
\begin{align}
& \begin{pmatrix}B_{1,1}(x)-B_{2,2}(x) &
B_{1,2}(x)+B_{2,1}(x) \\
 B_{1,2}(x)+B_{2,1}(x)& B_{2,2}(x)-B_{1,1}(x)\end{pmatrix} \no \\
&=\lim_{\substack{|z|\to \infty\\z\in\rho_+}} 2 \int_\bbR
d\lambda \, z^2(\lambda-z)^{-2}\Upsilon(\lambda,x,\alpha_0) 
\, \text{ for a.e. $x\in\bbR$}. \lb{5.1}
\end{align}
\end{theorem}
\begin{proof}
By \eqref{2.64},
\begin{equation}
\f{d}{dz}\ln(M(z,x,\alpha_0)) =\int_\bbR d\lambda\,
(\lambda-z)^{-2}\Upsilon(\lambda,x,\alpha_0). \lb{5.2}
\end{equation}
Next, suppose that $x\in\bbR$ is a left and right Lebesgue point of $B$.
By \eqref{4.103}, \eqref{4.104} one obtains
\begin{align}
&\f{d}{dz}\ln(M(z,x,\alpha_0)) \no \\
&\underset{\substack{|z|\to \infty\\z\in\rho_+}}{=}\f{1}{4}
\begin{pmatrix}B_{1,1}(x+0)-B_{2,2}(x+0) & B_{1,2}(x+0)+B_{2,1}(x+0) \\
B_{1,2}(x+0)+B_{2,1}(x+0)&
B_{2,2}(x+0)-B_{1,1}(x+0)\end{pmatrix}z^{-2} \lb{5.3} \\
& \quad \;\,\, +\f{1}{4}
\begin{pmatrix}B_{1,1}(x-0)-B_{2,2}(x-0) & B_{1,2}(x-0)+B_{2,1}(x-0) \\
B_{1,2}(x-0)+B_{2,1}(x-0)&
B_{2,2}(x-0)-B_{1,1}(x-0)\end{pmatrix}z^{-2} +\oh\big(z^{-2}\big) \no
\end{align}
and hence 
\begin{align}
& \quad \, \f{1}{2}\begin{pmatrix}B_{1,1}(x+0)-B_{2,2}(x+0) &
B_{1,2}(x+0)+B_{2,1}(x+0) \\
 B_{1,2}(x+0)+B_{2,1}(x+0)& B_{2,2}(x+0)-B_{1,1}(x+0)\end{pmatrix} \no \\
&+\f{1}{2}\begin{pmatrix}B_{1,1}(x-0)-B_{2,2}(x-0) &
B_{1,2}(x-0)+B_{2,1}(x-0) \\
 B_{1,2}(x-0)+B_{2,1}(x-0)& B_{2,2}(x-0)-B_{1,1}(x-0)\end{pmatrix} \no \\
&=\lim_{\substack{|z|\to \infty\\z\in\rho_+}} 2 \int_\bbR
d\lambda \, z^2(\lambda-z)^{-2}\Upsilon(\lambda,x,\alpha_0). \lb{5.2a}
\end{align}
Since a.e. $x\in\bbR$ is a Lebesgue point of $B$, one concludes
\eqref{5.1}.
\end{proof}

In the case $m=1$, a trace formula for Dirac-type operators, using Krein
spectral shift functions and exponential representations of Herglotz
functions, was discussed in \cite{Ti95}. This circle of ideas was first
introduced in connection with trace formulas of Schr\"odinger operators
in \cite{GS96} (see also \cite{GHS95}, \cite{GHSZ95}, \cite{Ry99},
\cite{Ry99a} in the scalar case $m=1$. The  corresponding case of trace
formulas for matrix-valued Schr\"odinger operators was introduced in
\cite{GH97} (see also \cite{CGHL00}). 

Analogous trace formulas can be drived for all higher-order coefficients
$M_k(x,\alpha_0)$ in \eqref{4.103} (see, e.g., \cite{GHSZ95} in connection
with scalar Schr\"odinger operators). 

A comparison of the trace formula (3.20) in
\cite{CGHL00} for Schr\"{o}dinger operators with its Dirac-type
counterpart \eqref{5.1} reveals characteristic differences. While in the
Schr\"{o}dinger case the trace formula directly involves the potential
coefficient $Q(x)$, $M_1(x,\alpha_0)$ differs markedly from a constant
multiple of $B(x)$, and consequently, the Dirac-type trace formula
\eqref{5.1} does not directly involve $B(x)$ but certain linear
combinations of $B_{j,k}(x)$. This is related to the fact that
$M(z,x_0,\alpha_0)$ (or equivalently, $\Upsilon(\lambda,x_0,\alpha_0)$),
in general, does not uniquely determine $B$ a.e. In fact, there exists a
typical ambiguity concerning the coefficients of $D$ related to unitary
gauge-transformations of $D$. In the case $m=1$ this ambiguity is
well-known and discussed, e.g., in \cite{GL66}, \cite[Sect.~I.10]{LS75},
\cite[Ch.~7]{LS91}. These gauge transformations leave the spectrum of
$D$ invariant and suggest that we focus our attention on certain normal
forms of $D$ in connection with inverse spectral problems for Dirac-type
operators.

\begin{lemma} \lb{l5.2}
Assume Hypothesis~\ref{h2.4}. Then $D=J\f{d}{dx}-B$ is unitarily
equivalent to $\wti D$, where $\wti D$ in $L^2(\bbR)^{2m}$ is of the
normal form
\begin{equation}
\wti D=J\f{d}{dx}-\wti B= \begin{pmatrix} -\wti B_{1,1} & -I_m
\f{d}{dx}- \wti B_{1,2} \\[1mm]
 I_m \f{d}{dx}-\wti B_{1,2}&  \wti B_{1,1} \end{pmatrix}. \lb{5.4}
\end{equation}
Here $\wti B=\wti B^*$ a.e. and 
\begin{align}
\wti B_{1,1}&=-(1/2)\Im\big(U_{1,1}^{-1}[(B_{1,2}+B_{2,1})
-i(B_{1,1}-B_{2,2})] U_{2,2}\big)= \wti B_{1,1}^*, \lb{5.5} \\
\wti B_{1,2}&=(1/2)\Re\big(U_{1,1}^{-1}[(B_{1,2}+B_{2,1})
-i(B_{1,1}-B_{2,2})] U_{2,2}\big)=\wti B_{1,2}^*, \lb{5.6}
\end{align}
with $U_{j,j}\in\bbC^{m\times m}$, $j=1,2$, satisfying the
first-order system of ordinary differential equations
\begin{align}
iU_{j,j}^\prime(x)&=-(1/2)\big((-1)^j(B_{1,1}(x)
+B_{2,2}(x))+i(B_{1,2}(x)-B_{2,1}(x)) \big)U_{j,j}(x), \no \\ 
&\hspace*{6cm} \text{for a.e. $x\in\bbR$}, \quad j=1,2. \lb{5.7}
\end{align}
\end{lemma}
\begin{proof}
We start with the unitary transformation $V$ in $L^2(\bbR)^{2m}$ defined
by
\begin{equation}
V=\f1{\sqrt{2}}\begin{pmatrix} i I_m & I_m \\ I_m & i I_m \end{pmatrix}, 
\quad V^{-1}=\f1{\sqrt{2}}\begin{pmatrix} -i I_m & I_m \\ I_m & -i I_m
\end{pmatrix}, \lb{5.8}
\end{equation}
which maps $D$ to $D_1$, where
\begin{align}
D_1&=V^{-1}DV= i \begin{pmatrix} I_m\f{d}{dx} & 0 \\
0 & -I_m\f{d}{dx} \end{pmatrix} \no \\
& \quad -\f12 \begin{pmatrix}B_{1,1}+B_{2,2}-i(B_{1,2}-B_{2,1}) &
B_{1,2}+B_{2,1}-i(B_{1,1}-B_{2,2}) \\ B_{1,2}+B_{2,1}+i(B_{1,1}-B_{2,2}) &
B_{1,1}+B_{2,2}+i(B_{1,2}-B_{2,1})\end{pmatrix}. \lb{5.9}
\end{align}
Next, we introduce the unitary operator $U$ in $L^2(\bbR)^{2m}$  defined
by
\begin{equation}
U=\begin{pmatrix}U_{1,1} & 0 \\ 0 & U_{2,2}\end{pmatrix}, \lb{5.10}
\end{equation}
where the unitary $m\times m$ matrices $U_{j,j}\in\bbC^{m\times m}$ are
solutions of the first-order system \eqref{5.7}. Since by hypothesis
$B_{j,k}\in L^1_{\loc}(\bbR)^{m\times m}$, $1\leq j,k \leq 2$, the
solutions of equation \eqref{5.7} are well-defined and
$U_{j,j}\in\AC_{\loc}(\bbR)^{m\times m}$, $j=1,2$. One computes
\begin{equation}
\what D=U^{-1} D_1 U= \begin{pmatrix}i I_m \f{d}{dx} & -\what B_{1,2} 
\\[1mm] -\what B_{1,2}^* & -i I_m \f{d}{dx}\end{pmatrix}, \lb{5.11}
\end{equation}
where $\what B_{1,2}\in L^1_{\loc}(\bbR)^{m\times m}$ and
\begin{equation}
\what B_{1,2}(x)=(1/2) U^{-1}_{1,1}(x)\big(B_{1,2}(x)+B_{2,1}(x)
-i(B_{1,1}(x)-B_{2,2}(x)\big)U_{2,2}(x). \lb{5.12}
\end{equation}
Finally, defining $\wti D=V\what D V^{-1}$, one arrives at 
\eqref{5.4}--\eqref{5.6}.
\end{proof}

Thus, unitary invariants of $D$ (such as the spectrum, $\spec (D)$, of $D$
and its multiplicity) cannot determine $B$ in general but at best a
potential matrix of the type (normal form) $\wti B$ in \eqref{5.4}.  A
further restriction on the solvability of inverse spectral
problems for Dirac-type operators is mentioned in the following result.

\begin{lemma} \lb{l5.3}
Assume Hypotheses~\ref{h2.4} and let 
$\omega=\omega^*\in\bbC^{m\times m}$ be a constant self-adjoint $m\times
m$ matrix. Then $D=J\f{d}{dx}-B$ is unitarily equivalent to $\wti
D_\omega$ in $L^2(\bbR)^{2m}$, where
\begin{equation}
\wti D_\omega=J\f{d}{dx}-\wti B_\omega=\begin{pmatrix} -\wti
B_{\omega,1,1} & -I_m\f{d}{dx}-\wti B_{\omega,1,2} \\[1mm] 
I_m\f{d}{dx}-\wti B_{\omega,1,2} & \wti B_{\omega,1,1}\end{pmatrix},
\lb{5.13}
\end{equation}
with 
\begin{align} 
\wti B_{\omega,1,1}&=-(1/2)\Im\big(e^{i\omega}U_{1,1}^{-1}[(B_{1,2}
+B_{2,1})-i(B_{1,1}-B_{2,2})]
U_{2,2}e^{i\omega}\big)=\wti B_{\omega,1,1}^*,  \no \\
\wti B_{\omega,1,2}&=(1/2)\Re\big(e^{i\omega}U_{1,1}^{-1}[(B_{1,2}
+B_{2,1})-i(B_{1,1}-B_{2,2})] U_{2,2}e^{i\omega}\big)=\wti
B_{\omega,1,2}^*, \lb{5.14}
\end{align}
and with $U_{j,j}$, $j=1,2$, satisfying the first-order system
\eqref{5.7}. 
\end{lemma}
\begin{proof}
Define 
\begin{equation}
U_\omega=\begin{pmatrix}e^{i\omega}& 0 \\ 
0 & e^{-i\omega} \end{pmatrix}. \lb{5.15}
\end{equation}
Using the notation employed in the proof of Lemma~\ref{l5.2} one verifies
that
\begin{equation}
\wti D_\omega=V U_\omega (VU)^{-1}DVU(VU_\omega)^{-1}. \lb{5.16}
\end{equation}
\end{proof}

In particular, choosing $\omega=(\pi/2)I_m$ effects the sign change 
$\wti B\to -\wti B$, with $\wti B$ given by \eqref{5.5}, \eqref{5.6}.

For detailed discussions of various normal forms for Dirac-type operators
we refer to \cite{GL66}, \cite{HJKS91}, \cite[Ch.~9]{LS75},
\cite[Ch.~7]{LS91} in the case $m=1$ and to \cite{Ga68}, \cite{LM00},
\cite{Ma99}, \cite[p.~193--195]{Ma86}, \cite{Ma65} in the general
matrix-valued case. Perhaps it should be noted that if $D$ is in its
normal form $\wti D$ as in \eqref{5.4}, ${\wti D}^2$ turns into a
$2m\times 2m$ matrix-valued Schr\"odinger operator under appropriate
regularity assumptions on $\wti B$. Details on this fact and the relation
between the $M$-matrices of $\wti D$ and ${\wti D}^2$ can be found in
Section~3 of \cite{GKM00}.

Next, we turn to Borg-type theorems, one of the principal topics of this
paper. In 1946 Borg \cite{Bo46} proved, among a variety of other inverse
spectral theorems, the following result.

\begin{theorem}[\cite{Bo46}] \lb{t5.4}
Assume $q\in L^2_{\loc} (\bbR)$ to be real-valued and periodic and let
\begin{equation}
h=-d^2/dx^2+q \lb{5.16a}
\end{equation}
be the associated self-adjoint Schr\"odinger operator in
$L^2(\bbR)$. Moreover, suppose that $\spec (h)=[e_0,\infty)$ for some
$e_0\in\bbR$. Then
\begin{equation}
q(x)=e_0 \, \text{ for a.e. $x\in\bbR$}. \lb{5.17}
\end{equation}
\end{theorem}

The analog of Theorem~\ref{t5.4} for Dirac-type operators (in the case
$m=1$) was proven by Giacheti and Johnson \cite{GJ84} in 1984 (see also
\cite{Ge89}, \cite{Ge91}, \cite{GSS91} in the special case where $p$ is
constant and \cite{GG93} in the case where $p,q\in L^2 (\bbR)$ are
real-valued and periodic).

\begin{theorem}[\cite{GJ84}] \lb{t5.4a}
Assume $p,q\in L^\infty (\bbR)$ to be real-valued and periodic and let
\begin{equation}
d=\begin{pmatrix}-p &-\f{d}{dx}-q\\ \f{d}{dx}-q & p\end{pmatrix}
\lb{5.17a}
\end{equation}
be the associated self-adjoint Dirac-type operator in
$L^2(\bbR)^2$. Moreover, suppose that $\spec (d)=\bbR$. Then
\begin{equation}
p(x)=q(x)=0 \, \text{ for a.e. $x\in\bbR$}. \lb{5.17b}
\end{equation}
\end{theorem}

Traditionally, uniqueness results such as Theorems~\ref{t5.4} and
\ref{t5.4a} are called Borg-type theorems. (However, this terminology is
not uniquely adopted and hence a bit unfortunate. Indeed, inverse spectral
results on finite intervals recovering the potential coefficient(s) from
several spectra, were also pioneered by Borg in his celebrated paper
\cite{Bo46}, and hence are also coined Borg-type theorems in the
literature, see, e.g., \cite{Ma94}, \cite{Ma99}, \cite{Ma99a}.) 

A quick and natural proof of Theorem~\ref{t5.4}, based on a trace formula
for $q$, was presented in \cite{CGHL00}. This strategy of proof was then
applied to the case of matrix-valued Schr\"odinger operators and the
corresponding matrix-valued analog of Theorem~\ref{t5.4} was also proved 
in \cite{CGHL00} along these lines. A closer examination of the proof of
Theorem~\ref{t5.4} shows that periodicity of $q$ is not the crucial element
in the proof of the 
uniqueness result \eqref{5.17}. The key ingredient (besides
$\spec (h)=[e_0,\infty)$) is clearly
the fact that for all $x\in\bbR$,
\begin{equation}
\xi(\lambda,x)=1/2 \, \text{ for a.e. } \lambda\in\spec_{\ess}(h)
\lb{5.18}
\end{equation}
($\spec_{\ess}(\dott)$ the essential spectrum), where $\xi(\cdot,x)$ is
defined by
\begin{equation}
\xi(\lambda,x)=\lim_{\varepsilon\to
0}\pi^{-1}\Im(\ln(g(\lambda+i\varepsilon,x))) \, \text{ for a.e.
$\lambda\in\bbR$}, \lb{5.19}
\end{equation}
and $g(z,x)$ denotes Green's function (i.e., the integral kernel of the
resolvent) of $h$ on the diagonal,
\begin{equation}
g(z,x)=(h-z)^{-1}(x,x). \lb{5.20}
\end{equation}
Completely analogous considerations apply to the Dirac-type case.

Real-valued periodic potentials are known to satisfy \eqref{5.18} but so
are certain classes of real-valued quasi-periodic and almost-periodic
potentials $q$ (see, e.g., \cite{CJ87}, \cite{Cr89}, \cite{DS83},
\cite{Jo82}, \cite{JM82}, \cite{Ko84}, \cite{Ko87a}, 
\cite{KK88}, \cite{KS88}, \cite{SY95}). In particular, the class of
real-valued algebro-geometric finite-gap  potentials $q$ (a subclass of
the set of real-valued quasi-periodic potentials) is a prime example
satisfying
\eqref{5.18} without necessarily being periodic. Traditionally,
potentials $q$ satisfying \eqref{5.18} are called \textit{reflectionless}
(see \cite{Cr89}, \cite{DS83}, \cite{KK88}, \cite{SY95}).
Again the analogous notions apply to the Dirac-type case (cf., e.g.,
\cite{CJ87}, \cite{GJ84}, \cite{Jo87}).

Taking this circle of ideas as the point of departure for our derivation
of Borg-type results for Dirac-type operators, we now use the
reflectionless situation described in \eqref{5.18}, actually, its proper
analog for Dirac-type systems, as the model for the subsequent definition.

\begin{definition}\lb{d5.5}
Assume Hypothesis~\ref{h2.1} with $A=I_{2m}$, and let 
$\alpha_0=(I_m\; 0)\in\bbC^{m\times 2m}$. Then $B$ is called
{\it reflectionless} if for all $x\in\bbR$,
\begin{equation}
\Upsilon(\lambda,x,\alpha_0)= (1/2) I_{2m} \, \text{  for a.e.\
$\lambda\in\spec_{\ess}(D)$}. \lb{5.21}
\end{equation}
\end{definition}

Since hardly any confusion can arise, we will also call the Dirac-type
operator $D$ reflectionless if \eqref{5.21} is satisfied.  

Given Definition~\ref{d5.5}, we turn to a Borg-type uniqueness theorem
and formulate the analog of Theorem~\ref{t5.4} for (reflectionless)
Dirac-type operators. 

\begin{theorem}\lb{t5.6}
Assume Hypothesis \ref{h2.1} with $A=I_{2m}$, and let 
$\alpha_0=(I_m\; 0)\in\bbC^{m\times 2m}$. If for all $x\in\bbR$,
$\Upsilon(\lambda,x,\alpha_0)=C$ is a constant $2m\times 2m$ matrix for
a.e.\ $\lambda\in\bbR$,  especially, if $B$ is reflectionless and
$\spec(D)=\bbR$, then
\begin{equation}
B_{1,1}(x)=B_{2,2}(x), \quad B_{1,2}(x)=-B_{2,1}(x) \, \text{ for a.e.
$x\in\bbR$}. \lb{5.22}
\end{equation}
In particular, if $D$ is assumed to be in its normal form \eqref{5.4},
that is, of the type $\wti D = J\f{d}{dx} - \wti B$, then
\begin{equation}
\wti B (x)=0 \text{ for a.e. $x\in\bbR$}. \lb{5.23}
\end{equation}
\end{theorem}
\begin{proof}
The fact that $\int_\bbR d\lambda\, (\lambda-z)^{-2}=0$ for all
$z\in\bbC\backslash\bbR$, that a.e.~$x\in\bbR$ is a Lebesgue point of
$B$, and the trace formula \eqref{5.1}, imply \eqref{5.22}. Together with
Lemma \ref{l5.2} this yields \eqref{5.23}.
\end{proof}

The analog of Theorem~\ref{t5.6} for matrix-valued Schr\"odinger
operators was recently proved in \cite{CGHL00}.

In the remainder of the section we will show that the 
case of periodic $B$ is covered by Theorem~\ref{t5.6} under appropriate 
uniform multiplicity assumptions on $\spec(D)$. In order to handle Floquet
theoretic aspects of periodic Dirac-type operators $D$, we adopt the 
following assumptions until the end of this section.

\begin{hypothesis} \lb{h5.7}
In addition to Hypothesis~\ref{h2.1} assume $A=I_{2m}$ and suppose that
$B$ is periodic, that is, there is an $\omega>0$ such that
$B(x+\omega)=B(x)$ for a.e.~$x\in\bbR$.
\end{hypothesis}

The following result has been proven in \cite[Theorem~4.6]{CGHL00}.

\begin{theorem} [\cite{CGHL00}, Theorem~4.6] \lb{t5.8}
Assume Hypothesis~\ref{h5.7} and let $\alpha_0=(I_m\; 0)\in\bbC^{2m\times
m}$. If $D$ has uniform spectral multiplicity $2m$, then for all
$x\in\bbR$  and all $\lambda\in\spec(D)^o$,
\begin{equation}
M_+(\lambda+i0,x,\alpha_0)=M_-(\lambda+i0,x,\alpha_0)^*
=M_-(\lambda-i0,x,\alpha_0).
\lb{5.24}
\end{equation}
In particular, $M_-(z,x,\alpha_0)$ is the analytic continuation of 
$M_+(z,x,\alpha_0)$ {\rm (}and vice versa{\rm )} through $\spec(D)^o$.
\end{theorem}
\noindent Here $A^o$ denotes the open interior of a set $A\subseteq\bbR$.

Strictly speaking, Theorem~4.6 in \cite{CGHL00} was proved for
matrix-valued Schr\"odinger operators. But the proof extends line by line
to the corresponding Dirac-type situation and was predominantly
formulated in terms of Hamiltonian systems notation (rather than
Schr\"odinger operator specifics) in order to be applicable to the
present context. In particular, the spectrum, $\spec(H)$, of the
Schr\"odinger operator $H$ should be replaced by that of $D$, the point
spectrum, $\spec_p(H^D_{x_0})$, of the Dirichlet Schr\"odinger operator
$H^D_{x_0}$ with a Dirichlet boundary condition at the point $x_0$ should
simply be replaced by the set $\{\lambda\in\bbR\,|\,
\det(\phi_1(\lambda,x_0+\omega,x_0,\alpha_0))=0\}$, etc.

\begin{theorem} \lb{t5.9}
Suppose Hypothesis~\ref{h5.7} and let $\alpha_0=(I_m\; 0)
\in\bbC^{2m\times m}$. If $D$ has uniform spectral multiplicity
$2m$, then $D$ is reflectionless and for all $x\in\bbR$ and all
$\lambda\in\spec(D)^o$,
\begin{equation}
\Upsilon(\lambda,x,\alpha_0)= (1/2) I_{2m}. \lb{5.30}
\end{equation}
\end{theorem}
\begin{proof}
This is clear from \eqref{2.620} and \eqref{5.24}, which imply
\begin{equation}
M(\lambda+i0,x,\alpha_0)=-M(\lambda+i0,x,\alpha_0)^*. \lb{5.27}
\end{equation}
\end{proof}

Theorems~\ref{t5.8} and \ref{t5.9} extend to more general
situations  (not necessarily periodic ones) as is clear from the 
corresponding results in \cite{CJ87}, \cite{GJ84}, \cite{GKT96},
\cite{Ko84}, \cite{Ko87a}, \cite{KK88}, \cite{SY95} in 
the scalar case $m=1$ (replacing the phrase ``for all
$\lambda\in\spec(D)^o$'' by  ``for~a.e.~$\lambda\in\spec(D)^o$'', etc.).
For the corresponding  matrix-valued Schr\"odinger operator case we refer
to \cite{KS88}. 

\begin{corollary} \lb{c5.10}
Assume Hypothesis~\ref{h5.7}. If $D$ has uniform spectral multiplicity
$2m$ and $\spec(D)=\bbR$, then
\begin{equation}
B_{1,1}(x)=B_{2,2}(x), \quad B_{1,2}(x)=-B_{2,1}(x) \, \text{ for a.e.
$x\in\bbR$}. \lb{5.31}
\end{equation}
In particular, if $D$ is assumed to be in its normal form $\wti D =
J\f{d}{dx} - \wti B$, with $\wti B$ given by \eqref{5.4}, then
\begin{equation}
\wti B (x)=0 \, \text{ for a.e. $x\in\bbR$}. \lb{5.32}
\end{equation}
\end{corollary}

\begin{remark} \lb{r5.11}
The assumption of uniform (maximal) spectral multiplicity $2m$ in
Corollary~\ref{c5.10} is an essential one. Otherwise, one can easily 
construct nonconstant potentials $B$ such that the associated operator
$D$ has overlapping  band spectra and hence spectrum the whole real line. 
Also self-adjointness of $B$ is crucial for Corollary~\ref{c5.10} to hold
(cf.~the corresponding discussion in Remark~4.2 of \cite{CGHL00} in the
context of Schr\"odinger operators). 
\end{remark}

The analog of Corollary~\ref{c5.10} for periodic matrix-valued
Schr\"odinger operators was first proved by Depres \cite{De95} and
recently rederived using such a trace formula approach in \cite{CGHL00}.

We note that all results presented in this paper also apply to
matrix-valued finite-difference Hamiltonian systems.  We refer the 
reader to \cite{CGR01} in this direction.

Finally, Borg-type uniqueness theorems for Hamiltonian systems are just a
beginning. There is  a natural extension of Borg's Theorem~\ref{t5.4} to 
self-adjoint periodic Schr\"{o}dinger, respectively, Dirac-type 
operators with one gap, respectively, two gaps in their spectrum. In the
case of (scalar) Schr\"odinger operators, such an extension is due to
Hochstadt \cite{Ho65}  and the resulting potential $q$ becomes twice the
elliptic Weierstrass function. In the case of Dirac-type operators (with
$m=1$ and vanishing diagonal coefficients in $B$) such an extension
involving elliptic functions can be found in \cite{Ge89}, \cite{Ge91},
\cite{GSS91} (see also \cite{GW98}). Extensions to matrix-valued versions
(i.e., for 
$m\geq 2$) are currently under active investigations.

\vspace*{3mm}
\noindent {\bf Acknowledgements.}
We would like to thank Suzanne Collier, Helge Holden, Konstantin Makarov,
Fedor Rofe-Beketov,  Alexei Rybkin, Lev Sakhnovich, and Barry Simon for
helpful discussions and  many hints regarding the literature, and
especially, Don Hinton,  Boris Levitan, Mark Malamud, and Alexander
Sakhnovich for repeated correspondence on various parts of  the material
in this paper. \\
S.~C. would like to thank the Mathematics Department of the
University of Missouri-Columbia for the great hospitality extended to
him during his 2000/2001 sabbatical when this work was completed.


\end{document}